\documentclass[12pt,twoside]{article}
\usepackage{amsfonts,amssymb,amsmath,bm,color,hyperref}
\usepackage{a4wide}
\definecolor{GreenYellow}{cmyk}{0.15,0,0.69,0}
\definecolor{Yellow}{cmyk}{0,0,1,0}
\definecolor{PaleYellow}{cmyk}{0,0,0.5,0}
\definecolor{Goldenrod}{cmyk}{0,0.10,0.84,0}
\definecolor{Dandelion}{cmyk}{0,0.29,0.84,0}
\definecolor{Apricot}{cmyk}{0,0.32,0.52,0}
\definecolor{Peach}{cmyk}{0,0.50,0.70,0}
\definecolor{Melon}{cmyk}{0,0.46,0.50,0}
\definecolor{YellowOrange}{cmyk}{0,0.42,1,0}
\definecolor{Orange}{cmyk}{0,0.61,0.87,0}
\definecolor{BurntOrange}{cmyk}{0,0.51,1,0}
\definecolor{Bittersweet}{cmyk}{0,0.75,1,0.24}
\definecolor{RedOrange}{cmyk}{0,0.77,0.87,0}
\definecolor{Mahogany}{cmyk}{0,0.85,0.87,0.35}
\definecolor{Maroon}{cmyk}{0,0.87,0.68,0.32}
\definecolor{BrickRed}{cmyk}{0,0.89,0.94,0.28}
\definecolor{Red}{cmyk}{0,1,1,0}
\definecolor{OrangeRed}{cmyk}{0,1,0.50,0}
\definecolor{RubineRed}{cmyk}{0,1,0.13,0}
\definecolor{WildStrawberry}{cmyk}{0,0.96,0.39,0}
\definecolor{Salmon}{cmyk}{0,0.53,0.38,0}
\definecolor{CarnationPink}{cmyk}{0,0.63,0,0}
\definecolor{Magenta}{cmyk}{0,1,0,0}
\definecolor{VioletRed}{cmyk}{0,0.81,0,0}
\definecolor{Rhodamine}{cmyk}{0,0.82,0,0}
\definecolor{Mulberry}{cmyk}{0.34,0.90,0,0.02}
\definecolor{RedViolet}{cmyk}{0.07,0.90,0,0.34}
\definecolor{Fuchsia}{cmyk}{0.47,0.91,0,0.08}
\definecolor{Lavender}{cmyk}{0,0.48,0,0}
\definecolor{Thistle}{cmyk}{0.12,0.59,0,0}
\definecolor{Orchid}{cmyk}{0.32,0.64,0,0}
\definecolor{DarkOrchid}{cmyk}{0.40,0.80,0.20,0}
\definecolor{Purple}{cmyk}{0.45,0.86,0,0}
\definecolor{Plum}{cmyk}{0.50,1,0,0}
\definecolor{Violet}{cmyk}{0.79,0.88,0,0}
\definecolor{RoyalPurple}{cmyk}{0.75,0.90,0,0}
\definecolor{BlueViolet}{cmyk}{0.86,0.91,0,0.04}
\definecolor{Periwinkle}{cmyk}{0.57,0.55,0,0}
\definecolor{CadetBlue}{cmyk}{0.62,0.57,0.23,0}
\definecolor{CornflowerBlue}{cmyk}{0.65,0.13,0,0}
\definecolor{MidnightBlue}{cmyk}{0.98,0.13,0,0.43}
\definecolor{NavyBlue}{cmyk}{0.94,0.54,0,0}
\definecolor{RoyalBlue}{cmyk}{1,0.50,0,0}
\definecolor{Blue}{cmyk}{1,1,0,0}
\definecolor{Cerulean}{cmyk}{0.94,0.11,0,0}
\definecolor{Cyan}{cmyk}{1,0,0,0}
\definecolor{ProcessBlue}{cmyk}{0.96,0,0,0}
\definecolor{SkyBlue}{cmyk}{0.62,0,0.12,0}
\definecolor{Turquoise}{cmyk}{0.85,0,0.20,0}
\definecolor{TealBlue}{cmyk}{0.86,0,0.34,0.02}
\definecolor{Aquamarine}{cmyk}{0.82,0,0.30,0}
\definecolor{BlueGreen}{cmyk}{0.85,0,0.33,0}
\definecolor{Emerald}{cmyk}{1,0,0.50,0}
\definecolor{JungleGreen}{cmyk}{0.99,0,0.52,0}
\definecolor{SeaGreen}{cmyk}{0.69,0,0.50,0}
\definecolor{Green}{cmyk}{1,0,1,0}
\definecolor{ForestGreen}{cmyk}{0.91,0,0.88,0.12}
\definecolor{PineGreen}{cmyk}{0.92,0,0.59,0.25}
\definecolor{LimeGreen}{cmyk}{0.50,0,1,0}
\definecolor{YellowGreen}{cmyk}{0.44,0,0.74,0}
\definecolor{SpringGreen}{cmyk}{0.26,0,0.76,0}
\definecolor{OliveGreen}{cmyk}{0.64,0,0.95,0.40}
\definecolor{RawSienna}{cmyk}{0,0.72,1,0.45}
\definecolor{Sepia}{cmyk}{0,0.83,1,0.70}
\definecolor{Brown}{cmyk}{0,0.81,1,0.60}
\definecolor{Tan}{cmyk}{0.14,0.42,0.56,0}
\definecolor{Gray}{cmyk}{0,0,0,0.50}
\definecolor{LightGray}{cmyk}{0,0,0,0.05}
\definecolor{Black}{cmyk}{0,0,0,1}
\definecolor{White}{cmyk}{0,0,0,0}
\numberwithin{equation}{section}
\newcommand{\hess}{\operatorname{Hess}}
\newcommand{\bmtheta}{\overline{\theta}}

\newcommand{\gx}{\mv g_{\vv x}}
\newcommand{\ux}{\mv u_{\vv x}}
\newcommand{\yx}{\mv y_{\vv x}}
\newcommand{\yxd}{\mv y_{\vv x'}}

\newcommand{\rg}{\mv r_{\mv g}}
\newcommand{\ru}{\mv r_{\mv u}}
\newcommand{\rf}{\mv r_{\mv y}}

\newcommand{\nz}{\smallsetminus\{0\}}
\newcommand{\df}[1]{\stackrel{\mathrm{def}}{#1}}
\newcommand{\GL}{\mathrm{GL}}
\newcommand{\pd}[2]{\partial_{{#2}}{#1}}
\newcommand{\ddd}{d}
\newcommand{\ccc}{m}
\newcommand{\mysubsection}[1]{\refstepcounter{subsection}\subsection*{\thesubsection~~#1}}
\newtheorem{theorem}{Theorem}[section]
\newtheorem{theoremKM}{Theorem KM}

\newtheorem{corollary}[theorem]{Corollary}
\newtheorem{lemma}[theorem]{Lemma}
\newtheorem{definition}[theorem]{Definition}
\newtheorem{remark}[theorem]{Remark}
\newtheorem{example}[theorem]{Example}

\newtheorem{problem}[theorem]{Problem}

\newtheorem{conjecture}[theorem]{Conjecture}
\def\P{\mathbb{P}}
\def\R{\mathbb{R}}

\def\Q{\mathbb{Q}}
\def\Z{\mathbb{Z}}
\def\N{\mathbb{N}}

\newcommand{\Rp}{\R^+}    
\renewcommand{\r}{\rho}

\def\cA{\mathcal{A}}

\def\cH{\mathcal{H}}
\def\cB{\mathcal{B}}
\def\cM{\mathcal{M}}

\def\cG{\mathcal{G}}
\def\cC{\mathcal{C}}

\def\cS{{\mathcal{S}}}

\def\cR{\mathcal{R}}

\def\cV{\mathcal{V}}
\newcommand{\ra}{R_{\alpha}}

\newcommand{\ba}{\beta_\alpha}
\newcommand{\al}{\alpha}

\newcommand{\we}{\wedge}
\newcommand{\We}{\mbox{\Large$\wedge$}}

\newcommand{\ie}{{\it i.e.}}

\newcommand{\ve}{\varepsilon}

\newcommand{\vv}[1]{{\mathbf{#1}}}
\newcommand{\mv}[1]{{\bm#1}}
\newcommand{\pdist}{d_p}

\newcommand{\Gr}{\operatorname{Gr}}

\newcommand{\dist}{\operatorname{dist}}

\newcommand{\diag}{\operatorname{diag}}

\newcommand{\proofend}{\hspace*{\fill}\raisebox{-2ex}{$\boxtimes$}}
\renewcommand{\tilde}{\widetilde}
\newcommand{\dt}{\cdot}
\newcommand{\codim}{\operatorname{codim}}
\newcommand{\vol}{\operatorname{Vol}}
\newcommand{\V}{\R^{n+1}}
\newcommand{\Jarnik}{Jarn\' \i k}
\begin{document}

\title{\bf Rational points near manifolds and\\ metric Diophantine approximation}

\author{\sc Victor Beresnevich%
 \footnote{EPSRC Advanced Research Fellow, grant no.EP/C54076X/1} ~{\small (York)}}

\date{{\small\it Dedicated to Maurice Dodson}}

\maketitle

\vspace*{-6ex}

\begin{abstract}
This work is motivated by problems on simultaneous Diophantine approximation on manifolds, namely, establishing Khintchine and \Jarnik{} type theorems for submanifolds of $\R^n$. These problems have attracted a lot of interest since Kleinbock and Margulis proved a related conjecture of Alan~Baker and V.G.~Sprind\v zuk. They have been settled for planar curves but remain open in higher dimensions. In this paper, Khintchine and \Jarnik{} type divergence theorems are established for arbitrary analytic non-degenerate manifolds regardless of their dimension. The key to establishing these results is the study of the distribution of rational points near manifolds -- a very attractive topic in its own right. Here, for the first time, we obtain sharp lower bounds for the number of rational points near non-degenerate manifolds in dimensions $n>2$ and show that they are ubiquitous (that is uniformly distributed).
\end{abstract}

{\footnotesize

}

\vspace*{0ex}

\noindent\textbf{Contents}

{\parskip=1ex

1~~~ Introduction \dotfill\pageref{rpnm}

2~~~ Diophantine approximation on manifolds
\dotfill\pageref{further_motivation}

3~~~ Some auxiliary geometry \dotfill\pageref{aux}

4~~~ Detecting rational points near a manifold
\dotfill\pageref{detecting}

5~~~ Integer points in `random' parallelepipeds
\dotfill\pageref{iprp}

6~~~ The proof of main result: Theorem~\ref {t:01}
\dotfill\pageref{proofs}

7~~~ Further theory for curves \dotfill\pageref{ext}

8~~~ Final comments \dotfill\pageref{final}

}

\newpage


\section{Introduction}\label{rpnm}

Let $\cM$ be a bounded smooth manifold in $\R^n$. Given $Q>1$ and
$\varepsilon>0$, let
$$
N(Q,\varepsilon)\ =\ \#\Big\{\vv p/q\in\Q^n:\ 1\le q\le Q,\ \dist(\vv
p/q,\cM)\le\varepsilon\Big\},
$$
where $\#S$ is the cardinality of a set $S$, $\vv p\in\Z^n$,
$q\in\Z$, $\dist(\vv r,\cM)=\inf_{\vv y\in \cM}|\vv r-\vv y|$ and
$|\cdot|$ is the Euclidean norm on $\R^n$. Thus, $N(Q,\varepsilon)$
counts rational points with bounded
denominator lying `$\varepsilon$-near' $\cM$. The following intricate problem will be our main concern.

\begin{problem}\label{p1}
Estimate $N(Q,\varepsilon)$ for a `generic' smooth manifold $\cM$.
\end{problem}

\noindent Our study of Problem~\ref{p1} is motivated by open problems on simultaneous Diophantine approximation on manifolds -- see \S\ref{further_motivation}. However, the interest to the distribution of rational points
near manifolds is not limited to these problems -- see, e.g., \cite{Elkies-2000, Mazur-04:MR2058289}. In this paper a sharp lower bound on $N(Q,\varepsilon)$ is established when $\varepsilon$ is bounded below by some naturally occurring function of $Q$. To begin with, we briefly review the state of the art.


\emph{Planar curves.} The first general estimates for $N(Q,\varepsilon)$
are due to Huxley \cite{Huxley-96:MR1420620, Huxley-1994-rational_points}. In particular, he
proved that for any curve $\cM$ in $\R^2$ with curvature
bounded between positive constants, \ $N(Q,\varepsilon)\ll \varepsilon
Q^{3+\theta}$ for $\varepsilon\gg Q^{-2}$, where $\theta>0$ is arbitrary and
``$\ll$'' is the Vinogradov symbol. Huxley's
estimate was the only general result until Vaughan and
Velani remarkably removed the $\theta$-term from Huxley's estimate \cite{Vaughan-Velani-2007}. On the other hand, Dickinson, Velani
and the author \cite{Beresnevich-Dickinson-Velani-07:MR2373145}
obtained the complementary bound $N(Q,\varepsilon)\gg
\varepsilon Q^{3}$ for $\varepsilon\gg Q^{-2}$. Consequently, the theory for planar curves is reasonably complete.

\emph{Higher dimensions.} Very little is known. Effectively, there
are only rather crude bounds on $N(Q,\varepsilon)$ obtained via
Khintchine's transference principle \cite{Bernik-03:MR2163817} and
estimates for topological products of planar curves \cite[\S4.4.2,
\S5.4.4]{BernikDodson-1999}. In this paper we investigate the
distribution of rational points near arbitrary analytic
non-degenerate submanifold of $\R^n$ for all $n>1$. Analytic
non-degenerate manifolds are natural to consider as they run through Diophantine approximation and beyond. Recall that a connected analytic submanifold $\cM$ of $\R^n$
is \emph{non-degenerate}\/ if $\cM$ is not contained in a proper
affine subspace of $\R^n$. If $\cM$ is immersed by an analytic map
$\overline\xi=(\xi_1,\dots,\xi_n):U\to\R^n$ defined on a ball $U\subset\R^d$
then $\cM$ is non-degenerate if and only if the functions
$1,\xi_1,\dots,\xi_n$ are linearly independent over $\R$.

Throughout $m=\codim\cM\ge1$. Then we have the following obvious `volume based'
\begin{equation}\label{e:001}
\textbf{Heuristic estimate:}\qquad\qquad N(Q,\varepsilon) \asymp
\varepsilon^{\ccc} Q^{n+1},\hspace*{0.2\textwidth}
\end{equation}
where $\asymp$ means both $\ll$ and $\gg$. In order to gain some insight into when the heuristic estimate (\ref{e:001}) could potentially be
true we now consider the following two counterexamples.

\begin{example}\rm\label{examp1+}
Let $\cM=\{(x_1,\dots,x_n)\in\R^n:x_1^2+x_2^2=3\}$.
Obviously, $\cM$ is non-degenerate. It is readily verified that
$\cM\cap\Q^n=\varnothing$. Further, if $\varepsilon=o(Q^{-2})$ and $Q$ is
large enough, the rational points contributing to
$N(Q,\varepsilon)$ must lie on $\cM$, resulting in $N(Q,\varepsilon)=0$ for $\ve=o(Q^{-2})$.
This example can be extended to submanifolds of any codimension by using Pyartli's slicing technique \cite{Pyartli-1969}. The next example is of a different nature.
\end{example}

\begin{example}\rm\label{examp1}
Let
$\cM=\{(x_1,\dots,x_{d-1},x_d,x_d^2,\dots,x_d^{m+1})\in\R^n:\max\limits_{1\le
i\le d}|x_i|<1\}$, where $d\ge 2$. Clearly $\cM$ is non-degenerate and bounded.
Given a positive integer $q\le Q$, the rational points $\vv p/q$
with $\vv p=(p_1,\dots,p_{d-1},0,\dots,0)\in\Z^n$ obviously lie on
$\cM$. The number of such points is $\asymp Q^{d}$, thus
implying $N(Q,\varepsilon)\gg Q^{d}$ regardless of the size of
$\varepsilon$. The latter is significantly larger than the heuristic estimate (\ref{e:001})
unless $\varepsilon\gg Q^{-(m+1)/m}$.
\end{example}

In this paper we shall show that the condition $\varepsilon\gg
Q^{-(m+1)/m}$ is sufficient to prove the heuristic lower bound for
$N(Q,\varepsilon)$. Also we shall see in \S\ref{ext} that this
condition can be significantly relaxed when $\cM$ is a curve. The
results will be presented in a form convenient for the applications
in metric Diophantine approximation that we have in mind -- see
\S\ref{further_motivation}. Furthermore, the form of their
presentation reveals the distribution of rational points in
question, which is far more delicate than simply counting.

We will naturally and non-restrictively work with manifolds $\cM$
locally. Then, in view of the Implicit Function Theorem, this allows
us to represent $\cM$ by Monge parameterisations. Therefore without loss
of generality, we can assume that
\begin{equation}\label{e:002}
 \cM:=\big\{(x_1,\dots,x_\ddd,f_1(\vv x),\dots,f_\ccc(\vv x))\in\R^n:\vv x=(x_1,\dots,x_\ddd)\in
 U\big\}\,,
\end{equation}
where $U$ is an open subset of $\R^{\ddd}$ and $\vv
f=(f_1,\dots,f_{\ccc}):U\to\R^{\ccc}$ is a map. Here and elsewhere
$\ddd=\dim\cM$ and $\ccc=\codim\cM$. The distribution of rational points near the manifold (\ref{e:002}) is then conveniently described in terms of the set
$$
  \cR^\delta(Q,\psi,B) \ := \ \left\{ (q,\vv a,\vv b)\in \N\times\Z^\ddd\times\Z^\ccc \ :
\begin{array}{l}
 \vv a/q\in B,\ \delta Q< q\le Q \\[0.2ex]
 |qf_l(\vv a/q)-\vv b|_\infty\le\psi\\[0.2ex]
 \gcd(q,\vv a,\vv b)=1
\end{array}
  \right\},
$$
where $Q>1$, $\psi\ge0$, $\delta\ge0$, $B\subset U$ and
$|\cdot|_\infty$ denotes the supremum norm. Also define
$$
\Delta^{\delta_0}(Q,\psi,B,\rho):=\bigcup_{(q,\vv a,\vv b)\in\cR^{\delta_0}(Q,\psi,B)}
B\big(\vv a/q,\,\rho\big),
$$
where $B(\vv x,\rho)$ denotes a ball centred at
$\vv x$ of radius $\rho$. Roughly speaking, the set $\Delta^{\delta_0}(Q,\psi,B,\rho)$ indicates which part of the manifold can be covered by balls of radius $\asymp\rho$ centered at the rational points of interest. The following key result of this paper shows that this part is substantial for a suitable choice of parameters.
In what follows $\mu_{\ddd}$ denotes $\ddd$-dimensional Lebesgue measure.

\begin{theorem}\label{t:01}
Let the manifold $(\ref{e:002})$ be analytic and non-degenerate and
let $B_0\subset U$ be a compact ball. Then there are
absolute positive constants $k_0$, $\rho_0$ and $\delta_0$ depending
on $B_0$ only with the following property. For any ball $B\subset
B_0$ there are positive constants $C_0=C_0(B)$ and $Q_0=Q_0(B)$ such
that for all $Q\ge Q_0$ and all $\psi$ satisfying
\begin{equation}\label{e:003}
    C_0Q^{-1/\ccc }< \psi < C_0^{-1}
\end{equation}
we have
\begin{equation}\label{e:004}
\mu_{\ddd}\left(\Delta^{\delta_0}(Q,\psi,B,\rho)\cap B\,\right) \ \ge \ k_0 \,
\mu_{\ddd}(B) \, ,
\end{equation}
where $\rho:=\rho_0\times(\psi^{\ccc }Q^{\ddd+1})^{-1/\ddd}.$
\end{theorem}

\begin{corollary}\label{corollary:01}
Let $\cM$ and $B_0$ be as in Theorem~\ref{t:01}. Then, there are
constants $\delta_0$ and $k_1>0$ such that for any ball $B\subset
B_0$ there exist $Q_0>0$ and $C_0>0$ such that for all $Q\ge
Q_0$ and all $\psi$ satisfying $(\ref{e:003})$ we have that
\begin{equation}\label{e:005}
N^{\delta_0}(Q,\psi,B)\ge k_1 \psi^{\ccc} Q^{\ddd+1}\mu_{\ddd}(B).
\end{equation}
\end{corollary}

\noindent\textit{Proof of Corollary~$\ref{corollary:01}$.} For any
$\vv r\in\R^{\ddd}$ we obviously have that $\mu_{\ddd}(B(\vv
r,\rho)\cap B)\le V_{\ddd}\rho^{\ddd}$, where $V_{\ddd}$ is the volume of a $\ddd$-dimensional ball of radius 1. Therefore, the
r.h.s.\! of (\ref{e:004}) ~(\emph{throughout r.h.s.\!\! means right
hand side}) is bounded above by
$N^{\delta_0}(Q,\psi,B)V_{\ddd}\rho^{\ddd}$. By (\ref{e:004}), we
get that $N^{\delta_0}(Q,\psi,B)\ge
V_{\ddd}^{-1}\rho^{-\ddd}k_0\mu_{\ddd}(B)$. Substituting the value
of $\rho$ from Theorem~\ref{t:01} into the last inequity completes
the proof.\proofend

\begin{remark}\rm
Clearly, every rational point $(\vv a/q,\vv b/q)$
arising from $\cR^{\delta_0}(Q,\psi,B)$ lies within the distance
$\varepsilon=\psi(\delta_0Q)^{-1}$ from $\cM$. Thus,
$
N(Q,\varepsilon)\ge N^{\delta_0}(Q,\varepsilon\delta_0Q,B_0).
$
By Corollary~\ref{corollary:01}, we get the lower bound $
N(Q,\delta)\gg \varepsilon^{\ccc} Q^{n+1} $ valid for $\varepsilon\gg
Q^{-(m+1)/m}$ consistent with (\ref{e:001}).
\end{remark}

\begin{remark}\rm
In the case of hypersurfaces $\ccc=1$. Therefore, the
condition $\varepsilon\gg Q^{-(m+1)/\ccc}$ transforms into $\varepsilon\gg
Q^{-2}$. This is the same as for planar curves
\cite{Beresnevich-Dickinson-Velani-07:MR2373145}. It tells us that
rational points with denominator $q\le Q$ can get
\textsc{const}$\times Q^{-2}$ close to an arbitrary analytic
non-degenerate hypersurface. In fact, in view of
Example~\ref{examp1+} this is generically best possible!
\end{remark}

\begin{remark}\rm
In the case of planar curves the lower bound (\ref{e:005}) has
already been established in
\cite[Theorem~6]{Beresnevich-Dickinson-Velani-07:MR2373145}.
However, in that paper the constant $k_1$ happens to dependent on
$B$, while in this paper $k_1$ is uniform.
\end{remark}

\section{Diophantine approximation on manifolds}\label{further_motivation}

In this section we apply Theorem~\ref{t:01}
to simultaneous Diophantine approximation on manifolds.
Traditionally, problems on the proximity of rational points to points in $\R^n$ assume finding optimal relations between the accuracy of approximation and the `height' of approximating rational points $\vv p/q$. In our case, the latter is measured by $q$ while the former is measured by $\psi/q$. Therefore, throughout this section $\psi:\N\to\Rp$ will be regarded as a decreasing function referred to as an \emph{approximation function}, where $\Rp=(0,+\infty)$. Given $\tau>0$, the approximation function $q\mapsto q^{-\tau}$ will be denoted by $\psi_\tau(q)$.

The point $\vv y\in\R^n$ is called
\emph{$\psi$-approximable}\/ if there are infinitely many $q\in\N$
satisfying
\begin{equation}\label{e:006}
    \|q\vv y\|<\psi(q)\,,
\end{equation}
where $\|q\vv y\|$ denotes the distance of $q\vv y$ from $\Z^n$ with
respect to the sup-norm $|\cdot|_\infty$. Throughout, $\cS_n(\psi)$
denotes the set of $\psi$-approximable points in $\R^n$.

\medskip

By Dirichlet's theorem (see, e.g., \cite{Schmidt-1980}), $\cS_n\big(\psi_{1/n}\big)=\R^n$. The points $\vv y\in\R^n$ such that $\vv y\not\in\cS_n(\psi_\tau)$ for any $\tau>1/n$ are called \emph{extremal}. A relatively easy consequence of the Borel-Cantelli lemma is that almost all points in $\R^n$ are extremal -- see, e.g., \cite{BernikDodson-1999}. The property of extremality is fundamental in Diophantine approximation. For example, Roth's celebrated theorem establishes nothing but the extremality of irrational algebraic numbers. Within this paper we will be dealing with problems that go back to the profound conjecture of Mahler \cite{Mahler-1932b} that almost all points on the Veronese curves $(x,\dots,x^n)$ are extremal. The problem was studied in depth for over 30 years and eventually settled in full by Sprind\v zuk in 1964 (see \cite{Sprindzuk-1969-Mahler-problem}) who also stated the following general conjecture \cite{Sprindzuk-1980-Achievements}:

\medskip

\noindent\textbf{Conjecture (Sprind\v zuk)\,:} \emph{Any analytic
non-degenerate submanifold of $\R^n$ is extremal.}

\medskip

\noindent Formally a differentiable manifold $\cM\subset\R^n$ is called \emph{extremal} if almost all points of $\cM$ (with respect to the induced Lebesgue measure on $\cM$) are extremal. For $n=2$ the conjecture is a consequence of Schmidt's theorem \cite{Schmidt-64:MR0171753} and for $n=3$ it has been proved by Bernik and the author \cite{Beresnevich-Bernik-96:MR1387861}. The full conjecture (with the analyticity assumption dropped) has been established by Kleinbock and Margulis in the tour de force \cite{Kleinbock-Margulis-98:MR1652916} and later re-established in \cite{Beresnevich-02:MR1905790} using different techniques. The work of Kleinbock and Margulis has also dealt with the far more delicate multiplicative case known as the Baker-Sprind\v zuk conjecture and led to a surge of activity that led to establishing the extremality of various classes of manifolds and sets -- see, for example, \cite{Kleinbock-03:MR1982150, Kleinbock-04:MR2094125, Kleinbock-Lindenstrauss-Weiss-04:MR2134453, Kleinbock-Tomanov-07:MR2314053}.

The following two major problems now arise (see, e.g.,
\cite[\S1]{Beresnevich-Dickinson-Velani-07:MR2373145} or
\cite[\S6]{Beresnevich-Bernik-Kleinbock-Margulis-02:MR1944505}):

\begin{problem}\label{p2}
To develop a Khintchine type theory for $\cS_n(\psi)\cap\cM$.
\end{problem}

\begin{problem}\label{p3}
To develop a Hausdorff measure/dimension theory for
$\cS_n(\psi)\cap\cM$.
\end{problem}
The goal of Problem~\ref{p2} is a metric theory of
$\cS_n(\psi)\cap\cM$ with $\psi$ being a general approximation
function, not just $\psi_\tau(q)=q^{-\tau}$
associated with extremality. The goal of Problem~\ref{p3} is to determine
the `size' of $\cS_n(\psi)\cap\cM$ via Hausdorff measure and dimension.

Before we proceed with the more detailed discussion of the above problems, it is worth mentioning that there are dual versions
of Problems~\ref{p2} and \ref{p3}. In the dual case the approximating objects are rational hyperplanes rather than rational points. The problems in the dual case are much more tractable and progress has been
significantly better. In particular, the dual version of Problem~\ref{p2}
has been fully settled \cite{Beresnevich-02:MR1905790,
Beresnevich-Bernik-Kleinbock-Margulis-02:MR1944505,
Bernik-Kleinbock-Margulis-01:MR1829381} and very deep answers
regarding the dual version of Problem~\ref{p3} found
\cite{Beresnevich-Bernik-Dodson-02:MR2069553,
Beresnevich-Dickinson-Velani-06:MR2184760, Bernik-1983a,
DickinsonDodson-2000a, DodsonRynneVickers-1989b}. However, as we shall see, Problems~\ref{p2}
and \ref{p3} (non-dual) have more or less been understood only in
$\R^2$.

\mysubsection{Khintchine type theory}\label{ktt}%

Let $\cM\subset\R^n$ be a manifold. If for any  approximation function $\psi:\N\to\Rp$ such that
\begin{equation}\label{e:007}
\sum_{q\in\Z}\psi(q)^n
\end{equation}
converges \emph{almost no} point on $\cM$ is $\psi$-approximable then $\cM$ is called of \emph{Khintchine type for convergence}. In turn, $\cM$ is called of \emph{Khintchine type for divergence} if
for any approximation function $\psi$ such that the sum
(\ref{e:007}) diverges \emph{almost all} points on $\cM$ are
$\psi$-approximable. This terminology represents a zero-one law and
has been introduced in \cite{BernikDodson-1999} to acknowledge the
fundamental contribution of Khintchine who discovered this beautiful
law in the case $\cM=\R^n$ \cite{Khintchine-1924, Khintchine-1926}.
We now discuss the state of the art for proper submanifolds of
$\R^n$.

\medskip

\textit{Planar curves $(n=2)$.} The story has begun with
the pioneering work \cite{Bernik-1979} of Bernik who showed that the parabola $(x,x^2)$ is of Khintchine type for convergence. Subsequently,
working towards a conjecture of Alan Baker, Mashanov
has established a multiplicative
analogue of Bernik's result \cite{Mashanov-87:MR888593}. There has been no progress with planar
curves since then, until Dickinson, Velani and the author
have shown that any
$C^{(3)}$ non-degenerate planar curve is of Khintchine type for
divergence \cite{Beresnevich-Dickinson-Velani-07:MR2373145} and subsequently Vaughan and Velani
have established that any $C^{(2)}$
non-degenerate planar curve is of Khintchine type for convergence \cite{Vaughan-Velani-2007}.
See also \cite{Badziahin-Levesley-07:MR2347267,
Beresnevich-Vaughan-Velani-08-Inhom,
Beresnevich-Velani-07:MR2285737} for further progress.

\medskip

\emph{Higher dimensions $(n>2)$.} In this case the Khintchine type theory also exists but is rather bizarre. Bernik \cite{Bernik-73:MR0337863, Bernik-77:MR0480402} has shown that the manifolds in $\R^{mk}$ given as the cartesian product of $m$ non-degenerate curves in $\R^k$ are of Khintchine type for convergence if $m\ge k$ and for divergence if $k=2$ and $m\ge 4$. Dodson, Rynne and Vickers \cite{DodsonRynneVickers-1991a, Dodson-Rynne-Vickers-1996} have found Khintchine type manifolds satisfying certain curvature conditions. However, these conditions significantly constrain the dimension of the manifolds and completely rule out curves. For example, the Khintchine type manifolds of \cite{DodsonRynneVickers-1991a, Dodson-Rynne-Vickers-1996}
assume that $d=\dim\cM\ge \max\{2,\sqrt{2n}-\frac32\}$ for convergence and $d\ge\frac34(n+5)$ \& $n\ge19$ for divergence. Thus, the simplest example of a Khintchine type manifold for divergence could only be an 18-dimensional surface in $\R^{19}$. It should be noted that Dodson, Rynne and Vickers established their divergence Khintchine type theorem in the quantitative form. Assuming a condition on $\psi$ which implies that $\cS_n(\psi)=\R^n$, Harman \cite{Harman-03:MR1979906} has obtained a quantitative result for Veronese curves and manifolds that are known to be of Khintchine type for convergence. Recently Gorodnik and Shah \cite{Gorodnik-Shah-08} have obtained a Khintchine type theorem for the quadratic varieties $x_1^2\pm\dots\pm x_d^2=1$ with the approximating rational points being of a special type. The Khintchine type theory for curves in dimensions $n>2$ is simply non-existent. However, in view of Pyartli's slicing technique \cite{Pyartli-1969}, curves underpin the whole theory. The following result of this paper covers arbitrary non-degenerate analytic curves as well as arbitrary non-degenerate analytic submanifolds of $\R^n$:

\begin{theorem}\label{t:02}
For any $n\ge2$ any non-degenerate analytic submanifold of\/ $\R^n$
is of Khintchine type for divergence.
\end{theorem}

\noindent\textit{Classical case.} In order to illustrate the
statement of Theorem~\ref{t:02}, let us consider the following
classical problem on rational approximations to consecutive powers
of a real number. That is, we consider the inequality
\begin{equation}\label{e:008}
    \max\big\{\,\|qx\|,\|qx^2\|,\dots,\|qx^n\|\,\big\} <
    \psi(q).
\end{equation}
Since the consecutive powers of $x$ are real analytic functions of
$x$ which, together with $1$, are linearly independent over $\R$,
Theorem~\ref{t:02} implies the following
\begin{corollary}\label{cor:01}
Given any monotonic $\psi:\N\to\Rp$ such that the sum
$(\ref{e:007})$ diverges, for almost all $x\in\R$ inequality
$(\ref{e:008})$ has infinitely many solutions $q\in\N$.
\end{corollary}


\noindent In 1925 Khintchine \cite{Khintchine-1925} established such a statement in the special case when $\psi(q)=cq^{-1/n}$ with arbitrary but fixed $c>0$. The latter has been generalised by R.C.~Baker \cite{Baker-1976} to smooth manifolds but the same class of approximation functions. Corollary~\ref{cor:01} is thus the first improvement on that result of Khintchine in the period of over 80 years. It obviously contains Khintchine's result and is believed to be best possible. In fact, a folk conjecture suggests that for almost all $x\in\R$ there are only finitely many $q\in\N$ satisfying
(\ref{e:008}) provided that the sum (\ref{e:007}) converges.

\mysubsection{Hausdorff dimension and measure theory}\label{hdrt}%

Problem~\ref{p3} throws up a few surprises. For example, unlike the dual case the
dimension of $\cS_n(\psi)\cap\cM$ happens to depend on the arithmetic properties of $\cM$. To grasp the ideas consider the following
popular example. Let $\cC_r$ be the circle $x^2+y^2=r$. It is easily
verified that if $r\in\N$, $\tau>1$ and $\psi(q)=\psi_\tau(q)=q^{-\tau}$ then all the rational
points implicit in (\ref{e:006}) must lie on $\cC_r$ for
sufficiently large $q$. For the unit circle $\cC_1$ these points are
parameterised by Pythagorean triples and well understood. As a
result
\begin{equation}\label{e:009}
\dim\cS_2(\psi_\tau)\cap\cC_1=\frac{1}{\tau+1}\qquad\text{for $\tau>1$},
\end{equation}
where $\dim$ stands for Hausdorff dimension. The fact (\ref{e:009})
has been established in two complementary papers by Melnichuk
\cite{Melnichuk-1979a} and Dickinson \& Dodson
\cite{DickinsonDodson-2001a}. On the other hand, it is easily seen
that $\cC_3\cap\Q^2=\varnothing$. Consequently
\begin{equation}\label{e:010}
\dim\cS_2(\psi_\tau)\cap\cC_3=0\qquad\text{for $\tau>1$}.
\end{equation}
Thus, scaling $\cC_1$ by $\sqrt 3$ completely changes the character of the set of $\psi_\tau$-approximable points lying on it. Luckily,
this cannot happen if $\tau< 1$. In fact, as shown in
\cite{Beresnevich-Dickinson-Velani-07:MR2373145}
\begin{equation}\label{e:011}
 \dim\cS_2(\psi_\tau)\cap\cC=\frac{2-\tau}{\tau+1}\quad\text{ when }
1/2\le\tau<1
\end{equation}
for all $C^{(3)}$ curves $\cC$ in $\R^2$ non-degenerate everywhere
except possibly on a set of Hausdorff dimension $\le \frac{2-\tau}{\tau+1}$. The Hausdorff
dimension of $\cS_2(\psi)\cap\cC$ has also been found in
\cite{Beresnevich-Dickinson-Velani-07:MR2373145} for general
approximation functions $\psi$. Furthermore, an analogue of
\Jarnik{}'s theorem \cite{Jarnik-1931} has been established in
\cite{Beresnevich-Dickinson-Velani-07:MR2373145} and
\cite{Vaughan-Velani-2007} which provides a complete picture of the
$s$-dimensional Hausdorff measure of $\cS_2(\psi)\cap\cC$ -- see
\cite{Beresnevich-Dickinson-Velani-07:MR2373145,
Vaughan-Velani-2007} for details.

\medskip

\noindent\textit{Higher dimensions.} Khintchine's transference
principle \cite{Schmidt-1980} can be used to deduce bounds on
$\dim\cS_n(\psi_\tau)\cap\cM$ from the much better
understood dual case.
Although the bounds obtained this way are rather crude, until
recently nothing else was known. In
\cite{Drutu-05:MR2195121} Drutu established a comprehensive theory for non-degenerate rational quadrics in $\R^n$ when
the approximating rational points lie on quadrics. In particular, her results
include (\ref{e:009}) and (\ref{e:010}) as two special cases. More
recently Budarina and Dickinson
\cite{Budarina-Dickinson-07:MR2397137} have investigated
$\cS_n(\psi_\tau)\cap\cM$ for hypersurfaces $\cM$ in $\R^n$ parameterised
by the forms $x_1^d+\dots+x_{n-1}^d$ of degree $d<\log n$, the
exponent $\tau$ being large and the approximating rational points
being lying on $\cM$. However, except for planar curves, the approximating rational points always lie on the
manifold. In view of this, Theorem~\ref{t:03} appears to be the first general result concerning Problem~\ref{p3} in dimensions
$n>2$.

Let $\cH^s$ denote $s$-dimensional Hausdorff measure. In order
to state the result we now introduce the \emph{exponent of $\psi$}
also known as the lower order of $1/\psi$ at infinity:
$$
\tau(\psi):=\liminf_{q\to\infty}\frac{-\log \psi(q)}{\log q}.
$$

\begin{theorem}\label{t:03}
Let $\cM$ be a non-degenerate analytic submanifold of\/ $\R^n$,
$\ddd=\dim\cM$ and $\ccc=\codim\cM$. Thus, $d+m=n$. Let $\psi:\N\to\Rp$ be a
monotonic function such that $q\psi(q)^{\ccc} \to\infty$ as
$q\to\infty$. Then for any
$s\in\big(\tfrac{\ccc}{\ccc+1}{\ddd},\ddd\big)$
\begin{equation}\label{e:012}
 \cH^s(\cS_n(\psi)\cap\cM)=\infty \qquad \text{if} \qquad \sum_{q=1}^\infty \
 q^n\Big(\frac{\psi(q)}q\Big)^{s+\ccc }=\infty.
\end{equation}
Consequently if $\tau=\tau(\psi)$ satisfies $1/n<\tau<1/\ccc $ then
\begin{equation}\label{e:013}
\dim \cS_n(\psi)\cap\cM \ge s_0:=\frac{n+1}{\tau+1}-\ccc .
\end{equation}
\end{theorem}

We shall see in \S\ref{ext} that for non-degenerate analytic curves ($d=1$) Theorem~\ref{t:03} holds for $s\in(d/2;d)$. It is also possible to obtain the version of Theorem~\ref{t:03} that would incorporate generalised Hausdorff measures. We opt to omit further details which can be easily recovered using the ideas of \cite[\S8.1]{Beresnevich-Dickinson-Velani-07:MR2373145} where the case $n=2$ is considered.

\mysubsection{Proof of Theorems~\ref{t:02} and
\ref{t:03}}\label{proof23}

The proof below generalises the arguments given in \S\S3,6,7 of \cite{Beresnevich-Dickinson-Velani-07:MR2373145} to higher dimensions.

\medskip

\noindent\textbf{Note 1:} {Within Theorem~$\ref{t:03}$ it suffices
to establish $(\ref{e:012})$ for $(\ref{e:013})$ follows from
$(\ref{e:012})$.}

\medskip

\noindent\textit{Proof.}
By the definition of $\tau(\psi)$, for any $\ve>0$ there are infinitely many $q$ such that $\psi(q)\ge q^{-\tau-\ve}$. Since $\psi$ is monotonic, $\psi(2^t)\ge 2^{-(t+1)(\tau+\ve)}$ for $t\in\Z$ satisfying $2^t \le q \le 2^{t+1}$. Therefore, there are infinitely many $t\in\N$ such that $\psi(2^t)\ge 2^{-(t+1)(\tau+\ve)}$. Hence, on taking $s=\frac{n+1}{\tau+1+\ve}-\ccc $ with $\ve>0$, one verifies that $ 2^{t(n+1)}(\psi(2^t)2^{-t})^{s+\ccc }\ge2^{-(n+1)}. $ The latter holds for infinitely many $t$ and implies that $\sum_{t=1}^\infty \ 2^{t(n+1)}(\psi(2^t)2^{-t})^{s+\ccc }=\infty$. Due to the monotonicity of $\psi$ this further implies that the sum in (\ref{e:012}) diverges and therefore, by (\ref{e:012}), $ \cH^s(\cS_n(\psi)\cap\cM)=\infty$. By the definition of Hausdorff dimension, we deduce that $\dim \cS_n(\psi)\cap\cM \ge s=\frac{n+1}{\tau+1+\ve}-\ccc$, whence (\ref{e:013}) readily follows.
\proofend

\medskip

\noindent\textbf{Note 2:} {The condition
\begin{equation}\label{e:014}
    \lim_{q\to\infty}q\psi(q)^{\ccc} =\infty,
\end{equation}
which is a part of Theorem~\ref{t:03}, can be
assumed in the proof of Theorem~\ref{t:02}.}

\medskip

\noindent\textit{Proof.}
To verify (\ref{e:014}) consider the monotonic function
$\psi_1(q)=\max\{q^{-2/(2n-1)},\psi(q)\}$. Then the divergence of
(\ref{e:007}) implies $\sum_{q=1}^\infty\psi_1(q)^n=\infty$.
Obviously
$\cS_n(\psi_1\,)=\cS_n(\psi)\cup\cS_n(2/(2n-1))$. Since
$2/(2n-1)>1/n$ and every non-degenerate submanifold of $\R^n$ is
extremal we obviously have that the set $\cM\cap\cS_n(2/(2n-1))$ has
zero measure on $\cM$. Hence $\cM\cap\cS_n(\psi_1\,)$ and
$\cM\cap\cS_n(\psi)$ are of the same measure and $\psi$ can be
replaced with $\psi_1$, which satisfies (\ref{e:014}).
\proofend

\medskip

\noindent\textbf{Note 3:} In view of the metric nature of
Theorems~\ref{t:02} and \ref{t:03} it is enough to consider a
sufficiently small neighborhood of an arbitrary point on $\cM$.
Therefore, by the Implicit Function Theorem, without loss of
generality we can assume that $\cM$ is of the Monge form
(\ref{e:002}) and that the functions $f_1,\dots,f_{\ccc}$ are
Lipschitz; that is, for some $c_1\ge1$
\begin{equation}\label{e:015}
  \max_{1\le l\le \ccc }|f_l(\vv x)-f_l(\vv x')|\le c_1|\vv x-\vv x'|_\infty\quad\text{for all
  }\vv x,\vv x'\in U.
\end{equation}

\noindent\textbf{Note 4:} Let $\cS_{\vv f}(\psi)$ be the set of $\vv
x\in U$ such that $(\vv x,\vv f(\vv x))\in\cS_n(\psi)$. Obviously,
$\cS_{\vv f}(\psi)$ is the orthogonal projection of
$\cS_n(\psi)\cap\cM$ onto $\R^{\ddd}$. By (\ref{e:015}), $\cS_{\vv
f}(\psi)$ and $\cS_n(\psi)\cap\cM$ are related by a bi-Lipschitz map
and therefore $\cS_{\vv f}(\psi)$ is of full Lebesgue measure in $U$
if and only if $\cS_n(\psi)\cap\cM$ is of full induced Lebesgue
measure on $\cM$ -- see \cite[\S1.5.1]{BernikDodson-1999}. Further,
recall that $d$-dimensional Lebesgue measure is comparable to
$\cH^d$. Therefore, to prove Theorem~\ref{t:02} it suffices to show
that for every compact ball $B_0$ in $U$
\begin{equation}\label{e:016}
 \cH^d(\cS_{\vv f}(\psi)\cap B_0)=\cH^d(B_0) \quad\text{if}\quad
 \sum_{q=1}^\infty \psi(q)^{n}=\infty.
\end{equation}
Similarly one can show that Theorem~\ref{t:03} follows on showing
that
\begin{equation}\label{e:017}
 \cH^s(\cS_{\vv f}(\psi)\cap B_0)=\cH^s(B_0) \quad\text{if}\quad
 \sum_{q=1}^\infty q^n\left(\frac{\psi(q)}{q}\right)^{s+\ccc }=\infty
\end{equation}
holds for every compact ball $B_0$ in $U$ and
$s\in(\ccc\ddd/(\ccc+1),\ddd)$. Note that for $s<d$,
$\cH^s(B_0)=\infty$. Also note that in the case $s=d$, (\ref{e:017})
is simply (\ref{e:016}).

\bigskip

\noindent \textbf{Upshot:} on establishing (\ref{e:017}) for
$s\in\big(\frac{\ccc}{\ccc+1}\ddd,\ddd\big]$ and $\psi$ satisfying
(\ref{e:014}) we prove Theorems~\ref{t:02} and \ref{t:03}.

\bigskip

\noindent\emph{Ubiquitous systems}. In what follows we will use the ubiquitous systems technique. The notion of ubiquity introduced below is equivalent to that of \cite{Beresnevich-Dickinson-Velani-06:MR2184760} in the setting that is now to be described. Let $B_0$ be a ball in $\R^{\ddd}$ and $\cR:=(\ra)_{\al\in J}$ be a family of points $\ra$ in $B_0$ (usually called \emph{resonant points}) indexed by a countable set $J$. Let $\beta:J\to \Rp:\alpha\mapsto\ba$ be a function on $J$, which attaches a `weight' $\ba$ to points $\ra$. For $t\in\N$ let $J(t):=\{\al \in J:\ba\le 2^t\}$ and assume $J(t)$ is always finite.

\begin{definition}\label{rs}\label{US}\rm
Let $\rho: \Rp \to\Rp$ be a function such that
$\lim_{t\to\infty}\rho(t)=0$. The system $(\cR;\beta)$ is called
{\em locally ubiquitous in $B_0$ relative to $\rho$} if there is an
absolute constant $k_0>0$ such that for any ball $B\subset B_0$
\begin{equation}\label{e:018}
\liminf_{t\to\infty}\,\,\mu_{\ddd}\Big(\,\bigcup_{\al\in J(t)}
B\big(\ra,\r(2^t)\big)\cap B\Big) \ \ge \  k_0 \, \mu_{\ddd}(B) \, .
\end{equation}
\end{definition}
Here as before $\mu_{\ddd}$ denotes Lebesgue measure in
$\R^{\ddd}$ and $B(\vv x,r)$ denotes the ball in $\R^{\ddd}$ centred
at $\vv x$ of radius $r$. The function $\rho$ is referred to as
\emph{ubiquity function}.

Given a function $\Psi:\Rp\to\Rp$, let
$$
\Lambda_\cR(\Psi) \ := \ \{\vv x\in B_0:|\vv x-\ra|_\infty<\Psi(\ba)
\ \mbox{holds for\ infinitely\ many\ }\al \in J \} \,.
$$
The following lemma follows from
Corollaries~2, 4 and 5 from
\cite{Beresnevich-Dickinson-Velani-06:MR2184760}. In the case
$\ddd=1$ a simplified proof of Lemma~\ref{l:01} is given in
\cite[Theorems~9~and~10]{Beresnevich-Dickinson-Velani-07:MR2373145}, see also \cite{Beresnevich-Velani-09}.

\begin{lemma}\label{l:01}
Let $\Psi:\Rp\to\Rp$ be a monotonic function such that for some
$\lambda<1$, $\Psi(2^{t+1})\le \lambda \Psi(2^t)$ holds for $t$
sufficiently large. Let $(\cR,\beta)$ be a locally ubiquitous system
in $B_0$ relative to $\rho$. Then for any $s \in (0,\ddd\,]$
\begin{equation}\label{e:019}
 \cH^s\big(\Lambda_\cR(\Psi)\big) \
= \
\cH^s(B_0)\qquad\text{if}\qquad\sum_{t=1}^{\infty}\frac{\Psi(2^t)^s}{\rho(2^t)^{\ddd}}\
= \ \infty \, .
\end{equation}
\end{lemma}

\bigskip

\noindent\textit{Proof of Theorem~\ref{t:02} and \ref{t:03}.} Recall
again that our goal is to establish (\ref{e:017}) for
$s\in(\ccc\ddd/(\ccc+1),\ddd]$ and approximation functions $\psi$
satisfying (\ref{e:014}), where $B_0$ is an arbitrary
non-empty compact ball in $U$. Therefore, for the rest of this section we
fix such a $B_0$. Also recall that the map $\vv f$ which arises from
(\ref{e:002}) satisfies the Lipschitz condition (\ref{e:015}). We
can also assume that $\lim_{q\to\infty}\psi(q)=0$ as otherwise
$\cS_n(\psi)=\R^n$ and there is nothing to prove.

\medskip

We first construct a ubiquitous system relevant to our main goal.
Let $\rho_0$ and $\delta_0$ be the same as in Theorem~\ref{t:01}.
Define the ubiquity function $\rho(q)=\rho_0\times(\psi(q)^{\ccc}
q^{\ddd+1})^{-1/\ddd}$ and the sequence $\cR:=\{\vv a/q\}_{(q,\vv
a)\in J}$ of resonant points in $B_0$, where
$$
 J:= \big\{(q,\vv a)\in\N\times\Z^{\ddd}\,: \, \vv a/q\in B_0,\
 \displaystyle \max_{1\le l\le \ccc }\|qf_l(\vv a/q)\|\le\tfrac{1}{2}\psi(q) \big\}.
$$
For $\alpha=(q,\vv a)\in J$ define $\beta_\alpha:=q$. We prove the
following
\begin{lemma}\label{l2.3}
Assume that Theorem~\ref{t:01} holds. Then, with $B_0$, $\cR$,
$\beta$ and $\rho$ as above, the system $(\cR,\beta)$ is locally
ubiquitous in $B_0$ relative to $\rho$.
\end{lemma}

\noindent\textit{Proof.} First of all, by (\ref{e:014}), $\rho(q)\to0$ as $q\to\infty$. We now verify (\ref{e:018}) for
the specific choice of $\cR$, $\beta$ and $\rho$ we have made.
Obviously $J(t)$ consists of $(q,\vv a)\in J$ such that $q\le
Q:=2^t$. Fix an arbitrary ball $B\subset B_0$ and consider the union
in (\ref{e:018}). This union contains
\begin{equation}\label{e:020}
\bigcup_{\delta_0 Q \le q\le Q }\ \bigcup_{\vv
a\in\Z^{\ddd}\,:\,(q,\vv a)\in J} B\big(\vv a/q,\r(Q)\big)\cap B\
\supset \Delta^{\delta_0}\big(Q,\tfrac12\psi(Q),B,\rho(Q)\big)\cap B\,,
\end{equation}
where $\Delta^{\delta_0}(\,\cdot\,)$ is the set defined in \S\ref{rpnm}
and appearing in Theorem~\ref{t:01}. By (\ref{e:014}) and the
assumption $\lim_{q\to\infty}\psi(q)=0$, conditions (\ref{e:003})
are met for sufficiently large $Q$ and therefore, by
Theorem~\ref{t:01}, the $\mu_{\ddd}$-measure of the sets in
(\ref{e:020}) is at least $k_0\mu_{\ddd}(B)$. Therefore
(\ref{e:018}) is fulfilled and the proof is complete. \proofend

\medskip

In the next two statements we establish a relation between
$\Lambda_\cR(\Psi)$ and $\cS_{\vv f}(\psi)$ and an analogue of
(\ref{e:017}) in terms of $\Lambda_\cR(\Psi)$.
\begin{lemma}\label{l2.5}
Let $\Psi(q)=\psi(q)/(2c_1q)$, where $c_1$ arises from (\ref{e:015})
and let $B_0$, $\cR$, $\beta$ and $\rho$ be as in Lemma~\ref{l2.3}.
Then $\Lambda_\cR(\Psi)\subset \cS_{\vv f}(\psi)$.
\end{lemma}

\noindent\textit{Proof.} Assume that $\vv
x=(x_1,\dots,x_{\ddd})\in\Lambda_{\cR}(\Psi)$. Then
\begin{equation}\label{e:021}
|\vv x-\vv a/q|_\infty<\Psi(q)=\psi(q)/(2c_1q)
\end{equation}
for infinitely many $(q,\vv a)\in\N\times\Z^{\ddd}$ such that
\begin{equation}\label{e:022}
    \max_{1\le l\le \ccc }|qf_l(\vv a/q)-b_l|\le\tfrac{1}{2}\psi(q)
\end{equation}
for some $\vv b=(b_1,\dots,b_{\ccc} )\in\Z^{\ccc} $. By the
triangular inequality,
\begin{equation}\label{e:023}
\begin{array}[b]{rcl}
  |f_l(\vv x)-b_l/q| & \le & |f_l(\vv x)-f_l(\vv a/q)|+|f_l(\vv a/q)-b_l/q|
   \\[1ex]
   & \stackrel{\eqref{e:015}}{\le} & c_1|\vv x-\vv a/q|_\infty+|f_l(\vv a/q)-b_l/q| \\[1ex]
   & \stackrel{\eqref{e:021}\&(\ref{e:022})}{<} & c_1\cdot\psi(q)/(2c_1q)+\frac12\psi(q)/q= \psi(q)/q.
\end{array}
\end{equation}
Since (\ref{e:021}) and (\ref{e:023}) hold for infinitely many $q$, we have that $(\vv x,\vv f(\vv x))\in\cS_n(\psi)$; that is $\vv x$
belongs to $\cS_{\vv f}(\psi)$. Therefore, $\Lambda_\cR(\Psi)\subset
\cS_{\vv f}(\psi)$. \proofend

\begin{lemma}\label{l2.4}
Assume that Theorem~\ref{t:01} holds. Let $\Psi(q)=\psi(q)/(2c_1q)$,
where $c_1$ arises from (\ref{e:015}) and let $B_0$, $\cR$, $\beta$
and $\rho$ be as in Lemma~\ref{l2.3}. Then
\begin{equation}\label{e:024}
\cH^s\big(\Lambda_\cR(\Psi)\big)=\cH^s(B_0)\quad\text{if}\quad
\sum_{q=1}^\infty q^n\left(\frac{\psi(q)}{q}\right)^{s+\ccc
}=\infty.
\end{equation}
\end{lemma}

\noindent\emph{Proof.} Since $\psi$ is decreasing, $\Psi(2^{t+1})\le
\lambda \Psi(2^t)$ with $\lambda=1/2$. Further, using the explicit
form for $\Psi$ and $\rho$ verify that
$$
\sum_{t=1}^\infty\frac{\Psi(2^t)^s}{\rho(2^t)^{\ddd}}\asymp
\sum_{t=1}^\infty\frac{\psi(2^t)^s2^{-st}}{\psi(2^t)^{-\ccc
}2^{-(\ddd+1)t}} \asymp
\sum_{t=1}^\infty\left(\frac{\psi(2^t)}{2^t}\right)^{s+\ccc
}2^{(n+1)t}\,.
$$
In view of the monotonicity of $\psi$ the latter sum diverges if and
only if $\sum\limits_{q=1}^\infty
q^n\left(\frac{\psi(q)}{q}\right)^{s+\ccc }$ diverges. Hence, by
Lemmas~\ref{l:01} and \ref{l2.3}, we get (\ref{e:024}). \proofend

\bigskip

We are now able to complete the proof of Theorems~\ref{t:02} and \ref{t:03}. Recall that we have to establish (\ref{e:017}). Let $\Psi$, $B_0$, $\cR$, $\beta$ and $\rho$ be as in Lemma~\ref{l2.4}. By the monotonicity of $\cH^s$, $\cH^s(\cS_{\vv f}(\psi)\cap B_0)\le\cH^s(B_0)$. Therefore, to establish (\ref{e:017}) it suffices to show that $\cH^s(\cS_{\vv f}(\psi)\cap B_0)\ge\cH^s(B_0)$ provided that the sum in (\ref{e:017}) diverges. In view of Lemma~\ref{l2.5} this follows from (\ref{e:024}) and the proof of Theorems~\ref{t:02} and \ref{t:03} modulo Theorem~\ref{t:01} is thus complete. \proofend

\section{Some auxiliary geometry}\label{aux}

The distance of a rational point from a manifold is conveniently studied using the notion of projective distance (due to H. and J. Weyl \cite{Weyl-Weyl-38:MR1503422}) which involves exterior and interior products. These classical and well established topics are now briefly recalled. The overview below is mostly taken from \cite{Schmidt-1980} and \cite{Whitney}. We will use the standard embedding of $\R^n$ into the real projective space $\P^n$. Given $\vv x=(x_1,\dots,x_n)\in\R^n$, the point $\mv x=(\lambda,\lambda x_1,\dots,\lambda x_n)\in\V$ with $\lambda\not=0$ will be referred to as the \emph{homogeneous coordinates} of $\vv x$.

\mysubsection{Exterior product and projective
distance}\label{sec:ext}

Throughout $\We^p(\R^{n+1})$ denotes the $p$-th exterior
power of $\V$ and ``$\we$'' denotes the exterior product. If $p\le
n+1$ and $\mv e_0,\dots,\mv e_n$ is a basis of $\V$, then the
multivectors
\begin{equation}\label{e:025}
\mv e_I:=\We_{i\in I}\ \mv e_i,\qquad I\in C(n+1,p)
\end{equation}
form a basis of $\We^p(\V)$, where $C(n+1,p)$ denotes the set of all
subsets of $\{0,\dots,n\}$ of cardinality $p$. The following well
known formula (see \cite[p.\,38]{Whitney}) expresses the exterior
product of vectors $\mv x_i=\sum_{j=0}^nx_{i,j}\ \mv e_j\in \V$ \ \
$(1\le i\le p)$ in terms of the basis (\ref{e:025})\,:
\begin{equation}\label{e:026}
\We_{i=1}^p\ \mv x_i=\sum_{I=\{i_1<\dots<i_p\}\in C(n+1,p)}
\det\Big(x_{j,i_k}\Big)_{1\le j,k\le p}\ \mv e_I.
\end{equation}
Recall that the exterior product is
\emph{alternating}, that is $\mv u\we\mv v=-\mv v\we\mv u$ so that $\mv v\we\mv v=0$. Further,
let $\mv u\dt\mv v$ denote the standard inner product of $\mv u,\mv
v\in \V$. Then, there is a uniquely defined \emph{inner product}\/
on $\We^p(\V)$ such that
\begin{equation}\label{e:027}
(\mv v_1\we\dots\we\mv v_p)\dt(\mv u_1\we\dots\we\mv u_p)= \det\big(
\mv v_i\dt\mv u_j \big)_{1\le i,j\le p}
\end{equation}
for any $\mv v_1,\dots,\mv v_p,\mv u_1,\dots,\mv u_p\in \V$.
Furthermore, if $\mv e_0,\dots,\mv e_n$ is an orthonormal basis then
so is (\ref{e:025}). Often (\ref{e:027}) is referred to as the
\emph{Laplace identity} \cite[p.105]{Schmidt-1980}. The Euclidean
norm on $\We^p(\V)$ induced by (\ref{e:027}) will be denoted by
$|\cdot|$. By (13) in \cite[p.\,49]{Whitney},
\begin{equation}\label{e:028}
|\mv u\we\mv v|\le|\mv u|\,|\mv v|\qquad\text{if $\mv u$ or $\mv v$
is decomposable.}
\end{equation}
Recall that a multivector $\mv u$ is \emph{decomposable}\/ if $\mv
u=\mv u_1\we\dots\we\mv u_p$ for some $\mv u_1,\dots,\mv u_p\in\V$.
Finally, given $\vv x,\vv y\in\R^n$,
$$
    \pdist(\vv x,\vv y)\ = \ \frac{|\mv x\we\mv y|}{|\mv x|\, |\mv
    y|}
$$
is called the \emph{projective distance}\/ between $\vv x$ and $\vv
y$. Obviously $\pdist(\vv x,\vv y)$ is well defined. It is known
that $\pdist(\vv x,\vv y)=\sin\varphi(\mv x,\mv y)$, where
$\varphi(\mv x,\mv y)$ denotes the acute angle between $\mv x$ and
$\mv y$ -- see (\ref{e:037}) below. In particular, this angular
property of $\pdist$ implies that $\pdist(\vv x,\vv y)$ is a metric.
Furthermore, $\pdist$ is locally comparable to the euclidean norm since
\begin{equation}\label{e:029}
\pdist(\vv x,\vv y) \ \le\ |\vv x-\vv y| \ \le \ \sqrt{1+|\vv
x|^2}\,\sqrt{1+|\vv y|^2} \, \pdist(\vv x,\vv y)
\end{equation}
for all $\vv x,\vv y\in\R^n$. To see that (\ref{e:029}) is true take
$\mv x=(1,x_1,\dots,x_n)$ and $\mv y=(1,y_1,\dots,y_n)$. Then the
l.h.s.\! of (\ref{e:029}) (\emph{l.h.s.\!\! means left hand side})
is proved as follows
$$
 \pdist(\vv x,\vv y)= \frac{|\mv x\we\mv y|}{|\mv x|\,|\mv y|}=\frac{|(\mv x-\mv y)\we\mv
y|}{|\mv x|\,|\mv y|}\stackrel{\eqref{e:028}}{\le}\frac{|\mv x-\mv
y|\,|\mv y|}{|\mv x|\,|\mv y|}
    =  \displaystyle\frac{|\mv x-\mv y|}{|\mv x|}\le
|\mv x-\mv y|=|\vv x-\vv y|.
$$
On the other hand, $|\vv x-\vv y|\stackrel{}{\le} |\mv x\we\mv y|=
\sqrt{1+|\vv x|^2}\,\sqrt{1+|\vv y|^2}\,\pdist(\vv x,\vv y)$, where
the first inequality is a consequence of (\ref{e:026}).

\mysubsection{Interior product and Hodge duality}

In what follows ``\ $\dt$\ '' will denote the \emph{interior
product} of multivectors. For $\mv u\in\We^p(\V)$ and $\mv
v\in\We^q(\V)$ the latter is defined as follows. Assume that $p\ge
q$ and consider the two linear functions on $\We^{p-q}(\V)$ given by
$$
\mv x\mapsto\mv u\dt(\mv v\we\mv x)\qquad\text{and}\qquad \mv
x\mapsto(\mv x\we\mv v)\dt\mv u.
$$
 Since $\We^{p-q}(\V)$ is Euclidean there are unique
$(p-q)$-vectors, which will be denoted by $\mv u\dt\mv v$ and $\mv
v\dt\mv u$, such that $
 (\mv u\dt\mv v)\dt\mv x=\mv u\dt(\mv v\we\mv x)\quad\text{and}\quad
 \mv x\dt(\mv v\dt\mv u)=(\mv x\we\mv v)\dt\mv u
$ for all $\mv x\in\We^{p-q}(\V)$. The multivectors $\mv u\dt\mv v$
and $\mv v\dt\mv u$ are called the \emph{interior products}\/ of
$\mv u$ and $\mv v$, and $\mv v$ and $\mv u$ respectively. It is
easily seen that $\mv v\dt\mv u=(-1)^{q(p-q)}\mv u\dt\mv v$ and that
in the case $p=q$ the interior product is simply the inner product
(\ref{e:027}). The definition of interior product readily implies
that
\begin{equation}\label{e:030}
\mv a\dt(\mv b\we\mv c)=(\mv a\dt\mv b)\dt\mv
c\qquad\text{and}\qquad (\mv c\we\mv b)\dt \mv a=\mv c\dt(\mv
b\dt\mv a)
\end{equation}
if $\mv a\in\We^p(\V)$, $\mv b\in\We^q(\V)$, $\mv c\in\We^r(\V)$
with $p\ge q+r$ -- see (5)$+$(6) in \cite[p.\,43]{Whitney}.

Let $\vv e_0,\dots,\vv e_n$ be the standard basis of $\V$ and $\mv
i=\vv e_0\we\vv e_1\we\ldots\we\vv e_n\in\We^{n+1}(\V)$. By ``\
${}^\perp$\ '' we will denote the Hodge star operator which is
defined by
\begin{equation}\label{e:031}
 \mv u^\perp:=\mv i\dt\mv u\,.
\end{equation}
Note that the multivector $\mv u\in\We^p(\V)$ is decomposable if and
only if $\mv u^\perp\in\We^{n+1-p}(\V)$ is decomposable -- see
Lemma~11A in \cite[p.\,48]{Whitney}. The map (\ref{e:031}) is
obviously linear. Also
\begin{equation}\label{e:032}
 (\mv v^\perp)^\perp= (-1)^{(n+1-p)p}\mv v\qquad\text{for any $\mv v\in\We^p(\V)$. }
\end{equation}
The latter, know as the Hodge duality, follows from (2) in
\cite[p.\,49]{Whitney} but can also be easily verified for basis
vectors and then extended by linearity. Obviously  $\mv v\mapsto\mv
v^\perp$ is a one-to-one correspondence between $\We^p(\V)$ and
$\We^{n+1-p}(\V)$. Also, an easy consequence of (\ref{e:030}) and
(\ref{e:032}) is that the Hodge operator is an isometry, that is
$|\mv v^\perp| = |\mv v|$ for any $\mv v\in\We^p(\V)$. Also the
Hodge operator conveniently relates the interior and exterior
products. Indeed, let $\mv u\in\We^p(\V)$ and $\mv v\in\We^q(\V)$.
Then using (\ref{e:030}) readily gives
\begin{equation}\label{e:033}
\mv v^\perp\dt\mv u=(\mv v\we\mv u)^\perp \qquad \text{if}\qquad
p+q\le n+1.
\end{equation}
Since the Hodge operator is an isometry, this relation implies that
\begin{equation}\label{e:034}
    |\mv v^\perp\dt\mv u|=|\mv v\we\mv u|\qquad\text{if}\qquad p+q\le n+1.
\end{equation}

\mysubsection{Relations between multivectors and subspaces of
$\V$}\label{mvss}

Throughout, $\cV(\mv v_1,\dots,\mv v_r)$ denotes the vector space
spanned by vectors $\mv v_1,\dots,\mv v_r$. Also, given a
multivector $\mv w\in\We(\V)$, let $\cV(\mv w)$ be the linear
subspace of $\V$ given by
$$
\cV(\mv w):=\{\mv x\in\V:\mv w\we\mv x=\mv 0\}.
$$

\begin{lemma}\label{l:02}
If $\mv u_1,\dots,\mv u_p\in \V$ are linearly independent, then
$\cV(\mv u_1\we\dots\we\mv u_p)=\cV(\mv u_1,\dots,\mv u_p)$.
Furthermore if $\mv u,\mv v\in\We^p(\V)$ are non-zero decomposable
multivectors, then \ $\cV(\mv u)=\cV(\mv v)$ \ \
$\Longleftrightarrow$ \ \ $\mv u=\theta\mv v$ for some
$\theta\not=0$.
\end{lemma}
For details see Lemma~6B and Lemma~6C in
\cite[pp.\,104--105]{Schmidt-1980}. Lemma~\ref{l:02} gives a
one-to-one correspondence between non-zero decomposable $p$-vectors
taken up to a constant multiple and linear subspaces in $\V$ of
dimension $p$. The latter is known as a \emph{Grassmann manifold}\/
and will be denoted by $\Gr_p(\V)$. Thus $\Gr_p(\V)$ is embedded
into $\P(\We^p(\V))$ and so is equipped with a natural topology
induced from $\P(\We^p(\V))$ with respect to which it is obviously
compact. Naturally, through the above correspondence the elements of
$\Gr_p(\V)$ can be thought of as unit decomposable $p$-vectors taken
up to sign.

The following lemma gives a convenient way of expressing orthogonal
subspaces via the Hodge operator and justifies the notation
for the operator that we use within this paper. In what follows
$W^\perp$ denotes the linear subspace of $\V$ orthogonal to
$W\subset\V$.

\begin{lemma}\label{l:03}
Let $\mv u\in\We^p(\V)$ be a non-zero decomposable multivector. Then
\begin{equation}\label{e:035}
\cV(\mv u^\perp)=\cV(\mv u)^\perp=\{\mv v\in\V:\mv u\dt\mv v=0\}.
\end{equation}
\end{lemma}

\noindent\textit{Proof.} Take any orthogonal basis $\mv
e_1,\dots,\mv e_p$ of $\cV(\mv u)$ such that $\mv u=\mv
e_1\we\dots\we\mv e_p$. This is possible in view of
Lemma~\ref{l:02}. If $\mv v\in\V$ is orthogonal to $\cV(\mv u)$
then, using (\ref{e:027}) it is easy to see that $\mv u\dt(\mv
v\we\mv x)=0$ for any decomposable $\mv x\in\We^{p-1}(\V)$. On the
other hand, if $\mv v\in\V$ is not orthogonal to $\cV(\mv u)$, say
$\mv e_1\dt\mv v\not=0$, then, by (\ref{e:027}), $\mv u\dt(\mv
v\we\mv e_2\we\dots\we\mv e_p)=\mv e_1\dt\mv v\not=0$. The upshot is
that $\mv u\dt(\mv v\we\mv x)$ vanishes identically for all $\mv
x\in\We^{p-1}(\V)$ if and only if $\mv v\in\cV(\mv u)^\perp$. By the
definition of interior product, this precisely means that $\mv
u\dt\mv v=0$ if and only if $\mv v\in\cV(\mv u)^\perp$. The latter
establishes the r.h.s.\! of (\ref{e:035}). Finally, by
(\ref{e:034}), $\mv u\dt\mv v=0$ if and only if $\mv u^\perp\we\mv
v=0$. The latter implies the l.h.s.\! of (\ref{e:035}). \proofend

\begin{lemma}\label{l:04}
Let $\mv u\in\We^p(\V)$ and $\mv v\in\We^q(\V)$ be decomposable.
Then $\cV(\mv u)\cap\cV(\mv v)=\varnothing$ if and only if $\mv
u\we\mv v\not=0$. Consequently, if $\mv u\we\mv v\not=0$ then
$\cV(\mv u)\oplus\cV(\mv v)=\cV(\mv u\we\mv v)$. Also if $p\ge q$
and $\mv u\dt\mv v\not=\vv0$ then $\cV(\mv u\dt\mv v)=\cV(\mv
u)\cap\cV(\mv v^\perp)$.
\end{lemma}

\noindent\textit{Proof.} The condition $\cV(\mv u)\cap\cV(\mv
v)=\varnothing$ means that the sum $\cV(\mv u)+\cV(\mv v)$ is
direct, which is equivalent to $\mv u\we\mv v\not=0$. The equality
$\cV(\mv u)\oplus\cV(\mv v)=\cV(\mv u\we\mv v)$ is then a
consequence of Lemma~\ref{l:02}. Finally, by (\ref{e:033}), $\mv
u\dt\mv v=\pm(\mv u^\perp\we\mv v)^\perp$. Then, by
Lemmas~\ref{l:03} and \ref{l:04}, $ \cV(\mv u\dt\mv v)=\cV(\mv
u^\perp\we\mv v)^\perp=(\cV(\mv u^\perp)\oplus\cV(\mv
v))^\perp=\cV(\mv u^\perp)^\perp\cap\cV(\mv v)^\perp=\cV(\mv
u)\cap\cV(\mv v^\perp)$. \proofend

The following lemma is easily established using the Laplace identity
(\ref{e:027}).

\begin{lemma}\label{l:06}
Let $\mv u\in\We^p(\V)$ and $\mv v\in\We^q(\V)$ be decomposable and
$p+q\le n+1$. If\/ $\cV(\mv u)\perp\cV(\mv v)$ then $|\mv u\we\mv
v|=|\mv u|\,|\mv v|$.
\end{lemma}

\mysubsection{Multivectors and projections}\label{alm}

There are various relations between exterior/interior product and
projections of vectors in $\V$ onto subspaces. The properties we are
particularly interested in are summarized as

\begin{lemma}\label{l:07}
Let $\mv u\in \V$, $\mv v\in\We^p(\V)$ with $1\le p\le n$ be
decomposable and let $\pi$ denote the orthogonal projection from
$\V$ onto $\cV(\mv v)$. Then
\begin{equation}\label{e:036}
|\mv v\we\mv u|=|\mv v|\cdot|\mv u-\pi\mv u|\qquad\text{and}\qquad
|\mv v\dt\mv u|=|\mv v|\cdot|\pi\mv u|.
\end{equation}
Furthermore, $|\mv v|^{2}\pi\mv u=\pm\,\mv v\dt(\mv v\dt\mv u)$,
where the sign is either $+$ or $-$.
\end{lemma}

\noindent\textit{Proof.} Fix an orthogonal basis $\mv v_1,\dots,\mv
v_p$ of $\cV(\mv v)$ such that $\mv v=\mv v_1\we\dots\we\mv v_p$.
Let $\mv u'=\mv u-\pi\mv u$. Obviously $\mv v_1,\dots,\mv v_p,\mv
u'$ is an orthogonal system. Also, since $\pi\mv u\in\cV(\mv v)$, by
Lemma~\ref{l:02}, $\mv v\we\pi\mv u=\vv0$. Therefore,
$\mv v\we\mv u=\mv v\we\mv u'$. Now applying (\ref{e:027}) gives
$$
 |\mv v\we\mv u|^2  =  |\mv v_1\we\dots\we\mv v_p\we\mv u'|^2
 \stackrel{\eqref{e:027}}{=} |\mv v_1|^2\dots|\mv v_p|^2\,|\mv
u'|^2 \stackrel{\eqref{e:027}}{=}|\mv v|^2\,|\mv u-\pi\mv u|^2.
$$
This establishes the l.h.s.\! of (\ref{e:036}). Further, notice that
$\mv u-\pi\mv u$ is the orthogonal projection of $\mv u$ onto
$\cV(\mv v^\perp)=\cV(\mv v)^\perp$. Therefore, the r.h.s.\! of
(\ref{e:036}) follows on applying (\ref{e:034}) to the l.h.s.\! of
(\ref{e:036}), when $\mv v$ is replaced by $\mv v^\perp$. The final
identity of the lemma is very well known and easy when $p=1$. We
consider $p\ge2$. First, notice that $\mv u\we\pi\mv u=\mv u\we(\mv
u-\mv u')=-\mv u\we\mv u'$ and that $\mv v\dt\mv u'=0$ -- see
Lemma~\ref{l:03}. Therefore, $ (\mv v\dt\mv u)\dt\pi\mv u
\stackrel{\eqref{e:030}}{=}  \mv v\dt(\mv u\we\pi\mv u)\ =\ -\mv
v\dt(\mv u\we\mv u')
 = \mv v\dt(\mv u'\we\mv u)\ \stackrel{\eqref{e:030}}{=}\ (\mv v\dt\mv
u')\dt\mv u\ =\ 0. $ Hence, by Lemma~\ref{l:03}, $\pi\mv
u\perp\cV(\mv v\dt\mv u)$. Also, since $\pi$ is the projection onto
$\cV(\mv v)$, we have that $\pi\mv u\perp\cV(\mv v)^\perp=\cV(\mv
v^\perp)$. Therefore, $\pi\mv u\perp\cV(\mv v\dt\mv u)+\cV(\mv
v^\perp)$. By Lemma~\ref{l:04}, the space $\cV(\mv v\dt\mv u)$ is a
subspace of $\cV(\mv v)$ and so is orthogonal to $\cV(\mv v^\perp)$.
Then, the sum $\cV(\mv v\dt\mv u)+\cV(\mv v^\perp)$ is direct and,
by Lemma~\ref{l:04}, it is equal to $\cV(\mv v^\perp\we(\mv v\dt\mv
u))$. The latter space is readily seen to have codimension 1.
Theretofore, the relation $\pi\mv u\perp\cV(\mv v\dt\mv u)+\cV(\mv
v^\perp)$ implies that $\pi\mv u\|\big(\mv v^\perp\we(\mv v\dt\mv
u)\big)^\perp\stackrel{\eqref{e:033}}{=}\pm\mv v\dt(\mv v\dt\mv u).$
Finally, since the Hodge operator is an isometry,
$$
|\mv v|^2\cdot|\pi\mv u|= |\mv v^\perp|\cdot|\mv v|\cdot|\pi\mv u|
\stackrel{\eqref{e:036}}{=} |\mv v^\perp|\cdot|\mv v\dt\mv u|
\stackrel{\text{{\rm Lemma~\ref{l:06}}}}{=} |\mv v^\perp\we(\mv v\dt\mv u)|
\stackrel{\eqref{e:034}}{=}|\mv v\dt(\mv v\dt\mv u)|
$$
and the identity $|\mv v|^{2}\pi\mv u=\pm\,\mv v\dt(\mv v\dt\mv u)$
now readily follows. \proofend

\bigskip

Given two lines $\ell_1$ and $\ell_2$ in $\V$ through the origin,
let $\varphi(\ell_1,\ell_2)$ denote the acute angle between $\ell_1$
and $\ell_2$. Further, given a linear subspace $L$ of $\V$ of
dimension $p$ and a line $\ell$ through the origin, the angle
$\varphi(\ell,L)$ between $L$ and $\ell$ is defined to be $
\inf_{\ell'\in L}\varphi(\ell,\ell')\,, $ where the infimum is taken
over over lines $\ell'$ in $L$ through the origin. It is well known
that $\varphi(\ell,L)$ is the angle between $\ell$ and the
orthogonal projection of $\ell$ onto $L$. Thus, if $\mv u$ is a
directional vector of $\ell$ and $\pi$ denotes the orthogonal
projection onto $L$ then $\sin\varphi(\ell,L)=|\mv u|^{-1}|\mv
u-\pi\vv u|$. Further, if $\mv v\in\We^p(\V)$ is a Grassmann
representative of $L$, that is $L=\cV(\mv v)$, then, by
Lemma~\ref{l:07},
\begin{equation}\label{e:037}
\sin\varphi(\ell,L)=\frac{|\mv v\we\mv u|}{|\mv v|\,|\mv
u|}\stackrel{\eqref{e:034}}{=}\frac{|\mv v^\perp\dt\mv u|}{|\mv
v|\,|\mv u|}.
\end{equation}

The following lemma is a consequence of the fact that the angle
between a line $\ell$ and a plane $L_1$ is not bigger than the angle
between this line $\ell$ and any other plane $L_2\subset L_1$.

\begin{lemma}\label{l:08}
Let  $\mv v\in\We^p(\V)$ be a non-zero decomposable multivector and
$\mv u\in\V$. Then for any non-zero $\mv w\in\cV(\mv v)$
 $$\frac{|\mv w\dt\mv u|}{|\mv w|}\le \frac{|\mv v\dt\mv u|}{|\mv
v|}.$$
\end{lemma}

\noindent\textit{Proof.} In view of (\ref{e:034})
\begin{equation}\label{e:038}
\dfrac{|\mv w\dt\mv u|}{|\mv w|}\le \dfrac{|\mv v\dt\mv u|}{|\mv
v|}\qquad\Longleftrightarrow\qquad
 \dfrac{|\mv v^\perp\we\mv u|}{|\mv v|\,|\mv u|}\ \ge \
\dfrac{|\mv w^\perp\we\mv u|}{|\mv w|\,|\mv u|}.
\end{equation}
Obviously $L_2:=\cV(\mv v^\perp)\subset L_1:=\cV(\mv w^\perp)$. Let
$\ell:=\cV(\mv u)$. Therefore, by (\ref{e:037}), the l.h.s. of
(\ref{e:038}) is equivalent to $\sin\varphi(\ell,L_2)\ge
\sin\varphi(\ell,L_1)$. The latter is obvious in view of the fact
that $L_2\subset L_1$. The proof is thus complete. \proofend

\section{Detecting rational points near a manifold}\label{detecting}

In this section we describe the mechanism for investigating the
distribution of rational points near manifolds.

\mysubsection{Local geometry near a manifold}

Let $\cM$ be a $C^{(2)}$ manifold of the Monge form (\ref{e:002}).
For $\vv x=(x_1,\dots,x_\ddd)\in U$ let $\vv y=\vv y(\vv x)$ be the point
$(\vv x,\vv f(\vv x))\in\cM$. We
will use the lifting of $\cM$ into $\V$ given by
\begin{equation}\label{e:039}
    \mv y(\vv x)=(1,\vv y(\vv x))=(1,\vv x,\vv f(\vv x))
\end{equation}
which represents the projective embedding of $\vv y(\vv x)$.
Further, consider the following maps:
  \begin{equation}\label{e:040}
  \mv g:U\to\We^{\ccc }(\V)\ :\ \vv x\mapsto
  \big(\mv y(\vv x)\we\pd{\mv y(\vv x)}{1}\we\dots\we\pd{\mv y(\vv
  x)}{\ddd}\big)^\perp
  \end{equation}
and
\begin{equation}\label{e:041}
  \mv u:U\to\We^{\ddd}(\R^{n+1})\ :\ \vv x\mapsto \big(\mv y(\vv x)\we\mv g(\vv
  x)\big)^\perp,\hspace*{17.4ex}
\end{equation}
where $\partial_i:=\partial/\partial x_i$. Since $\vv y(\vv x)$ is
of the Monge form the vectors $\mv y(\vv x),\pd{\mv y(\vv
x)}{1},\dots,\pd{\mv y(\vv x)}{\ddd}$ are linearly independent, thus
giving $\mv g(\vv x)\not=\vv 0$. Also, by Lemma~\ref{l:03}, $\mv
y(\vv x)\perp\cV(\mv g(\vv x))$. Therefore, $\mv y(\vv x)\we\mv
g(\vv x)\not=0$ further implying $\mv u(\vv x)\not=\vv 0$.

\emph{Convention.} In order to simplify notation, we will write $\gx$, $\ux$ and $\yx$ for $\mv
g(\vv x)$, $\mv u(\vv x)$ and $\mv y(\vv x)$ respectively and drop
the subscript $\vv x$ whenever there is no risk of confusion. It is
useful to keep in mind the following geometric nature of $\mv g$ and
$\mv u$. The homogeneous equations $\mv g\dt\mv z=0$ and $\mv
u\dt\mv z=0$ with respect to $\mv z=(z_0,z_1,\dots,z_n)$ define the
tangent and transversal planes to $\cM$ respectively. Furthermore,
$|\mv g\dt\mv r|\,|\mv g|^{-1}|\mv r|^{-1}$ (\emph{resp}. $|\mv
u\dt\mv r|\,|\mv u|^{-1}|\mv r|^{-1}$) is the projective distance of
$\mv r$ from the tangent (\emph{resp}. transversal) plane -- see
(\ref{e:037}).

\begin{lemma}\label{l:09}
For every $\vv x\in U$ we have that $\V=\cV(\mv g)\oplus\cV(\mv
u)\oplus\cV(\mv y)$ is a decomposition of\/ $\V$ into pairwise
orthogonal subspaces.
\end{lemma}

\noindent\textit{Proof. } Recall the convention that $\mv g=\gx$, $\mv u=\ux$ and $\mv y=\yx$. Fix an $\vv x\in U$. Let $\mv t:=\pd{\mv
y(\vv x)}{1}\we\dots\we\pd{\mv y(\vv x)}{\ddd}$. Then, by
(\ref{e:040}), $\mv g=(\mv y\we\mv t)^\perp$. Then, using
Lemmas~\ref{l:02} and \ref{l:03}, we get that $\cV(\mv g)=\cV((\mv
y\we\mv t)^\perp)=\cV(\mv y\we\mv t)^\perp\subset \cV(\mv y)^\perp$.
It follows that $\cV(\mv g)\perp\cV(\mv y)$. It is similarly
established that $\cV(\mv u)\perp\cV(\mv y)$ and $\cV(\mv
g)\perp\cV(\mv u)$. Thus, the subspaces $\cV(\mv g)$, $\cV(\mv u)$
and $\cV(\mv y)$ are pairwise orthogonal and so their sum is direct.
Moreover, using Lemma~\ref{l:02} one readily finds the dimension of
each of the subspaces, resulting in $ \dim\cV(\mv g)\oplus\cV(\mv
u)\oplus\cV(\mv y) = n+1$. Therefore, $\V=\cV(\mv g)\oplus\cV(\mv
u)\oplus\cV(\mv y)$. \proofend

\medskip

Lemma~\ref{l:09} provides a natural choice for local coordinates
akin to the Frenet frame. The following Lemma~\ref{l:10} estimates the
projective distance of a point $\vv r\in\R^n$ from $\vv y\in \cM$ in
terms of the projective distance of $\vv r$ from the tangent and
transversal planes.

\begin{lemma}\label{l:10}
For any $\mv r\in\R^{n+1}$ and any $\vv x\in U$
\begin{equation}\label{e:042}
\frac{|\mv y\we\mv r|}{|\mv y|}\le\frac{|\mv g\dt\mv r|}{|\mv
g|}+\frac{|\mv u\dt\mv r|}{|\mv u|}.
\end{equation}
\end{lemma}

\noindent\textit{Proof}. Let $\mv r_{\mv g}$, $\mv r_{\mv u}$ and
$\rf$ be the orthogonal projections of $\mv r$ onto $\cV(\mv g)$,
$\cV(\mv u)$ and $\cV(\mv y)$ respectively. Then, by
Lemma~\ref{l:09}, $\mv r=\rg+\ru+\rf$ and therefore
$\mv r-\rf=\rg+\ru$. By Lemma~\ref{l:07},
\begin{equation}\label{e:043}
|\mv y\we\mv r|\cdot|\mv y|^{-1}=|\mv r-\rf|\le|\mv r_{\mv g}|+|\mv
r_{\mv u}|.
\end{equation}
Again, by Lemma~\ref{l:07}, $|\mv g\dt\mv r|=|\mv g|\cdot|\rg|$ and
$|\mv u\dt\mv r|=|\mv u|\cdot|\ru|$. Substituting $|\rg|$ and
$|\ru|$ from the latter equalities into (\ref{e:043}) gives
(\ref{e:042}).\proofend

\medskip

Lemma~\ref{l:10} is in general sharp as (\ref{e:042}) can be
reversed with some positive constant. Nevertheless, the distance of
$\mv r$ from $\cM$ rather than from a particular point $\vv y$ on
$\cM$ can be estimated in a more efficient way. This relies on the
fact that the tangent plane deviates from a $C^{(2)}$ manifold with
a quadratic error. A similar idea is explored by Elkies
\cite{Elkies-2000} in his algorithm for computing rational points
near manifolds. Before we state the next result, recall that
given a ball $B=B(\vv x,r)$ and $\lambda>0$, $\lambda B:=B(\vv
x,\lambda r)$ and $\overline B$ is the closure of $B$.

\begin{lemma}\label{l:11}
Let $\cM$ be a $C^{(2)}$ manifold of the form $(\ref{e:002})$ and
$B_0$ be a ball of radius $r_{B_0}<\infty$ such that $2\overline
B_0\subset U$. Then there is a constant $C>1$ depending on $B_0$
only satisfying the following property. For any $\mv r\in\V$ and
$\vv x\in B_0$ such that
\begin{equation}\label{e:044}
    \frac{|\gx\dt\mv r|}{|\gx|\,|\mv r|}<\delta\qquad\text{and}\qquad
    \frac{|\ux\dt\mv r|}{|\ux|\,|\mv r|}<\ve
\end{equation}
for some positive $\delta$ and $\ve$ satisfying
\begin{equation}\label{e:045}
    \ve^2\le \delta\le \ve\le\ve_0:=\frac{\min\{1,r_{B_0}\}}{2\ddd(n+1)(C+1)^2}
\end{equation}
there is a point $\vv x'\in 2B_0$ such that
\begin{equation}\label{e:046}
\frac{|\yxd\we\mv r|}{|\yxd|\, |\mv r|}\le
K\,\delta\,,\qquad\text{where $K=14(n+1)^3(C+1)^5\ddd^2$.}
\end{equation}
\end{lemma}

\noindent\textit{Proof of Lemma~\ref{l:11}}. Without loss of
generality we will assume that $|\mv r|=1$. Since $2\overline
B_0\subset U$, there is a constant $C>1$ such that
\begin{equation}\label{e:047}
   2 \overline B_0\subset[-C,C]^{\ddd}
\end{equation}
and
\begin{equation}\label{e:048}
    \sup_{\vv x\in 2\overline B_0} \max\Big\{|f_l(\vv x)|, \max_{1\le i\le \ddd}\ |\partial_{i}f_l(\vv x)|,
    \max_{1\le i,j\le \ddd} |\partial_{i}\partial_jf_l(\vv x)|\Big\} \le  C
\end{equation}
for $1\le l\le \ccc$, where $\partial_i$ means differentiating by
$x_i$ and the functions $f_l$ arise from (\ref{e:002}).

\medskip

\noindent\textbf{Step 1.} At this step we express $\mv r$ as a
linear combination of $\mv y$, $\partial_1\mv
y$,\dots,$\partial_{\ddd}\mv y$ plus an error term. Let $\rg$, $\ru$
and $\rf$ be the orthogonal projections of $\mv r$ onto $\cV(\mv
g)$, $\cV(\mv u)$ and $\cV(\mv y)$ respectively. By Lemma~\ref{l:07}
and the assumption $|\mv r|=1$, inequalities (\ref{e:044}) imply
that
\begin{equation}\label{e:049}
    |\mv r_{\mv g}|<\delta\qquad\text{and}\qquad |\mv r_{\mv
    u}|<\ve.
\end{equation}
Also, by Lemma~\ref{l:10}, inequalities (\ref{e:044}) imply that
$|\mv y|^{-1}|\mv y\we\mv r|<\delta+\ve$. By (\ref{e:027}), we have
the identity $|\mv y\we\mv r|^2=|\mv y|^2|\mv r|^2-|\mv y\dt\mv
r|^2$. Since $|\mv r|=1$, the latter implies
$$
0\le 1-\frac{|\mv y\dt\mv r|}{|\mv y|}\le 1-\left(\frac{|\mv y\dt\mv
r|}{|\mv y|}\right)^2\ = \ \left(\frac{|\mv y\we\mv r|}{|\mv
y|}\right)^2\le
(\delta+\ve)^2\stackrel{\eqref{e:045}}{\le}4\delta.
$$
The latter inequality together with the fact that $|\mv y|^{-1}|\mv
y\dt\mv r|=|\mv r_{\mv y}|$ implied by Lemma~\ref{l:07}, shows that
for some $\eta\in\{-1,1\}$
\begin{equation}\label{e:050}
\mv r_{\mv y}=\eta|\mv y|^{-1}\mv y+\mv
w_0\qquad\text{with}\qquad|\mv w_0|\le 4\delta.
\end{equation}
By (\ref{e:040}) and Lemma~\ref{l:02}, we see that the vectors $\mv
y=\mv y(\vv x)$, $\pd{\mv y(\vv x)}{1},\dots,\pd{\mv y(\vv
x)}{\ddd}$ form a basis of $\cV(\mv g^\perp)$. By Lemmas~\ref{l:02}
and \ref{l:09}, we have that $\cV(\mv u)\subset\cV(\mv g^\perp)$.
Therefore, since $\mv r_{\mv u}\in\cV(\mv u)$, there are real
numbers $\lambda_0,\dots,\lambda_{\ddd}$ such that
\begin{equation}\label{e:051}
 \mv r_{\mv u}=\lambda_0\mv y(\vv
 x)+\sum_{i=1}^{\ddd}\lambda_i\pd{\mv y(\vv x)}{i}.
\end{equation}
Since $\mv y$ is of the Monge form,
\begin{equation}\label{e:052}
 \mv r_{\mv u}=(\lambda_0,\lambda_1+x_1\lambda_0,\dots,\lambda_{\ddd}+x_{\ddd}\lambda_0,*,\dots,*),
\end{equation}
where $\,\,*\,\,$ stands for a real number. By (\ref{e:047}),
(\ref{e:052}) and the r.h.s.\! of (\ref{e:049}),
\begin{equation}\label{e:053}
 |\lambda_0|<\ve\qquad\text{and}\qquad |\lambda_i|<(C+1)\ve\qquad(1\le i\le \ddd).
\end{equation}
On plugging the expressions for $\mv r_{\mv y}$ and $\mv r_{\mv u}$
given by (\ref{e:050}) and (\ref{e:051}) into the identity $\mv
r=\mv r_{\mv g}+\mv r_{\mv u}+\mv r_{\mv y}$ and applying the
l.h.s.\! of (\ref{e:049}) we get
\begin{equation}\label{e:054}
\mv r=\lambda_0^*\mv y(\vv
 x)+\sum_{i=1}^{\ddd}\lambda_i\pd{\mv y(\vv x)}{i}+\mv
 w_1\qquad\text{with $|\mv w_1|\le 5\delta$, }
\end{equation}
where $\lambda_0^*=\eta|\mv y(\vv x)|^{-1}+\lambda_0$.

\bigskip

\noindent\textbf{Step 2.} At this step we define the point $\vv x'$.
By (\ref{e:047}) and (\ref{e:048}),  $|\mv y(\vv
x)|^{-1}\ge (n+1)^{-1}C^{-1}$. On the other hand, by (\ref{e:045})
and (\ref{e:053}),
$|\lambda_0|\le\tfrac12\,(n+1)^{-1}C^{-1}$. Therefore,
$|\lambda_0^*|\ge \tfrac12\,(n+1)^{-1}C^{-1}$ or equivalently
\begin{equation}\label{e:055}
|\lambda_0^*|^{-1}\le 2\,(n+1)C .
\end{equation}
Further, define $\lambda_i^*=\lambda_i/\lambda_0^*$ for
$i=1,\dots,\ddd$. Inequalities (\ref{e:053}) and (\ref{e:055})
imply that
\begin{equation}\label{e:056}
|\lambda_i^*|\le2\ve(n+1)(C+1)^2\qquad(1\le i\le \ddd).
\end{equation}
Dividing (\ref{e:054}) by $\lambda_0^*$ and applying (\ref{e:055})
to estimate the remainder term gives
\begin{equation}\label{e:057}
    \lambda_0^*{}^{-1}\mv r=\mv y(\vv
 x)+\sum_{i=1}^{\ddd}\lambda_i^*\pd{\mv y(\vv x)}{i}+\mv
 w_2\qquad\text{with \ \ $|\mv w_2|\le 10\delta(n+1)C$. }
\end{equation}
Now define $\vv x'=\vv x+\bm\lambda^*$, where
$\bm\lambda^*=(\lambda_1^*,\dots,\lambda_{\ddd}^*)$. Conditions
(\ref{e:045}) and (\ref{e:056}) ensure that $|\bm\lambda^*|\le
r_{B_0}$. Therefore, since $\vv x\in B_0$,  $\vv x'\in
2B_0$.

\bigskip

\noindent\textbf{Step 3.} At this step we verify (\ref{e:046}). By
(\ref{e:047}), (\ref{e:048}), (\ref{e:056}) and Taylor's
formula, we get
\begin{equation}\label{e:058}
\big|\mv y(\vv x')-\mv y(\vv x)-\sum_{i=1}^{\ddd}\lambda_i^*\pd{\mv y(\vv
x)}{i}\big|\le 4\ve^2(n+1)^3(C+1)^5\ddd^2 .
\end{equation}
Further, using (\ref{e:045}), (\ref{e:057}) and (\ref{e:058})
we get
\begin{equation}\label{e:059}
    |\yxd-\lambda_0^*{}^{-1}\mv r|\le\delta\Big(10(n+1)C+4(n+1)^3(C+1)^5\ddd^2\Big)\le K\delta.
\end{equation}
From (\ref{e:039}), $|\yxd|\ge1$. Therefore, using
$|\mv r|=1$, we obtain
$$
\dfrac{|\yxd\we\mv r|}{|\yxd|\, |\mv r|} \le |\yxd\we\mv
r|=|(\yxd-\lambda_0^*{}^{-1}\mv r)\we\mv r|
\stackrel{\eqref{e:028}}{\le}  |(\yxd-\lambda_0^*{}^{-1}\mv
r)|\cdot|\mv r|\stackrel{\eqref{e:059}}{\le} K\delta.
$$
This establishes (\ref{e:046}) and thus completes the proof of
Lemma~\ref{l:11}.\proofend

\mysubsection{Good ``cells'' near a manifold}\label{GCT}

Let $\psi_*$, $Q_*$ and $\kappa$ be positive parameters. In
practice, $Q_*$ and $\psi_*$ will be proportional to $Q$ and $\psi$
respectively. Further, for every $\vv x\in U$ consider the system
\begin{equation}\label{e:060}
         \frac{|\gx\dt\mv r|}{|\gx|} < \psi_*\,,\qquad
         \frac{|\ux\dt\mv r|}{|\ux|} < (\psi_*^{\ccc } Q_*)^{-\frac{1}{\ddd}}\,,\qquad
         \frac{|\yx\dt\mv r|}{|\yx|}  \le  \kappa\, Q_*,
\end{equation}
where $\mv r\in\V$.
Obviously the set of $\mv r\in\V$ satisfying (\ref{e:060}) is a
convex body symmetric about the origin. Then as a consequence of
Minkowski's theorem for convex bodies one has

\begin{lemma}\label{l:12}
Let $v_d$ denote the volume of a ball of diameter 1 in $\R^d$ and
$\kappa_0:=(v_dv_m)^{-1}$. Then, for any $\kappa\ge\kappa_0$, all
$\psi_*,Q_*>0$ and every $\vv x\in U$, there is an integer point $\mv
r\in\Z^{n+1}\nz$ satisfying $(\ref{e:060})$.
\end{lemma}

The convex body (\ref{e:060}) in $\V$ is essentially a set of homogeneous coordinates of points that lies in a certain
\emph{``cell''}\/ near $\vv y(\vv x)\in\cM$. Clearly, the bigger the $|\mv
r|$, the smaller the projective distance of $\mv r$ from the
tangent and transversal planes to $\cM$ (note however that $|\mv r|\ll Q$ in any case). Then, using
Lemma~\ref{l:11} one can efficiently estimate the distance of $\mv r$ from $\cM$. In order to give a formal statement we introduce the following sets.
Let $\cB_{\vv f}(Q_*,\psi_*,\kappa)$ be the set of $\vv x\in U$ such
that there is an $\mv r\in\Z^{n+1}\nz$ satisfying (\ref{e:060}).
Further, let $\cG_{\vv f}(Q_*,\psi_*,\kappa)=U\setminus \cB_{\vv
f}(Q_*,\psi_*,\kappa)$. We will restrict $\vv y$ to $\cG_{\vv f}(Q_*,\psi_*,\kappa)$ for some suitably chosen $\kappa$. This has the benefit of minimizing the distance of $\mv r$ from $\cM$.

\begin{theorem}\label{t:04}
Let $\cM$ be a $C^{(2)}$ submanifold given by $(\ref{e:002})$ and
let $B$ be a ball of radius $r_{B}<\infty$ such that $2\overline
B\subset U$. Then there is an explicit constant $c_0
>2$ such that for any choice of positive
numbers $\psi_*,Q_*,\kappa $ such that $\kappa<1$,
\begin{equation}\label{e:062}
Q_*\ge\max\Big\{\frac{c_0 }{\kappa^2 },\ \frac{c_0^2 }{\kappa^4
r_{B} }\Big\}
\end{equation}
and
\begin{equation}\label{e:063}
\kappa^{-\textstyle\frac{\ddd}{2n-\ddd}}\
Q_*^{-\textstyle\frac{\ddd+2}{2n-\ddd}}\ \le\ \psi_*\ \le\ 1
\end{equation}
we have the inclusion
$$
B\cap\cG_{\vv f}(Q_*,\psi_*,\kappa )\ \subset \ \Delta^{\delta_0}(Q,\psi,2B,\rho)\,,
$$
where $\rho:=c_0
\kappa^{-2}\big(\psi_*^{\ccc }Q_*^{\ddd+1}\big)^{-\frac1{\ddd}}$, \
$\psi=c_0^3\kappa^{-2}\psi_*$, \ $Q=c_0Q_*$ \ and \ $\delta_0=\kappa
c_0^{-2}$.
\end{theorem}

\medskip

Before establishing Theorem~\ref{t:04} we shall give a formal proof
of Lemma~\ref{l:12}.

\medskip

\noindent\textit{Proof of Lemma~\ref{l:12}.} Fix an arbitrary $\vv x\in
U$. Obviously our goal is to show that there is an $\mv
r\in\Z^{n+1}\nz$ satisfying (\ref{e:060}). Recall that
$B=\{\mv r\in\V:\text{(\ref{e:060}) holds}\}$ is a convex
body in $\V$ symmetric about the origin. By Lemmas~\ref{l:07}
and~\ref{l:09}, $B$ is the direct sum of $B_{\mv g}$, $B_{\mv u}$
and $B_{\mv y}$ where the latter are the orthogonal projections of $B$
onto the subspaces $\cV(\mv g)$, $\cV(\mv u)$ and $\cV(\mv y)$
respectively. Furthermore, $B_{\mv g}$ is a ball in
$\cV(\mv g)$ of radius $\psi_*$, $B_{\mv u}$ is a ball in $\cV(\mv
u)$ of radius $(\psi_*^{\ccc } Q_*)^{-\frac1{\ddd}}$, and $B_{\mv
y}$ is a ball in $\cV(\mv y)$ of radius $\kappa Q_*$. Since
$\dim\cV(\mv g)=\ccc $, $\dim\cV(\mv u)=\ddd$ and $\dim\cV(\mv y)=1$ (Lemma~\ref{l:02}),
\begin{equation}\label{e:064}
\begin{array}{rcl}
 \vol(B_{\mv g})&=& 2^{\ccc }\psi_*^{\ccc }v_{\ccc },\\
 \vol(B_{\mv u})&=&2^{\ddd}\Big((\psi_*^{\ccc }Q_*)^{-\frac1{\ddd}}\Big)^{\ddd}v_{\ddd},\\
 \vol(B_{\mv y})&=&2\kappa Q_*.
\end{array}
\end{equation}
Since the subspaces $\cV(\mv g)$, $\cV(\mv u)$ and $\cV(\mv y)$ are
orthogonal,  $\vol(B)=\vol(B_{\mv g})\times\vol(B_{\mv
u})\times\vol(B_{\mv y})$. The latter together with (\ref{e:064})
implies that $ \vol(B)=2^{n+1}\kappa v_{\ccc }v_{\ddd}. $ If
$\kappa>(v_{\ccc }v_{\ddd})^{-1}$ then $\vol(B)>2^{n+1}$ and, by
Minkowski's theorem for convex bodies \cite[\S4.1]{Schmidt-1980}, $B$
contains a non-zero integer point $\mv r=\mv r_\kappa$. This proves
the lemma when $\kappa>(v_{\ccc }v_{\ddd})^{-1}$. Finally notice
that the integer points $\mv r_\kappa$ with
$\kappa_0<\kappa<\kappa_0+1$ are contained in a bounded set.
Therefore there are only finitely many of these points. It follows
that there is a sequence $(\kappa_i)$ with $\kappa_i>\kappa_0$ and
$\kappa_i\to\kappa_0$ as $i\to\infty$ such that the points $\mv
r_{\kappa_i}$ are the same and equal to, say, $\mv r'$. This point
is easily seen to satisfy (\ref{e:060}) with $\kappa=\kappa_0$.
\proofend

\bigskip

\noindent\textit{Proof of Theorem~\ref{t:04}}. Since $2\overline
B\subset U$, there is a constant $C>1$ such that (\ref{e:047}) and
(\ref{e:048}) are fulfilled. We will assume that $\kappa<\kappa_0$
as otherwise, by Lemma~\ref{l:12}, there is nothing to prove. Let
$\psi_*$, $Q_*$ and $\kappa$ satisfy the conditions of
Theorem~\ref{t:04}. Take any $\vv x\in B\cap \cG_{\vv
f}(Q_*,\psi_*,\kappa)$. Our goal is to show that
\begin{equation}\label{e:065}
\text{$\vv x\in B(\vv a/q,\rho)$ \ for some \ $(q,\vv a,\vv b)\in
\cR^{\delta_0}(Q,\psi,2B)$.}
\end{equation}
The constant $c_0$ is defined to absorb various other constants
appearing in the proof. More precisely, we set
\begin{equation}\label{e:066}
  c_0  := \max\Big\{\,\ve_0^{-2};\quad \kappa_0+1;\quad
  16C^2(n+1)^4;\quad  6K(\kappa_0+1)(n+1)^2C^2\,\Big\}\,,
\end{equation}
where $\ve_0=\min\{1,r_{B}\}(4\ddd(n+1)C)^{-1}$ and
$K=14(n+1)^3(C+1)^5\ddd^2$ are the constants appearing in
Lemma~\ref{l:11} and $\kappa_0$ is as in Lemma~\ref{l:12}. By
Lemma~\ref{l:12}, $(\ref{e:060})_{\kappa=\kappa_0}$ has a solution
$\mv r=(r_0,r_1,\dots,r_n)\in\Z^{n+1}\nz$. Without loss of
generality we can assume that $\gcd(r_0,r_1,\dots,r_n)=1$ and that $r_0\ge0$. We set
$q=r_0$, $\vv a=(r_1,\dots,r_{\ddd})$ and $\vv
b=(r_{\ddd+1},\dots,r_{n})$. Obviously $\gcd(q,\vv a,\vv b)=1$. For
the rest of the proof we show that $(q,\vv a,\vv b)$ is the required
point, that is (\ref{e:065}) is satisfied for this choice of $(q,\vv
a,\vv b)$.

\bigskip

\noindent\textbf{Step 1 -- bounds on $|\mv r|$.} Let $\rg$, $\ru$
and $\rf$ be the orthogonal projections of $\mv r$ onto $\cV(\mv
g)$, $\cV(\mv u)$ and $\cV(\mv y)$. By
$(\ref{e:060})_{\kappa=\kappa_0}$ and Lemma~\ref{l:07},
\begin{equation}\label{e:067}
\text{$|\rg|<\psi_*$, \qquad $|\ru|<(\psi_*^{\ccc
}Q_*)^{-1/\ddd}$\qquad and\qquad $|\rf|\le\kappa_0Q_*$.}
\end{equation}
By Lemma~\ref{l:09}, $\rg$, $\ru$ and $\rf$ are pairwise orthogonal.
Therefore, $|\mv r|^2=|\mv r_{\mv g}|^2+|\mv r_{\mv u}|^2+|\mv
r_{\mv y}|^2$. The latter together with (\ref{e:067}) gives
\begin{equation}\label{e:068}
|\mv r|^2<\psi_*^2+(\psi_*^{\ccc }Q_*)^{-2/\ddd}+\kappa_0^2Q_*^2.
\end{equation}
Using the l.h.s.\! of (\ref{e:063}) and the fact that $\kappa <1$
one readily verifies that
\begin{equation}\label{e:069}
 (\psi_*^{\ccc }Q_*)^{-1/\ddd}<Q_*^{1/2}.
\end{equation}
By the r.h.s.\! of (\ref{e:063}),  $\psi_*<1$. Then
(\ref{e:068}) implies that $|\mv r|^2<1+Q_*+\kappa_0^2Q_*^2\le
(\kappa_0^2+1)Q_*^2$. The latter inequality is due to (\ref{e:062}). Hence
$|\mv r|<(\kappa_0+1)Q_*$. Further, notice that the fact that $\vv
x\in \cG_{\vv f}(Q_*,\psi_*,\kappa)$ ensures that $(\ref{e:060})$
does not have a solution in $\Z^{n+1}\nz$. This is only possible if
$|\mv y|^{-1}|\mv y\dt\mv r|\ge\kappa Q_*$. Therefore, $|\mv
r|\ge\kappa\, Q_*$, whence
\begin{equation}\label{e:070}
    \kappa \,Q_*\le|\mv r|\le (\kappa_0+1)Q_*.
\end{equation}

\noindent\textbf{Step 2 -- bounds on $|r_0|$. } We now show the first inequality of
the following relations:
\begin{equation}\label{e:071}
    |r_0|\ge\frac{\kappa \,Q_*}{2(n+1)C}\stackrel{\eqref{e:066}}{\ge}\frac{\kappa \,Q_*}{c_0}\,.
\end{equation}
Assume the contrary. Then, by (\ref{e:070}), there is an
$i_0\in\{1,\dots,n\}$ such that $|r_{i_0}|\ge\kappa (n+1)^{-1}Q_*$.
Let $\mv y=(1,y_1,\dots,y_n)$. Observe that the expression
$r_{i_0}-r_0y_{i_0}$ is one of the coordinates of $\mv y\we\mv r$ in
the standard basis. Therefore,
\begin{equation}\label{e:072}
|\mv y\we\mv r|\ge |r_{i_0}-r_0y_{i_0}|\ge|r_{i_0}|-|r_0y_{i_0}|\ge
\frac{\kappa Q_*}{n+1}-\frac{\kappa Q_*}{2(n+1)C}\times
C=\frac{\kappa Q_*}{2(n+1)}.
\end{equation}
Here we used the fact that $|y_{i_0}|\le C$ implied by (\ref{e:047})
and (\ref{e:048}). In order to derive a contradiction we now obtain
an upper bound for $|\mv y\we\mv r|$.  By Lemma~\ref{l:10} and
(\ref{e:060}), $\frac{1}{|\mv y|}|\mv y\we\mv r|\le
\psi_*+(\psi_*^{\ccc } Q_*)^{-1/\ddd}$. Further, by (\ref{e:063})
and (\ref{e:069}), we get $ \tfrac{1}{|\mv y|}|\mv y\we\mv
r|<1+Q_*^{1/2}< 2Q_*^{1/2}$. The latter together with (\ref{e:047})
and (\ref{e:048}) gives $|\mv y\we\mv r|< 2C(n+1)Q_*^{1/2}$.
Combining the latter with (\ref{e:072}) implies that
$Q_*^{1/2}<4C(n+1)^2/\kappa $. In view of (\ref{e:062}) and
(\ref{e:066}) the latter inequality is contradictory, thus
establishing (\ref{e:071}).

\bigskip

\noindent\textbf{Step 3 -- completion of the proof.} We will first
use Lemmas~\ref{l:11} with
\begin{equation}\label{e:073}
\delta=\frac{\psi_*}{\kappa\,
Q_*}\qquad\text{and}\qquad\ve=\frac{(\psi_*^{\ccc }
Q_*)^{-\frac1{\ddd}}}{\kappa \,Q_*}.
\end{equation}
Therefore, we assume that $\delta\le\ve$ and we begin by verifying
(\ref{e:044}) and (\ref{e:045}).

Obviously, (\ref{e:060}) and (\ref{e:070}) imply (\ref{e:044}).
Further, the l.h.s.\! of (\ref{e:063}) implies that $\ve^2\le\delta$
-- this is the first inequality of (\ref{e:045}). The second
inequality of (\ref{e:045}), that is $\delta\le\ve$, is simply
assumed. Finally, by (\ref{e:069}), $\ve\le (\kappa\,
Q_*^{1/2})^{-1}$. By (\ref{e:062}) and (\ref{e:066}),
$(\kappa\, Q_*^{1/2})^{-1}\le\ve_0$ and hence $\ve\le\ve_0$ --- this
shows the last inequality of (\ref{e:045}). Thus, Lemma~\ref{l:11}
is applicable and therefore, by (\ref{e:046}), there is a point $\vv
x'\in 2B$ such that $\pdist(\vv y_{\vv x'},\vv r)\le K\delta$, where
$\vv r=(r_1/r_0,\dots,r_n/r_0)$. Also, by Lemma~\ref{l:10} together
with (\ref{e:044}), we get $\pdist(\vv y_{\vv x},\vv r)\le 2\,\ve$.
Thus, using (\ref{e:073}) we obtain that
\begin{equation}\label{e:074}
\pdist(\vv y_{\vv x'},\vv r)\le K\,\frac{\psi_*}{\kappa
Q_*}\qquad\text{and}\qquad \pdist(\vv y_{\vv x},\vv r)\le
2\,\frac{(\psi_*^{\ccc } Q_*)^{-\frac1{\ddd}}}{\kappa Q_*}.
\end{equation}
We have shown the validity of (\ref{e:074}) under the assumption
that $\delta\le\ve$. However, note that (\ref{e:074}) also holds
when $\delta>\ve$. Indeed, we simply set $\vv x'=\vv x$. Then
(\ref{e:074}) is an easy consequence of (\ref{e:044}),
Lemma~\ref{l:10} and the fact that $K>2$.

By (\ref{e:047}) and (\ref{e:048}),
\begin{equation}\label{e:075}
|\vv y_{\vv x'}|\le nC\quad\text{and}\qquad |\vv y_{\vv x}|\le nC.
\end{equation}
Also, by (\ref{e:070}) and (\ref{e:071}),
\begin{equation}\label{e:076}
 |\vv r|\le\frac{|\mv r|}{|r_0|}\le \frac{(\kappa_0+1)Q_*}{\frac{\kappa Q_*}{2(n+1)C}}=
 \frac{2(n+1)(\kappa_0+1)C}{\kappa }\,.
\end{equation}
Recall that the euclidean and projective distances are locally
comparable -- see (\ref{e:029}). Then, by (\ref{e:075}),
(\ref{e:076}) and (\ref{e:029}), the l.h.s.\! of (\ref{e:074})
implies that
\begin{equation}\label{e:077}
\begin{array}[b]{rcl}
 |\vv r-\vv y_{\vv x'}| & \le &\displaystyle
 \Big(\frac{2(\kappa_0+1)(n+1)C}{\kappa}+1\Big) \ (nC+1)\
K\,\frac{\psi_*}{\kappa \,Q_*}\\[4ex]
 & \le & \displaystyle\frac{3K(\kappa_0+1)(n+1)^2C^2}{\kappa^2}\ \frac{\psi_*}{Q_*}\
 \le\ \frac{c_0}{2\kappa^2}\ \frac{\psi_*}{Q_*}\ < \
 \frac{c_0}{\kappa^2}\ \frac{\psi_*}{Q_*}
\end{array}
\end{equation}
and similarly the r.h.s.\! of (\ref{e:074}) implies that
\begin{equation}\label{e:078}
 |\vv r-\vv y_{\vv x}| < \frac{c_0}{\kappa^2}\ (\psi_*^{\ccc }
  Q_*^{\ddd+1})^{-\frac1{\ddd}}=\rho.
\end{equation}
Trivially, (\ref{e:078}) implies that $|\vv a/q-\vv x| < \rho$,
that is $\vv x\in B(\vv a/q,\rho)$ whence the l.h.s.\! of
(\ref{e:065}) holds. Also, by (\ref{e:062}), $\rho\le r_{B}$ and
therefore $\vv a/q\in 2B$ . Further, using the triangle inequality,
the Mean Value Theorem and (\ref{e:048}), we get
$$
\begin{array}{rcl}
 |f_l(\vv a/q)-b_l/q| & \le & |f_l(\vv a/q)-f_l(\vv x')|+|f_l(\vv x')-b_l/q| \\[1ex]
 & \le& C|\vv a/q-\vv x'|+|f_l(\vv x')-b_l/q|
  \le  C|\vv r-\vv y_{\vv x'}|\ \stackrel{\eqref{e:077}}{\le} \displaystyle C\, \frac{c_0}{\kappa^2}\
 \frac{\psi_*}{Q_*}
\end{array}
$$
This implies that $|qf_l(\vv a/q)-b_l|\stackrel{\eqref{e:070}}{<}
(\kappa_0+1)C\,c_0\kappa^{-2}\psi_*<c_0^3\kappa^{-2}\psi_*=\psi$.
Trivially, (\ref{e:070}) and (\ref{e:071}) give $\delta_0Q\le q\le
Q$. Thus, $(q,\vv a,\vv b)\in\cR^{\delta_0}(Q,\psi,2B)$
and the r.h.s.\! of (\ref{e:065}) is established. This completes the
proof of Theorem~\ref{t:04}.\proofend

\mysubsection{Uniform version of Theorem~\ref{t:04}}

Within Theorem~\ref{t:04} the constant $c_0$ depends on $B$. Now
restricting $B$ to lie in a compact ball $B_0\subset U$ gives the
following version of Theorem~\ref{t:04} in which $c_0$ is
independent of $B$.

\begin{theorem}\label{t:05}
Let $\cM$ be a $C^{(2)}$ submanifold given by $(\ref{e:002})$ and
let $B_0$ be a compact subset of $U$. Then there is a constant
$c_0=c_0(B_0)>1$ such that for any choice of positive numbers
$\psi_*,Q_*,\kappa $ satisfying $\kappa <1$, $(\ref{e:063})$ and
\begin{equation}\label{e:079}
Q_*\ge 4c_0^2\kappa^{-4}
\end{equation}
for any ball $B\subset B_0$ we have that
\begin{equation}\label{e:080}
\tfrac12B\cap\cG_{\vv f}(Q_*,\psi_*,\kappa )\ \subset \ \Delta^{\delta_0}(Q,\psi,B,\rho)\,,
\end{equation}
where $\rho:=c_0
\kappa^{-2}\big(\psi_*^{\ccc }Q_*^{\ddd+1}\big)^{-\frac1{\ddd}}$, \
$\psi=c_0\kappa^{-2}\psi_*$, \ $Q=c_0Q_*$ \ and \ $\delta_0=\kappa
c_0^{-1}$.
\end{theorem}

\noindent\textit{Proof.} Since $B_0\subset U$ and $U$ is open, for
every $\vv x\in B_0$ there is a ball $B_{\vv x}$ centred at $\vv x$
such that $2B_{\vv x}\subset U$. The collection of balls $\{B_{\vv
x}:\vv x\in B_0\}$ is obviously a cover of $B_0$. Since $B_0$ is
compact, there is a finite subcover $\cC=\{B_1,\dots,B_t\}$. Any
$B_i\in\cC$ satisfies the conditions of Theorem~\ref{t:04}. Let
$c_{0,i}$ be the constant $c_0$ arising from Theorem~\ref{t:04} when
$B=B_i$. Set
$$
 c_0=\frac{\max_{1\le i\le t}c_{0,i}^3}{\min\big\{1,\min_{1\le i\le t}r_{B_i}\big\}}\,,
$$
where $r_{B_i}$ is the radius of $B_i$. Let $\psi_*$, $Q_*$ and
$\kappa$ satisfy the conditions of Theorem~\ref{t:05}. Then, by the
choice of $c_0$ and by Theorem~\ref{t:04}, it is readily seen that
\begin{equation}\label{e:081}
B_0\cap\cG_{\vv f}(Q_*,\psi_*,\kappa )\ \subset \ \Delta^{\delta_0}(Q,\psi,U,\rho)\,,
\end{equation}
with $\rho:=c_0 \kappa^{-2}\big(\psi_*^{\ccc
}Q_*^{\ddd+1}\big)^{-\frac1{\ddd}}$, \ $\psi=c_0\kappa^{-2}\psi_*$,
\ $Q=c_0Q_*$ \ and \ $\delta_0=\kappa c_0^{-1}$. Now, let $B\subset
B_0$ be a ball. Trivially, if $\vv a/q\not\in B$ then $(1-\rho)B\cap
B(\vv a/q,\rho)=\varnothing$. By (\ref{e:079}), $\rho<1/2$.
Therefore, $\frac12B\subset(1-\rho)B$ and $\tfrac12B\cap B(\vv
a/q,\rho)=\varnothing$ if $\vv a/q\not\in B$. Therefore,
(\ref{e:080}) is implied by (\ref{e:081}) and the proof is
complete.\proofend

\section{Integer points in `random' parallelepipeds}\label{iprp}

\mysubsection{Main problem and result}\label{weight}

By Minkowski's theorem on linear forms, any parallelepiped $\Pi$ in
$\R^{k}$ symmetric about the origin contains a non-zero integer
point provided that the volume of $\Pi$ is bigger than $2^{k}$. The
latter condition is in general best possible, though $\Pi$ might contain
a non-zero integer point otherwise. Suppose $\Pi(\vv x)$ is a smooth
family of parallelepipeds of small volume, where $\vv
x\in B$, a ball in $\R^{\ddd}$. In this section we
consider the following

\begin{problem}\label{p4}
What is the probability that $\Pi(\vv x)$ contains a non-zero
integer point?
\end{problem}

\noindent As we shall see in \S\ref{proofs} answering  the question
of Problem~\ref{p4} is absolutely crucial to achieving our main goal
-- establishing Theorem~\ref{t:01}. To avoid ambiguity the
parallelepipeds $\Pi(\vv x)$ will be given by the system of
inequalities
\begin{equation}\label{e:082}
    \Big|\sum_{j=1}^k g_{i,j}(\vv x)\,a_j\Big|\le\theta_i\quad
    (1\le i\le k)\,,
\end{equation}
where $g_{i,j}:U\to\R$ are some functions of $\vv x$ defined on an
open subset $U$ of $\R^{\ddd}$, $a_1,\dots,a_k$ are real variables and
$\bmtheta=(\theta_1,\dots,\theta_k)$ is a fixed $k$-tuple of
positive numbers. We will naturally assume that the matrix $G(\vv
x):=(g_{i,j}(\vv x))_{1\le i,j\le k}$ is non-degenerate for every
$\vv x\in U$. Thus $G:U\to \GL_{k}(\R)$. The above
family of parallelepipeds $\Pi$ is therefore determined by the map
$G$ and the vector of parameters $\bmtheta$. Further, define the
set
$$
\cA(G,\bmtheta)\df=\{\vv x\in U:\exists\ \vv
a=(a_1,\dots,a_k)\in\Z^{k}\smallsetminus\{\vv0\}\text{ satisfying
(\ref{e:082})}\}.
$$
Problem~\ref{p4} restated in terms of $G$ and $\bmtheta$ can now
be formalized as follows: \emph{given a ball $B\subset U$, what is
the probability that a random $\vv x\in B$ belongs to
$B\cap\cA(G,\bmtheta)$?}

\medskip

In this section we introduce a characteristic of $G$ which enables
us to produce an effective bound on the measure of $\cA(G,\bmtheta)$
for arbitrary analytic maps $G$. The characteristic is computable for
various natural classes of $G$ and is indeed computable for the maps
$G$ arising from the applications we have in mind.

As before let $\bmtheta=(\theta_1,\dots,\theta_k)$ be the $k$-tuple
of positive numbers and let $\theta$ be given by
\begin{equation}\label{e:083}
\theta^{k}=\theta_1\cdots\theta_k.
\end{equation}
Thus, $\theta$ is the geometric mean value of
$\theta_1,\dots,\theta_k$. Given $\vv x\in U$ and a linear subspace
$V$ of $\R^k$ with $\codim V=r$, $1\le r<k$, we define the number
\begin{equation}\label{e:084}
 \Theta_{\bmtheta}(\vv x,V):=\min\left\{\theta^{-r}\prod_{i=1}^r\theta_{j_i}\,:\,
 \begin{array}{l} (j_1,\dots,j_r)\in C(k,r) \ \ \text{such
that}\\[0.5ex]
 V\oplus\cV\big(\vv g_{j_1}(\vv x),\dots,\vv g_{j_r}(\vv x)\big)=\R^k\,,
\end{array}\right\}
\end{equation}
where $\cV\big(\vv g_{j_1},\dots,\vv g_{j_r}\big)$ is a vector
subspace of $\R^k$ spanned by $\vv g_{j_1},\dots,\vv g_{j_r}$ and
$C(k,r)$ denotes the set of all subsets of $\{1,\dots,k\}$ of
cardinality $r$. Obviously, since $G(\vv x)\in\GL_{k}(\R)$, the set
in the r.h.s.\! of (\ref{e:084}) is not empty and thus
$\Theta_{\bmtheta}(\vv x,V)$ is well defined and positive. We will
be interested in the local behavior of $\Theta_{\bmtheta}(\vv x,V)$
in a neighborhood a point $\vv x_0$ by looking at
\begin{equation}\label{e:085}
\widehat\Theta_{\bmtheta}(\vv x_0,V):=\liminf_{\vv x\to\vv
x_0}\Theta_{\bmtheta}(\vv x,V)\qquad\text{and}\qquad
\widehat\Theta_{\bmtheta}(\vv
x_0):=\sup_V\widehat\Theta_{\bmtheta}(\vv x_0,V),
\end{equation}
where the latter supremum is taken over all linear subspaces $V$ of
$\R^k$ with $1\le\codim V<k$. The number
$\widehat\Theta_{\bmtheta}(\vv x_0)$ will be referred to as the
\emph{$\bmtheta$-weight of $G$ at $\vv x_0$}. The following
statement represents the main result of this section.

\begin{theorem}[Random parallelepipeds theorem]\label{t:06}
Let $U$ be an open subset of\/ $\R^{\ddd}$, $G:U\to\GL_{k}(\R)$ be
an analytic map and $\vv x_0\in U$. Then there is a ball $B_0\subset
U$ centred at $\vv x_0$ and constants $K_0,\alpha>0$ such that for
any ball $B\subset B_0$ there is a constant $\delta=\delta(B)>0$
such that for any $k$-tuple $\bmtheta=(\theta_1,\dots,\theta_k)$ of
positive numbers
\begin{equation}\label{e:086}
\mu_{\ddd}\Big(B\cap\cA(G,\bmtheta)\Big)\le K_0\,\Big(1+\sup_{\vv
x\in B}\widehat\Theta_{\bmtheta}(\vv
x)^\alpha/\delta^\alpha\Big)\,\theta^{\alpha}\,\mu_{\ddd}(B)\,.
\end{equation}
\end{theorem}

\mysubsection{Auxiliary statements}

We will derive Theorem~\ref{t:06} from a general result due to
Kleinbock and Margulis. This will require translating the problem
into the language of lattices. We proceed with further notation.
Given a lattice $\Lambda\subset \R^k$, let
$\delta(\Lambda):=\min_{\vv v\in\Lambda\nz}|\vv v|_\infty$.
Thus, $\delta$ is a map on the space of lattices. Then the set $\cA(G,\bmtheta)$ can be
straightforwardly rewritten using this $\delta$-map as follows:
$$
 \cA(G,\bmtheta)\, :=\,
 \big\{\,\vv x\in U\,:\,
 \delta(\diag(\bmtheta)^{-1}G(\vv x)\Z^k)\le1\big\},
$$
where $\diag(\bmtheta)$ denotes the diagonal $k\times k$ matrix with
$\bmtheta$ on the diagonal. In order to see this simply multiply the
$i$-th inequality of (\ref{e:082}) by $\theta_i^{-1}$. Then it is
readily seen that the fact $\vv x\in\cA(G,\bmtheta)$ is equivalent
to the existence of $\vv a\in\Z^k\nz$ such that
$|\diag(\bmtheta)^{-1}G(\vv x)\vv a|_\infty\le 1$. The latter is
obviously the same as saying that the lattice
$\diag(\bmtheta)^{-1}G\Z^k$ has a non-zero vector of norm $\le 1$,
that is $\delta(\diag(\bmtheta)^{-1}G(\vv x)\Z^k)\le 1$.

The map $\delta$ obviously satisfies the property that
$\delta(x\Lambda)=x\delta(\Lambda)$ for any lattice $\Lambda$ and
any positive scalar $x$. Therefore, multiplying
$\delta(\diag(\bmtheta)^{-1}G(\vv x)\Z^k)\le 1$ through by $\theta$
(see (\ref{e:083}) for the definition of $\theta$), we get the
equivalent inequality $\delta(g_{\vv t}G(\vv x)\Z^k)\le \theta$,
where $g_{\vv t}=\diag\{t_1,\dots,t_k\}$ and
\begin{equation}\label{e:087}
 t_i:=\theta/\theta_i\quad(1\le i \le k)\,.
\end{equation}
Note that $\det g_{\vv t}=1$. To sum up,
\begin{equation}\label{e:088}
\cA(G,\bmtheta)=\{\vv x\in U:\delta\big(h(\vv x)\Z^k\big) \le
\theta\}\,,\qquad\text{where $ h(\vv x)=g_{\vv t}G(\vv x) $}.
\end{equation}

As we have mentioned above the proof of Theorem~\ref{t:06} will be
based on a result due to Kleinbock and Margulis. In order to state
this result we recall various definitions from
\cite{Kleinbock-Margulis-98:MR1652916}. Let $U$ be an open subset of
$\R^{\ddd}$, $f:U\to\R$ be a continuous function and let
$C,\alpha>0$. The function $f$ is called {\it $(C,\alpha)$-good on
$U$}\/ if for any open ball $B\subset U$ the following is satisfied
\begin{equation}\label{e:089}
    \forall\,\ve > 0\qquad \mu_{\ddd}\Big\{\vv x\in B : |f(\vv x)| < \ve \,\sup_{\vv x\in B}|f(\vv x)|\,\Big\}\
    \le\
C\,\ve^\alpha \,\mu_{\ddd}(B).
\end{equation}
Given a $\lambda>0$ and a ball $B=B(\vv x_0,r)\subset\R^{\ddd}$
centred at $\vv x_0$ of radius $r$, $\lambda B$ will denote the ball
$B(\vv x_0,\lambda r)$. Further, $\cC(\Z^{k})$ will denote the set
of all non-zero complete sublattices of $\Z^k$. An integer lattice
$\Lambda\subset\Z^k$ is called \emph{complete}\/ if it contains all
integer points lying in the linear space generated by $\Lambda$.
Given a lattice $\Lambda\subset\R^k$ and a basis $\vv w_1,\dots,\vv
w_r$ of $\Lambda$, the multivector $\vv w_1\we\dots\we\vv w_r$
is uniquely defined up to sign since any two basis of
$\Lambda$ are related by a unimodular transformation. Therefore, the
following \emph{height} function on the set of non-zero lattices is
well defined:
\begin{equation}\label{e:090}
\|\Lambda\|\df=|\vv w_1\we\dots\we\vv w_r|_\infty\,,
\end{equation}
where $|\cdot|_\infty$ denotes the supremum norm on $\We(\R^k)$. The
following result due to Kleinbock and Margulis appears as Theorem
5.2 in \cite{Kleinbock-Margulis-98:MR1652916}.

\begin{theoremKM}
Let $\ddd,k\in\N$, $C,\alpha>0$ and $0<\rho<1$ be given. Let $B$ be a
ball in $\R^{\ddd}$ and $h:3^kB \to \GL_k(\R)$ be given. Assume that for
any $\Lambda \in \cC(\Z^{k})$
\begin{enumerate}
\item[$(i)$] the function $\vv x\mapsto \|h(\vv x)\Lambda\|$ is
$(C,\alpha)$-good\ on $3^{k}B$, and
\item[$(ii)$] $\sup_{\vv x\in B}\|h(\vv x)\Lambda\|\ge\rho$.
\end{enumerate}
Then there is a constant $N_{\ddd}$ depending on $\ddd$ only such that for
any $ \ve >0$ one has
$$
\mu_{\ddd}\Big\{\vv x\in B:\delta\big(h(\vv x)\Z^k\big) \le
\ve\Big\} \le kC (3^{\ddd}N_{\ddd})^{k} \left(\frac\ve \rho
\right)^\alpha \mu_{\ddd}(B).
$$
\end{theoremKM}

Before we proceed with the proof of Theorem~\ref{t:06}, let us
recall some auxiliary statements about $(C,\alpha)$-good functions.

\begin{lemma}[{{\sc Lemma 3.1 in \cite{Kleinbock-Margulis-98:MR1652916}}}]\label{l:13}
Let $U\subset\R^{\ddd}$ be open and $C,\alpha>0$. If $f_1,\dots,f_m$
are $(C,\alpha)$-good functions on $U$ and
$\lambda_1,\dots,\lambda_m\in\R$, then $\max_{i}|\lambda_if_i| $ is
a $(C,\alpha)$-good function on $U$.
\end{lemma}

\begin{lemma}[{{\sc Corollary 3.3 in \cite{Kleinbock-03:MR1982150}}}]\label{l:15}
 Let $\vv f=(f_1,\dots,f_m)$ be a real analytic map from a connected
 open subset $U$ of\/ $\R^{\ddd}$ to $\R^m$. Then for any point $\vv x_0\in U$
 there is a ball $B(\vv x_0)\subset U$ centred at\/ $\vv x_0$ and constants $C,\alpha>0$ such that
 any function
 $
  \alpha_0+\sum_{i=1}^m\alpha_if_i\quad\text{with}\quad
  \alpha_0,\dots,\alpha_m\in\R
 $
is $(C,\alpha)$-good on $B(\vv x_0)$.
\end{lemma}

Also for the purpose of establishing Theorem~\ref{t:06} we now prove
the following technical statement that translates the definition of
$\widehat \Theta_{\bmtheta}(\vv x_0)$ into the language of exterior
algebra. Within this section we refer to \S\ref{aux} assuming that
$n+1=k$.

\begin{lemma}\label{l:16}
Let $r\in\{1,\dots,k-1\}$ and $\vv x_0\in U$. Then for any ball
$B\subset U$ centred at $\vv x_0$ for any non-zero decomposable
multivector $\mv v\in\We^r(\R^k)$ there is a $J\in C(k,r)$ and $\vv
x\in B$ such that
\begin{equation}\label{e:091}
    \theta^{-r}\prod_{j\in J}\theta_{j}\le\widehat\Theta_{\bmtheta}(\vv x_0)
\end{equation}
and
\begin{equation}\label{e:092}
 \bigwedge_{j\in J}\vv g_j(\vv x)\,\dt\mv v\not=\vv 0.
\end{equation}
\end{lemma}

\noindent\textit{Proof.} Let $\mv v\in \We^r(\R^k)$ be a
decomposable multivector with $1\le r<k$. Define $V=\cV(\mv
v^\perp)$, a vector subspace of $\R^k$. By Lemma~\ref{l:02}, $\codim V=r$. Observe that for a fixed $\bmtheta$ the
function $\Theta_{\bmtheta}(\vv x,V)$ of $\vv x$ takes discrete
values. Then, using (\ref{e:085}) it is easy to see that for any
ball $B\subset U$ centred at $\vv x_0$ there is an $\vv x\in B$ such
that $\Theta_{\bmtheta}(\vv x,V)\le\widehat\Theta_{\bmtheta}(\vv
x_0)$. By the definition of $\Theta_{\bmtheta}(\vv x,V)$, there is a
$J=\{j_1,\dots,j_r\}\in C(k,r)$ satisfying (\ref{e:091}) such that $
 V\oplus\cV\big(\vv g_{j_1}(\vv x),\dots,\vv g_{j_r}(\vv
x)\big)=\R^k$, that is since $V=\cV(\mv v^\perp)$, $\cV(\mv v^\perp)\oplus
\cV\big(\vv g_{j_1}(\vv x),\dots,\vv g_{j_r}(\vv x)\big)=\R^k$, whence, by Lemma~\ref{l:04}, $\mv v^\perp\we\big(\We_{j\in
J}\vv g_{j}(\vv x)\big)\not=\vv 0$. Finally,
$$
\vv0\not=|\mv v^\perp\we\big(\bigwedge_{j\in J}\vv g_{j}(\vv
x)\big)|\stackrel{\eqref{e:034}}{=}|(\mv
v^\perp)^\perp\dt\bigwedge_{j\in J}\vv g_{j}(\vv
x)|\stackrel{\eqref{e:032}}{=}|\mv v\dt\bigwedge_{j\in J}\vv
g_{j}(\vv x)|\,.
$$
\proofend

\mysubsection{Proof of Theorem~\ref{t:06}}\label{sec:01}

By (\ref{e:083}) and (\ref{e:087}), we obviously have that
$\prod_{i=1}^kt_i=1$. Therefore, $\det g_{\vv t}=1$ and
\begin{equation}\label{e:093}
\det G(\vv x)=\det h(\vv x)\,,
\end{equation}
where $h$ is given by (\ref{e:088}). Therefore, $h(\vv x)$ is a
map from $U$ to $\GL_{k}(\R)$.

Our next goal is to verify conditions $(i)$ and $(ii)$ of Theorem~KM
for the specific choice of $h$ made by (\ref{e:088}). Fix a
$\Gamma\in\cC(\Z^k)$. Let $r=\dim\Gamma>0$. Fix a basis of $\Gamma$,
say $\vv w_1,\dots,\vv w_r\in\Z^k$. Then $h(\vv x)\vv w_1,\dots,
h(\vv x)\vv w_r$ is a basis of the lattice $h(\vv x)\Gamma$. By
definition (\ref{e:090}),
$$
\|h(\vv x)\Gamma\|=|h(\vv x)\vv w_1\wedge\dots\wedge h(\vv x)\vv
w_r|_\infty.
$$
Given an $l\in\{1,\dots,r\}$, it is readily seen that the
coordinates of $h(\vv x)\vv w_l$ are equal to $t_i\vv g_i(\vv x)\vv
w_l$, \ $i=\overline{1,k}$. Therefore, by (\ref{e:026}), for every
$I=\{i_1<\dots<i_r\}\subset\{1,\dots,k\}$ the $I$-coordinate of
\begin{equation}\label{e:094}
    h(\vv x)\vv w_1\wedge\dots\wedge h(\vv x)\vv w_r\in\We^r(\R^k)
\end{equation}
in the standard basis equals
\begin{equation}\label{e:095}
    \begin{array}[b]{rcl}
  \displaystyle \det\Big(t_{i_j}\vv
g_{i_j}(\vv x)\vv w_l\Big)_{1\le j,l\le r} & \stackrel{}{=} &
\displaystyle\Big(\prod_{j=1}^rt_{i_j}\Big)\det\Big(\vv g_{i_j}(\vv
x)\vv
w_l\Big)_{1\le j,l\le r} \\[2ex]
   & \stackrel{\eqref{e:027}}{=} & \displaystyle\Big(\prod_{j=1}^rt_{i_j}\Big) \Big(\We_{j=1}^r\vv g_{i_j}(\vv
   x)\Big)\dt\Big(\We_{l=1}^r\vv w_l\Big).
\end{array}
\end{equation}
Since $G$ is analytic, the coordinate functions of $\We_{j=1}^r\vv
g_{i_j}(\vv x)$ are analytic. Let $f_1,\dots,f_M$ be the collection
of these functions taken over all possible choices of $r$ and $I$.
Note that this is a finite collection of analytic functions.
Obviously, (\ref{e:095}) is a linear combination of $f_1,\dots,f_M$.
By Lemma~\ref{l:15}, there is a ball $B_{0}$ centred at $\vv x_0$
and positive $C$ and $\alpha$ such that (\ref{e:095}) (regarded as a
function of $\vv x$) is $(C,\alpha)$-good on $3^kB_{0}$ for any
choice of $r$ and $I$. If $B_{0}$ is sufficiently small then, by the
continuity of $G$, we can also ensure the conditions
\begin{equation}\label{e:096}
    |\det G(\vv x)|\ge\frac12\,|\det G(\vv x_0)|\qquad\text{for all }\vv x\in B_{0}
\end{equation}
and
\begin{equation}\label{e:097}
    \max_{1\le j\le k}\ \sup_{\vv x\in B_0}\big|\vv g_j(\vv
    x)\big|<\infty.
\end{equation}
Take any ball $B\subset B_{0}$. Since every coordinate function of
$h(\vv x)\Gamma$ is $(C,\alpha)$-good on $3^kB$, by
Lemma~\ref{l:13}, the map $\vv x\mapsto\|h(\vv x)\Gamma\|$ is
$(C,\alpha)$-good on $3^kB$. This verifies condition $(i)$ of
Theorem~KM. We proceed with establishing condition $(ii)$. This
splits into 2 cases.

\bigskip

\noindent\textbf{Case \,$r<k$ \,: } Let $C'(k,r)$ be the subset of
$C(k,r)$ consisting of $I=\{i_1<\dots<i_r\}\subset\{1,\dots,k\}$
such that
\begin{equation}\label{e:098}
    \theta^{-r}\prod_{j=1}^r\theta_{i_j}\le\widehat\Theta:=\sup_{\vv
x\in B}\widehat\Theta_{\bmtheta}(\vv x).
\end{equation}
It is readily seen that $C'(r,k)$ is non-empty. By (\ref{e:095})
and (\ref{e:098}), for any $I\in C'(r,k)$ we get that
$$
  \|h(\vv x)\Gamma\|\ge\displaystyle\Big(\prod_{i\in I}t_{i}\Big) \Big|\Big(\We_{i\in I}\vv g_{i}(\vv
   x)\Big)\dt\Big(\We_{l=1}^r\vv w_l\Big)\Big|
       \stackrel{\eqref{e:098}}{\ge}     \frac{1}{\rule{0ex}{2.5ex}\,\widehat\Theta\,} \Big|\Big(\We_{i\in I}\vv g_{i}(\vv
   x)\Big)\dt\Big(\We_{l=1}^r\vv w_l\Big)\Big|.
$$
Taking the supremum over all $\vv x\in B$ and then taking the
maximum over all $I\in C'(r,k)$ gives
\begin{equation}\label{e:099}
\sup_{\vv x\in B}\|h(\vv x)\Gamma\|\ge
\frac{1}{\rule{0ex}{2.5ex}\,\widehat\Theta\,}\ \max_{I\in C'(r,k)}\
\sup_{\vv x\in B}\Big|\Big(\bigwedge_{i\in I}\vv g_i(\vv
x)\Big)\cdot \Big(\bigwedge_{l=1}^{r}\vv w_l\Big)\Big|.
\end{equation}
Now, since $\vv w_l$ are integer points,
$\We_{l=1}^r\vv w_l$ has integer coordinates. Since $\vv
w_1,\dots,\vv w_r$ are linearly independent, $\We_{l=1}^r\vv w_l$ is
non-zero and therefore $|\We_{l=1}^r\vv w_l|\ge1$. Dividing the
r.h.s. of (\ref{e:099}) by $|\We_{l=1}^r\vv w_l|$ gives
\begin{equation}\label{e:100}
    \sup_{\vv x\in B}\|h(\vv x)\Gamma\|\ge
\frac{1}{\rule{0ex}{2.5ex}\,\widehat\Theta\,}\ \max_{I\in C'(r,k)}\
\sup_{\vv x\in B}\Big|\Big(\bigwedge_{i\in I}\vv g_i(\vv
x)\Big)\cdot \frac{\vv w_1\wedge\dots\wedge\vv w_r}{|\vv
w_1\wedge\dots\wedge\vv w_r|}\Big|.
\end{equation}
The multivector $ \mv u=|\vv w_1\wedge\dots\wedge\vv w_r|^{-1}\vv
w_1\wedge\dots\wedge\vv w_r $ is unit and decomposable. Thus, taking
the infimum in (\ref{e:100}) over all $\mv u\in\Gr_r(\R^{k})$,
that is over all unit decomposable $r$-vectors $\mv u$ taken up to
sign, gives
\begin{equation}\label{e:101}
    \sup_{\vv x\in B}\|h(\vv
x)\Gamma\|\ge \frac{1}{\rule{0ex}{2.5ex}\,\widehat\Theta\,}\
\inf_{\mv u\in\Gr_r(\R^{k})}\ \max_{I\in C'(r,k)}\ \sup_{\vv x\in
B}\Big|\Big(\bigwedge_{i\in I}\vv g_i(\vv x)\Big)\cdot \mv u \Big|.
\end{equation}

Our next goal is to show that the constant in the r.h.s.\! of
(\ref{e:101}) is positive. To this end, consider the following
functions of $\vv x\in B$ and $\mv u\in\We^r(\R^k)$ given by
\begin{equation}\label{e:102}
     M_{r,I}(\mv u,\vv x)=\Big|\Big(\bigwedge_{i\in I}\vv g_i(\vv x)\Big)\cdot
\mv u \Big|\qquad\text{and}\qquad  M_{B,r,I}(\mv u)=\sup_{\vv x\in
B} M_{r,I}(\mv u,\vv x).
\end{equation}
For every fixed $\vv x$ the function $M_{r,I}(\mv u,\vv x)$ is the
absolute value of a function linear in $\mv u$. Therefore, using
(\ref{e:097}) one readily gets that $M_{r,I}(\mv u,\vv x)$ is
uniformly continuous in $\mv u$. Henceforth, $M_{B,r,I}(\mv u)$ is
continuous. To prove this formally fix any $\mv
u_0\in\bigwedge^r(\R^k)$ and any $\ve>0$. Then there is an $\eta>0$
such that for all $\mv u\in\bigwedge^r(\R^k)$ satisfying $|\mv u-\mv
u_0|<\eta$
\begin{equation}\label{e:103}
 |M_{r,I}(\mv u,\vv x)-M_{r,I}(\mv u_0,\vv x)|<\ve/2\qquad\text{for all $\vv x\in
 B$}.
\end{equation}
By definition (\ref{e:102}), there is $\vv x_0\in B$ such that
$M_{B,r,I}(\mv u_0)<M_{r,I}(\mv u_0,\vv x_0)+\ve/2$. Therefore,
$$
M_{B,r,I}(\mv u_0)<M_{r,I}(\mv u_0,\vv
x_0)+\ve/2\stackrel{\eqref{e:103}}{\le} M_{r,I}(\mv u,\vv x_0)+\ve
\stackrel{\eqref{e:102}}{\le} M_{B,r,I}(\mv u)+\ve.
$$
Similarly we show the complementary inequality, namely that
$M_{B,r,I}(\mv u_0)>M_{B,r,I}(\mv u)-\ve$. Therefore,
$|M_{B,r,I}(\mv u)-M_{B,r,I}(\mv u_0)|<\ve$ for all $\mv u$
satisfying $|\mv u-\mv u_0|<\eta$. This proves the continuity of
$M_{B,r,I}(\mv u)$. Further, define
$$
 M_{B,r}(\mv u):=\max_{I\in C'(r,k)}M_{B,r,I}(\mv u).
$$
This is also a continuous function of $\mv u$ as the maximum of a
finite number of continuous functions. By Lemma~\ref{l:16} and the
definition of $M_{B,r}(\mv u)$, $M_{B,r}(\mv u)>0$ for
all decomposable multivectors $\mv u\in\bigwedge^r(\R^k)$. Since the
Grassmannian $\Gr_r(\R^{k})$ is compact and $M_{B,r}(\mv u)$ is
continuous, there is a $\mv u_{B,r}\in\Gr_r(\R^{k})$ such that
$$
\inf_{\mv u\in\Gr_r(\R^{k})}\ \ M_{B,r}(\mv u)=M_{B,r}(\mv
u_{B,r})>0.
$$
Thus, (\ref{e:101}) implies that
$$
    \sup_{\vv x\in B}\|h(\vv
x)\Gamma\|\ge \frac{1}{\rule{0ex}{2.5ex}\,\widehat\Theta\,}\
M_{B,r}(\mv u_{B,r})\ge
\frac{1}{\rule{0ex}{2.5ex}\,\widehat\Theta\,}\ M_B
$$
for any $\Gamma\in\cC(\Z^k)$ with $\dim\Gamma<k$, where
$M_B=\min\limits_{1\le r<k}M_{B,r}(\mv u_{B,r})>0$.

\bigskip

\noindent\textbf{Case $r=k$\,: } Now we assume that $\dim\Gamma=k$.
Since $\Gamma$ is complete, $\Gamma=\Z^{k}$ and
therefore the standard basis of $\R^k$, say $\vv e_1,\dots,\vv e_k$,
is also a basis of $\Gamma$. Therefore, (\ref{e:094}) is exactly
$\pm\det h(\vv x)$. Further,
$$
\sup_{\vv x\in B}\|h(\vv x)\Gamma\|=\sup_{\vv x\in B}|\det h(\vv
x)|\stackrel{\eqref{e:093}}{=}\sup_{\vv x\in B}|\det G(\vv x)|\
\stackrel{\eqref{e:096}}{\ge} \ \frac12\,|\det G(\vv x_0)|>0.
$$

\bigskip

\noindent\textbf{Final step. } The upshot of the above discussion is
that for any $\Gamma\in\cC(\Z^k)$
\begin{equation}\label{e:104}
\sup_{\vv x\in B}\|h(\vv
x)\Gamma\|\ge\min\Big\{\frac12,\,\frac12\,|\det G(\vv
x_0)|,\frac{M_B}{\rule{0ex}{2.5ex}\,\widehat\Theta\,}\Big\}=\rho>0.
\end{equation}
This verifies condition $(ii)$ of Theorem~KM. Further, using the
trivial inequality $\min\{|x|,|y|,|z|\}^{-1}\le
|x|^{-1}+|y|^{-1}+|z|^{-1}$ we deduce from (\ref{e:104}) that
\begin{equation}\label{e:105}
\textstyle\rho^{-\alpha} \le  2^{\alpha}+2^{\alpha}|\det G(\vv
x_0)|^{-\alpha}+\left(\frac{\,\widehat\Theta\,}{M_B}\right)^\alpha
  =
2^\alpha\Big(1+|\det G(\vv
x_0)|^{-\alpha}\Big)\left(1+\left(\frac{\,\widehat\Theta\,}{\delta}\right)^\alpha\right),
\end{equation}
where $\delta=\delta(B)$ is implied by (\ref{e:105}). By
(\ref{e:088}) and Theorem~KM (with $\ve=\theta$), we now obtain
(\ref{e:086}) with $K_0 = 2^\alpha kC (3^{\ddd}N_{\ddd})^{k}
\Big(2^{\alpha}+2^{\alpha}|\det G(\vv x_0)|^{-\alpha}\Big)$.
Obviously, $K_0$ is independent of $B$. The proof of
Theorem~\ref{t:06} is thus complete.

\mysubsection{Hierarchic families of parallelepipeds}

It is in general possible but not straightforward to give bounds on
the $\bmtheta$-weight of $G$. In this subsection we introduce a
condition on $G$ that enables us to give a clear-cut estimate for
$\widehat\Theta_{\bmtheta}(\vv x)$ and produce an interesting
corollary of Theorem~\ref{t:06}. Let $B$ be a ball in $U$. We will
say that $G$ is \emph{hierarchic on $B$} if for any vector subspace $V$
of $\R^k$ of $\codim V=r$ the set
\begin{equation}\label{e:106}
\Big\{\,\vv x\in B\,:\ V\oplus\cV\big(\vv g_1(\vv x),\dots,\vv
g_r(\vv x)\big)=\R^k\,\Big\}
\end{equation}
is dense in $B$.

\begin{lemma}\label{l:17}
If $G:U\to\GL_k(\R)$ is hierarchic on a ball $B_0\subset U$ then for
any $k$-tuple $\bmtheta=(\theta_1,\dots,\theta_k)$ of positive
numbers and any $\vv x_0\in U$
\begin{equation}\label{e:107}
    \widehat\Theta_{\bmtheta}(\vv x_0)\ \le\ \tilde\Theta:=\max_{1\le
    r\le k-1}\ \frac{\theta_1\cdots\theta_{r}}{\theta^r}\,.
\end{equation}
\end{lemma}

\noindent\textit{Proof.} Fix any $\vv x_0\in U$. In order to prove
(\ref{e:107}) it suffices to show that
$\widehat\Theta_{\bmtheta}(\vv x_0,V)\le\tilde\Theta$ for every
subspace $V$ of $\R^k$ with $\codim V=r\in\{1,\dots,k-1\}$. Since
the set (\ref{e:106}) is dense in $U$,
$\theta^{-r}\prod_{j=1}^r\theta_{j}$ belongs to the set in the
r.h.s. of (\ref{e:084}) for points $\vv x\in U$ arbitrarily close to
$\vv x_0$. This means that $\Theta(\vv x,V)\le
\theta^{-r}\prod_{j=1}^r\theta_{j}\le\Tilde\Theta$ for points $\vv
x$ arbitrarily close to $\vv x_0$. Therefore, by (\ref{e:085}),
$\widehat\Theta_{\bmtheta}(\vv x_0,V)\le\tilde\Theta$ and the proof
is complete. \proofend

\bigskip

The following example of hierarchic maps will be utilized to sharpen
Theorem~\ref{t:01} is \S\ref{ext}.

\begin{lemma}\label{l:18}
Let $G=\big(g^{(i-1)}_j\big)_{1\le i,j\le k}:U\to\GL_k(\R)$ be the Wronski matrix of analytic linearly independent over $\R$ functions $g_1,\dots,g_k:U\to\R$ defined on an interval $U$ in $\R$. Then $G$ is hierarchic on $U$.
\end{lemma}

\noindent\textit{Proof.} Recall the well known fact that $r$
analytic functions of a real variable are linearly dependent if and
only if their Wronskian is identically zero -- see, for example,
\cite{Bocher-00:MR1503482}. Let $\vv g=(g_1,\dots,g_k)$. Take any
non-trivial vector subspace $V$ of $\R^k$ with $\codim V=r\le k-1$.
We will verify that the set (\ref{e:106}) is dense in $U$ by
showing that its complement is countable. Let $\mv v_1,\dots,\mv
v_{r}$ be a basis of $V^\perp$. Define $\mv v:=\mv
v_1\we\dots\we\mv v_{r}$. Then, by Lemma~\ref{l:03},
$V=\cV(\mv v^\perp)$.

Let $S(V)$ denote the complement of the set (\ref{e:106}).
Obviously, the point $x$ belongs to $S(V)$ if and only if $V\cap
\cV(\vv g(x)\we\dots\we\vv g^{(r)}(x))\not=\varnothing$. By
Lemma~\ref{l:04}, this is equivalent to $\mv v^\perp\we(\vv
g(x)\we\dots\we\vv g^{(r)}(x))=0$ and, by (\ref{e:034}) and the fact
that $\mv v=\mv v_1\we\dots\we\mv v_{r}$, this further gives
\begin{equation}\label{e:108}
 (\mv v_1\we\dots\we\mv v_r)\dt(\vv g(x)\we\dots\we\vv g^{(r)}(x))=0.
\end{equation}
By the Laplace identity (\ref{e:027}), the latter is exactly the
Wronskian of the functions $\eta_i(x)=\vv g(x)\dt\mv v_i$. Since
$\mv v_1,\dots,\mv v_r$ are linearly independent vectors, the
functions $\eta_1,\dots,\eta_r$ ($1\le i\le r$) are linearly
independent over $\R$. Therefore, the Wronskian of
$\eta_1,\dots,\eta_r$ is not identically zero and, as an analytic
function, it can vanish only on a countable subset of $U$.
Therefore, the set $S(V)$ is at most countable and the proof is
complete.\proofend

\medskip

In view of Lemmas~\ref{l:17} and \ref{l:18} specializing
Theorem~\ref{t:06} to the Wronski matrix gives

\begin{theorem}\label{t:07}
Let $g_1,\dots,g_k$ be a collection of real analytic linearly
independent over $\R$ functions defined on an interval $U\subset\R$.
Let $x_0\in U$ be a point such that the Wronskian
$W(g_1,\dots,g_k)(x_0)\not=0$. Then there is an interval $I_0$
centred at $x_0$ and positive constants $K_0$ and $\alpha$
satisfying the following property. For any interval $J\subset I_0$
there is a constant $\delta=\delta(J)$ such that for any positive
$\theta_1,\dots,\theta_{k}$ the set
$$
\left\{ x\in J \, :\,
\begin{array}{l}
    \exists\ (a_1,\dots,a_k)\in\Z^{k}\nz \quad \text{satisfying} \\[0.2ex]
    |a_1g_1^{(i)}(x)+\dots+a_kg_k^{(i)}(x)|< \theta_i\quad \forall\ i=1,\dots,k
\end{array}
\right\}
$$
has Lebesgue measure at most
$
 K_0\big(1+(\delta^{-1}\widetilde\Theta)^{\alpha}\big)\theta^\alpha|J|,
$
where $|J|$ is the length
of $J$,
$$
\theta=(\theta_1\dots\theta_{k})^{1/k}\qquad\text{and}\qquad\widetilde\Theta:=\max_{1\le
    r\le k-1}\frac{\theta_1\cdots\theta_{r}}{\theta^r}.
$$
\end{theorem}

\bigskip

\noindent The following even more explicit estimate for $\tilde\Theta$ is now given.

\begin{lemma}\label{l:19}
Let $\theta_1\le\theta_2\le\dots\le \theta_{k-1}\le\theta_k$. Then
$\tilde\Theta\le(\theta_{k-1}/\theta_k)^{1/k}\le1$.
\end{lemma}

\noindent\textit{Proof.} By definition, there is an $r<k$ such that
$\tilde\Theta= \theta_1\cdots\theta_{r}/\theta^{r}$. Raise the
latter equation to the power $k$ and substitute
$\theta_1\dots\theta_k$ for $\theta^k$. This way we obtain
$$
\tilde\Theta^k=\frac{\theta_1^k\cdots\theta_{r}^k}{\theta_1^r\cdots\theta_k^r}=
\frac{\overbrace{\theta_1\dots\theta_1}^k\cdot\ldots\cdot\overbrace{\theta_{r}\dots\theta_r}^k}
{\ \
\underbrace{\theta_1\dots\theta_1}_r\cdot\ldots\cdot\underbrace{\theta_{k}\dots\theta_k}_r\
\ }\,.
$$
Obviously the numerator and the denominator of the above fraction
have the same number of multiples. Also, by the
conditions of the lemma, any multiple in the numerator is not bigger
than the corresponding multiple in the denominator in the same
place. This gives that $\tilde\Theta^k\le \theta_{r}/\theta_k$.
Furthermore, since $r<k$, $\theta_r\le\theta_{k-1}$ and
so $\tilde\Theta^k\le\theta_{k-1}/\theta_k$, whence the lemma
readily follows.\proofend

\section{\label{proofs}The proof of main result: Theorem~1.4}

\mysubsection{Localisation and outline proof}

Using standard covering arguments we establish the following lemma,
which \mbox{allows} us to impose a convenient condition on $B_0$
while establishing Theorem~\ref{t:01}.

\begin{lemma}\label{l:20}
Let $\cC$ be a collection of non-empty compact balls contained in $U$
such that $U= \bigcup_{B_0\in \cC}\tfrac12B^\circ_0$, where
$B^\circ_0$ denotes the interior of $B_0$. Then the validity of the
statement of Theorem~\ref{t:01} for all $B_0\in\cC$ implies the
validity of the statement of Theorem~\ref{t:01} for arbitrary
compact ball $B_0$ in $U$.
\end{lemma}

\noindent\textit{Proof}. Fix an arbitrary compact ball $B_0\subset
U$. Since $\{\frac12B^\circ:B\in\cC\}$ is an open cover of $B_0$,
there is a finite subcollection of $\cC$, say
$\cC_0=\{B_{0,1},\dots,B_{0,N}\}$, such that
\begin{equation}\label{e:109}
    \textstyle B_0\subset\bigcup_{i=1}^N\tfrac12B_{0,i}.
\end{equation}
We may assume that every ball in this subcollection is of positive
radius. For $i=1,\dots,N$ let $k_{0,i}$, $\rho_{0,i}$ and
$\delta_{0,i}$ be the constants $k_0$, $\rho_0$ and $\delta_0$
arising from Theorem~\ref{t:01} when $B_0=B_{0,i}$. Also let $r_0$
be the radius of the smallest ball in $\cC_0$. Clearly, $r_0$ is positive. Define
$$
\rho_0=\max_{1\le i\le N}\rho_{0,i},\ \ \
  \delta_0=\min_{1\le i\le N}\delta_{0,i},\ \ \
  k_0=\min_{1\le i\le N}k_{0,i},
$$
and take any ball $B\subset B_0$. Note that verifying (\ref{e:004})
for some suitable choice of $C_0$ and $Q_0$ would complete the proof
of Lemma~\ref{l:20}. This splits into 2 cases.

\bigskip

\noindent\textbf{Case (i):} Assume that $r(B)$, the radius of $B$,
satisfies $r(B)\le\tfrac14 r_0$. By (\ref{e:109}) and the inclusion
$B\subset B_0$, there is a $B_{0,i}\in\cC_0$ such that $\tfrac12
B_{0,i}\cap B\not=\varnothing$. Then, since $r(B_{0,i})\ge r_0$ and
$r(B)\le\tfrac14r_0$, $B\subset B_{0,i}$ and the
validity of (\ref{e:004}) becomes obvious.

\bigskip

\noindent\textbf{Case (ii):} Assume that $r(B)>\tfrac14 r_0$. In
this case the idea is to pack $B$ with sufficiently many disjoint
balls of radius $\le \tfrac14r_0$ and apply Case~(i) to each of
these balls. The formal procedure is as follows.

Let $\cC'=\{B_1,\dots,B_M\}$ be a maximal collection of pairwise
disjoint balls centred in $\tfrac12B$ of common radius
$r(B_i)=\tfrac18r_0$. The existence of $\cC'$ is readily seen.
Obviously $\cC'$ is non-empty and, by construction, any ball
$B_i\in\cC'$ is contained in $B$. Let $\vv x\in \tfrac12B$. By the
maximality of $\cC'$, the ball $B(\vv x,\tfrac18r_0)$ may not be
pairwise disjoint with all the balls in $\cC'$. Therefore, $\vv x\in
2B_i$ for some $B_i\in\cC'$. It follows that
$\tfrac12B\subset\bigcup_{i=1}^M 2B_i$. Hence,
\begin{equation}\label{e:110}
 2^{-\ddd}\mu_{\ddd}(B)=\mu_{\ddd}(\tfrac12B)\le \sum_{i=1}^M\mu_{\ddd}(2B_i)= 2^{\ddd}\sum_{i=1}^M\mu_{\ddd}(B_i).
\end{equation}
Since every $B_i\in\cC'$ is of radius $<\tfrac14r_0$, we are within
Case~(i). This means that there exist constants $C_{0,i}>0$ and
$Q_{0,i}>0$ such that for all $Q\ge Q_{0,i}$ and all $\psi$
satisfying the inequalities $C_{0,i}Q^{-1/\ccc }< \psi <
C_{0,i}^{-1}$
\begin{equation}\label{e:111}
\mu_{\ddd}\big(\Delta^{\delta_0}(Q,\psi,B_i,\rho)\cap
B_i\,\big) \ \ge \ k_0 \, \mu_{\ddd}(B_i) \, .
\end{equation}
Now define $ C_0=\max_{1\le i\le M}C_{0,i}$, $Q_0=\max_{1\le i\le
M}Q_{0,i}$. Then (\ref{e:111}) holds whenever (\ref{e:003}) is
satisfied and $Q>Q_0$. Using the disjointness of balls in $\cC'$ and
the fact that $\bigcup_{i=1}^MB_i\subset B$ we get from
(\ref{e:111}) that
$$
\begin{array}{rcl}
\displaystyle \mu_{\ddd}\big(\Delta^{\delta_0}(Q,\psi,B,\rho)\cap
B\,\big) & \ge & \sum_{i=1}^M
\mu_{\ddd}\big(\Delta^{\delta_0}(Q,\psi,B_i,\rho)\cap B_i\,\big) \\[2ex] & \stackrel{\eqref{e:111}}{\ge} & \sum_{i=1}^Mk_0 \, \mu_{\ddd}(B_i)
\  \stackrel{\eqref{e:110}}{\ge} \ 4^{-\ddd}k_0 \, \mu_{\ddd}(B)
\, .
\end{array}
$$
This shows (\ref{e:004}) with $k_0$ replaced by $4^{-d}k_0$ and thus
completes the proof.\proofend

\bigskip

\noindent\textit{Outline proof of Theorem~\ref{t:01}.}
Recall that $\cM$ is a non-degenerate analytic submanifold of $\R^n$ given by (\ref{e:002})
and $B_0$ be a compact ball in $U$. By
Lemma~\ref{l:20}, $B_0$ is assumed to be a sufficiently small
ball. The proof contains the following 3 steps.

\bigskip

\noindent\textbf{\phantom{ii}(i)} Firstly, to establish (\ref{e:004}) take any ball $B$ in $B_0$. In view of
Theorem~\ref{t:05}, namely inclusion (\ref{e:080}),
(\ref{e:004}) follows on showing that for sufficiently large $Q_*$
\begin{equation}\label{e:112}
\mu_{\ddd}(\tfrac12B\cap\cG(Q_*,\psi_*,\kappa))\gg\mu_{\ddd}(B).
\end{equation}

\medskip

\noindent\textbf{\phantom{i}(ii)} In order to establish
(\ref{e:112}), for each $\vv x\in B_0$ we circumscribe a
parallelepiped (\ref{e:082}) around the body (\ref{e:060}). This way
the complement of $\cG_{\vv f}(Q_*,\psi_*,\kappa)$ becomes embedded
into the set $\cA(G,\bmtheta)$ appearing in Theorem~\ref{t:06}, thus
giving
\begin{equation}\label{e:113}
\tfrac12B\setminus \cA(G,\bmtheta)\subset\cG_{\vv
f}(Q_*,\psi_*,\kappa)\cap \tfrac12B.
\end{equation}

\medskip

\noindent\textbf{(iii)} On applying Theorem~\ref{t:06} we will
obtain that
$\mu_d(\tfrac12B\cap\cA(G,\bmtheta))\le\tfrac12\mu_d(\tfrac12B)$. In
view of the embedding (\ref{e:113}) it will further imply
(\ref{e:112}) and complete the task.

We now proceed with the details of the proof.

\mysubsection{$G$ and $\bmtheta$}

Let $\mv g,\mv u,\mv y$ be given by (\ref{e:039})--(\ref{e:041}).
For $\cM$ is analytic, $\mv y$ is analytic. Further, the coordinate
functions of $\mv g$ and $\mv u$ are obviously polynomials of
analytic functions and thus are analytic.

\begin{lemma}\label{l:21}
Let $\mv g,\mv u,\mv y$ be as above. Then for every point $\vv
x_0\in U$ there is a ball $B_0\subset U$ centred at $\vv x_0$ and an
analytic map $G:B_0\to\GL_{n+1}(\V)$ with rows $\vv g_1,\dots,\vv
g_{n+1}$ such that for every $\vv x\in B_0$
\begin{equation}\label{e:114}
|\vv g_i(\vv x)|\le 1\quad\text{ for all }i=1,\dots,n+1
\end{equation}
and
\begin{equation}\label{e:115}
\begin{array}{rcl}
 \vv g_i(\vv x)\in\cV\big(\mv g(\vv x)\big)&\ &\text{for}\quad i=1,\dots,\ccc ,\\[1ex]
 \vv g_{i}(\vv x)\in\cV\big(\mv u(\vv x)\big)&\ &\text{for}\quad i=\ccc +1,\dots,n,\\[1ex]
 \vv g_{i}(\vv x)\in\cV\big(\mv y(\vv x)\big)&\ &\text{for}\quad i=n+1.
\end{array}
\end{equation}
\end{lemma}

\noindent\textit{Proof.} Fix any basis $\vv g_1(\vv x_0),\dots,\vv
g_{n+1}(\vv x_0)$ of $\V$ with $|\vv g_i(\vv x_0)|\le 1/2$ for all
$i=1,\dots,n+1$ such that $(\ref{e:115})_{\vv x=\vv x_0}$ is
satisfied. Define
\begin{equation}\label{e:116}
\begin{array}{rcrcl}
 \vv g_i(\vv x)&:=&\displaystyle\frac{1}{|\mv g(\vv x)|^2}\ \mv g(\vv x)\dt\big(\mv g(\vv x)\dt\vv g_i(\vv x_0)\big)&\ &\text{for}\ i=1,\dots,\ccc \,,\\[2ex]
  \vv g_{i}(\vv x)&:=&\displaystyle\frac{1}{|\mv u(\vv x)|^2}\ \mv u(\vv x)\dt\big(\mv u(\vv x)\dt\vv g_{i}(\vv x_0)\big)&\ &\text{for}\ i=\ccc +1,\dots,n\,,\\[2ex]
 \vv g_{i}(\vv x)&:=&\displaystyle\frac{1}{|\mv y(\vv x)|^2}\ \mv y(\vv x)\dt\big(\mv y(\vv x)\dt\vv g_{i}(\vv x_0)\big)&\ &\text{for}\ i=n+1.
\end{array}
\end{equation}
By Lemma~\ref{l:07}, $\vv g_i(\vv x)$ is (up to sign) the orthogonal
projection of $\vv g_i(\vv x_0)$ onto $\cV(\mv g(\vv x))$ for
$i\in\{1,\dots,\ccc \}$, onto $\cV(\mv u(\vv x))$ for $i\in\{\ccc
+1,\dots,n\}$ and $\cV(\mv y(\vv x))$ for $i=n+1$. Obviously, the
maps $\vv g_i$ given by (\ref{e:116}) are well defined and analytic.
Also, by continuity,
\begin{equation}\label{e:117}
    \bigwedge_{i=1}^{n+1}\vv g_i(\vv x)\to\pm\bigwedge_{i=1}^{n+1}\vv
g_i(\vv x_0)\qquad\text{as}\quad \vv x\to\vv x_0.
\end{equation}
Since $\vv g_1(\vv x_0),\dots,\vv g_{n+1}(\vv x_0)$ are linearly
independent, the r.h.s.\! of (\ref{e:117}) is non-zero. Therefore,
there is a neighborhood $B_0$ of $\vv x_0$ such that for all $\vv
x\in B_0$ the l.h.s.\! of (\ref{e:117}) is non-zero. This proves
that $G$ is non-degenerate. In view of the continuity of $\vv g_i$
and the condition $|\vv g_i(\vv x_0)|\le 1/2$ we have $|\vv g_i(\vv
x)|\le 1$ for all $i=1,\dots,n+1$ provided that $B_0$ is
small enough.\proofend

%

\begin{lemma}\label{l:22}
Let $G$ and $B_0$ arise from Lemma~\ref{l:21} and
$\psi_*,Q_*,\kappa$ be any positive numbers. Let
\begin{equation}\label{e:118}
\begin{array}{rcl}
  \theta_1  = \dots=\theta_{\ccc }&=&\psi_*, \\[1ex]
  \theta_{\ccc +1}=\dots=\theta_n&=&(\psi_*^{\ccc }Q_*)^{-1/\ddd},\\[1ex]
  \theta_{n+1}&=&\kappa Q_*.
\end{array}
\end{equation}
Then $(\ref{e:113})$ is satisfied, where $\bmtheta=(\theta_1,\dots,\theta_{n+1})$.
\end{lemma}

\noindent\textit{Proof.} Observe that (\ref{e:113}) is equivalent to
$\tfrac12B\setminus\cG_{\vv f}(Q_*,\psi_*,\kappa)\subset
\tfrac12B\cap\cA(G,\bmtheta)$. By definition, for every point $\vv
x\in \tfrac12B\setminus\cG_{\vv f}(Q_*,\psi_*,\kappa)$ there is a
non-zero integer solution $\mv r$ to the system (\ref{e:060}). Using
(\ref{e:114}), Lemma~\ref{l:08} and Lemma~\ref{l:21} in an obvious
manner implies that (\ref{e:082}) is satisfied when
$(a_1,\dots,a_k)$ is identified with $\mv r$. This exactly means
that $\vv x\in\cA(G,\bmtheta)$ and completes the proof. \proofend

\medskip

We now estimate the $\bmtheta$-weight of $G$ for the above
$G$ and $\bmtheta$. See \S\ref{weight} for its definition.

\begin{lemma}\label{l:23}
Let $\cM$ be a non-degenerate analytic manifold given
by $(\ref{e:002})$. Let $G$ and $B_0$ arise from Lemma~\ref{l:21}
and let $\bmtheta$ be given by $(\ref{e:118})$. Let $\kappa$,
$\psi_*$ and $Q_*$ satisfy the conditions of Theorem~\ref{t:05} and
let
\begin{equation}\label{e:119}
    C_*Q_*^{-1/\ccc }\le \psi_*\le C_*^{-1}
\end{equation}
for some $C_*>1$. Then for any $\vv x_0\in B_0$
\begin{equation}\label{e:120}
\widehat\Theta_{\bmtheta}(\vv x_0)\le (\kappa\, C_*)^{-1/(n+1)}.
\end{equation}
\end{lemma}

\noindent\textit{Proof.} By the definitions of $\theta$ and
$\bmtheta$, \ie{} by (\ref{e:083}) and (\ref{e:118}),
\begin{equation}\label{e:121}
\theta=\kappa^{1/(n+1)}.
\end{equation}
Further, using inequalities (\ref{e:119}) and the assumption $C_*>1$
it is readily seen that
\begin{equation}\label{e:122}
\theta_i\le1\qquad (1\le i\le n)
\end{equation}
Fix any point $\vv x_0\in B_0$ and any vector subspace $V$ of
$\R^{n+1}$ with $\codim V=r\in\{1,\dots,n\}$. Since $\cM$ is
non-degenerate, for every ball $B(\vv x_0)\subset B_0$ centred at
$\vv x_0$ there is a point $\vv x\in B(\vv x_0)$ such that $\mv
y=\mv y(\vv x)\not\in V^\perp$. That is $\cV(\mv y)\not\subset
V^\perp$. The latter is easily seen to be equivalent to $\cV(\mv
y)^\perp\not\supset V$. By Lemma~\ref{l:09} and by (\ref{e:115}), we
see that the first $n$ rows of $G$, which are simply the vectors
$\vv g_1(\vv x),\dots,\vv g_n(\vv x)$, form a basis of $\cV(\mv
y)^\perp$. Thus, $V\not\subset\cV(\vv g_1(\vv x),\dots,\vv g_n(\vv
x))$ and therefore
\begin{equation}\label{e:123}
\dim \big(V+\cV(\vv g_1(\vv x),\dots,\vv g_n(\vv
x))\big)>\dim\cV\big(\vv g_1(\vv x),\dots,\vv g_n(\vv x)\big)=n.
\end{equation}
The latter implies that the l.h.s.\! of (\ref{e:123}) is
equal to $n+1$. Hence there is a subcollection
$J=\{j_1<\ldots<j_r\}\subseteq \{1,\dots,n\}$ satisfying
$
V\oplus\cV\big(\vv g_{j_1}(\vv x),\dots,\vv g_{j_r}(\vv
x)\big)=\R^{n+1}.
$
By (\ref{e:084}),
\begin{equation}\label{e:124}
\begin{array}[b]{rcl}
\Theta_{\bmtheta}(\vv x,V) &\le& \theta^{-r}\prod_{i=1}^r\theta_{j_i}\
\stackrel{\eqref{e:121}}{=} \ \kappa^{-r/(n+1)}\prod_{i=1}^r\theta_{j_i}\\[2ex]
&\stackrel{\eqref{e:122}}{\le}& \kappa^{-r/(n+1)}\max_{1\le i\le r}\theta_{j_i}\\[2ex]
&\stackrel{\eqref{e:118}}{\le}& \kappa^{-r/(n+1)}\max\{\psi_*,(\psi_*^{\ccc }Q_*)^{-1/\ddd}\}\\[2ex]
&\stackrel{\eqref{e:119}}{\le}&\kappa^{-r/(n+1)}\max\{C_*^{-1},C_*^{-1/\ddd}\}\\[2ex]
&\stackrel{C_*>1,\,\kappa<1}{\le}&\kappa^{-1/(n+1)}C_*^{-1/(n+1)}=(\kappa
C_*)^{-1/(n+1)}.
\end{array}
\end{equation}
Recall that $B(\vv x_0)$ can be made arbitrarily small so that $\vv
x$ can be made arbitrary close to $\vv x_0$. Therefore, in view of
the definition (\ref{e:085}) of $\widehat\Theta_{\bmtheta}(\vv
x_0,V)$, (\ref{e:124}) implies that $\widehat\Theta_{\bmtheta}(\vv
x_0,V)\le (\kappa\, C_*)^{-1/(n+1)}$. Finally, since $V$ is
arbitrary non-trivial subspace of $\V$, we obtain (\ref{e:120}).
\proofend

\mysubsection{Completion of the proof of Theorem~\ref{t:01}}

We now proceed to the final phase of the proof of
Theorem~\ref{t:01}. Let $\vv x_0\in U$ be an arbitrary point. Let
$B_0$ be a ball centred at $\vv x_0$ arising from Lemma~\ref{l:21}.
We may assume without loss of generality that $B_0$ is compact.
Further, shrink $B_0$ if necessary to ensure that Theorem~\ref{t:06}
is applicable. Next, let $B$ be an arbitrary ball in $B_0$ and let
$\psi$ and $Q$ satisfy the conditions of Theorem~\ref{t:01}, where
$C_0$ and $Q_0$ are sufficiently large constants.

\medskip

Let $c_0=c_0(B_0)>1$ be the constant arising from Theorem~\ref{t:05}
and let $K_0$, $\alpha$ and $\delta=\delta(\frac12B)$ be the
constants arising from Theorem~\ref{t:06}. Set
\begin{equation}\label{e:125}
    \kappa= (4K_0)^{-\frac{n+1}{\alpha}}.
\end{equation}
Obviously, $0<\kappa<1$ and is independent of $B$. Define
\begin{equation}\label{e:126}
\psi_*=\kappa^2c_0^{-1}\psi,\quad Q_*=c_0^{-1}Q,\quad
\delta_0=\kappa c_0^{-1},\quad \rho=c_0
 \kappa^{-2}\big(\psi_*^{\ccc }Q_*^{\ddd+1}\big)^{-\frac1{\ddd}}.
\end{equation}
It is easily verified that $1/\ccc \le (\ddd+2)/(2n-\ddd)$.
Therefore, (\ref{e:003}) implies that
\begin{equation}\label{e:127}
    C_0Q^{-(\ddd+2)/(2n-\ddd)}< \psi < C_0^{-1}.
\end{equation}
Then, using (\ref{e:126}) and (\ref{e:127}) one readily verifies
(\ref{e:063}) and (\ref{e:079}) provided that $C_0$ and $Q_0$ are
sufficiently large. Therefore, Theorem~\ref{t:05} is applicable and
so (\ref{e:080}) is satisfied. Further, let
$\bmtheta=(\theta_1,\dots,\theta_{n+1})$ be given by (\ref{e:118})
and $G$ be as in Lemma~\ref{l:21}. Then, by Lemma~\ref{l:22} and
(\ref{e:080}), we obtain
\begin{equation}\label{e:128}
 \tfrac12B\setminus \cA(G,\bmtheta) \ \subset \ \tfrac12B\cap\cG_{\vv f}(Q_*,\psi_*,\kappa)
 \ \subset \ \Delta^{\delta_0}(Q,\psi,B,\rho)\cap B.
\end{equation}
Since Theorem~\ref{t:06} is applicable, by (\ref{e:086}), we get
\begin{equation}\label{e:129}
\mu_{\ddd}\Big(\tfrac12B\cap\cA(G,\bmtheta)\Big)\le
K_0\,\Big(1+\widehat\Theta_{\bmtheta}(\tfrac12B)^\alpha/\delta^\alpha\Big)\,\theta^{\alpha}\,\mu_{\ddd}(\tfrac12B)\,,
\end{equation}
where $\widehat\Theta_{\bmtheta}(\frac12B):=\sup_{\vv x\in
\frac12B}\widehat\Theta_{\bmtheta}(\vv x)$. By (\ref{e:003}) and
(\ref{e:126}), (\ref{e:119}) holds with
\begin{equation}\label{e:130}
    C_*=C_0\kappa^2c_0^{-1-1/\ccc}.
\end{equation}
Clearly $C_*>1$ if $C_0$ is sufficiently large. Then, by
Lemma~\ref{l:23}, we get
\begin{equation}\label{e:131}
\widehat\Theta_{\bmtheta}(\tfrac12B)\le (\kappa\, C_*)^{-1/(n+1)}
\le \delta
\end{equation}
provided that $C_0$ and respectively $C_*$ is sufficiently large.
Recall by (\ref{e:121}) that $\theta=\kappa^{1/(n+1)}$. Then, by
(\ref{e:129}) and (\ref{e:131}), we get that
\begin{equation}\label{e:132}
\mu_{\ddd}\Big(\tfrac12B\cap\cA(G,\bmtheta)\Big)\le 2K_0
\kappa^{\alpha/(n+1)}\,\mu_{\ddd}(\tfrac12B) \
\stackrel{\eqref{e:125}}{=} \ \tfrac12\mu_{\ddd}(\tfrac12B).
\end{equation}
Combining (\ref{e:132}) with (\ref{e:128}) gives
$
\mu_d\big(\Delta^{\delta_0}(Q,\psi,B,\rho)\cap B\big)\ge\tfrac12\mu_d(\tfrac12B)=2^{-d-1}\mu_d(B),
$
thus establishing (\ref{e:004}) with $k_0=2^{-d-1}$ and
$\rho_0=(c_0^{n+d+1}\kappa^{-2n})^{1/d}$. The latter constant is
easily deducted from (\ref{e:126}) and is absolute. This completes
the proof of Theorem~\ref{t:01}.

\section{Further theory for curves}\label{ext}

In this section we relax condition (\ref{e:003}) in the case of curves. Namely, the exponent $\frac1m=\frac1{n-1}$ will be replaced by $\frac 3{2n-1}$. The latter allows us to widen the range of $s$ Theorem~\ref{t:03} is applicable by the factor of $\frac n2$.

\mysubsection{Statement of results}

Given an analytic map $\mv y=(y_0,y_1,\dots,y_n):U\to\R^{n+1}$,
where $U\subset\R$ is an interval, let $W_{\mv y}(x)$ denote
the Wronskian of $y_0,y_1,\dots,y_n$.

\begin{theorem}\label{t:08}
Let $\ddd=1$ and the curve $(\ref{e:002})$ satisfies $W_{\mv
y}(x)\not=0$ for all $x\in U$, where $\mv y$ as in $(\ref{e:039})$.
Then Theorem~\ref{t:01} and consequently
Corollary~\ref{corollary:01} remain valid if $(\ref{e:003})$ is
replaced by
\begin{equation}\label{e:133}
    C_0Q^{-3/(2n-1) }< \psi < C_0^{-1}.
\end{equation}
\end{theorem}

Recall that the analytic curve $\cM$ is non-degenerate if and only if the functions $1,y_1,\dots,y_n$ are linearly independent over $\R$. Equivalently, $W_{\mv y}(x)$ is not identically zero. As $\mv y=(1,y_1,\dots,y_n)$ is analytic, the Wronskian $W_{\mv y}(x)$ is analytic too. The non-degeneracy of $\cM$ then implies that $W_{\mv y}(x)\not=0$ everywhere except possibly on a countable set consisting of isolated points. Therefore, the condition ``$W_{\mv y}(x)\not=0$ for all $x\in U$'' imposed in the statement of Theorem~\ref{t:08} is not particularly restrictive if compared to non-degeneracy.

\begin{theorem}\label{t:09}
Let $\cM$ be a non-degenerate analytic curve in\/ $\R^n$. Let
$\psi:\N\to\Rp$ be a monotonic function such that
$q\psi(q)^{(2n-1)/3} \to\infty$ as $q\to\infty$. Then for any
$s\in(\frac12,1)$
\begin{equation}\label{e:134}
 \cH^s(\cS_n(\psi)\cap\cM)=\infty \qquad \text{if} \qquad \sum_{q=1}^\infty \
 q^n\Big(\frac{\psi(q)}q\Big)^{s+n-1 }=\infty.
\end{equation}
and consequently if $\tau=\tau(\psi)$ satisfies $1/n<\tau<3/(2n-1) $
then
\begin{equation}\label{e:135}
\dim \cS_n(\psi)\cap\cM \ge \frac{n+1}{\tau+1}-(n-1) .
\end{equation}
\end{theorem}

The proof of Theorem~\ref{t:09} can be obtained by making  minor and indeed obvious modifications to the proof of Theorem~\ref{t:03}. Below we consider the proof of Theorem~\ref{t:08} only.

\mysubsection{Dual map}\label{dual}

Let the analytic map $\mv y$ be given by (\ref{e:039}). The map $\mv
z:U\to\R^n$ given by
\begin{equation}\label{e:136}
    \mv z(x)=\Big(\mv
y(x)\we\mv y'(x)\we\ldots\we\mv y^{(n-1)}(x)\Big)^\perp
\end{equation}
will be called \emph{dual} to $\mv y$. Obviously, every coordinate
function of $\mv z$ is a polynomial expression of coordinate
functions of $\mv y$ and their derivatives. Therefore,
$\mv z$ is analytic. The following statement describes $\mv z$ via a system of linear
differential equations.

\begin{lemma}\label{l:25}
Let $\mv y$ and $\mv z$ be as above. Then
\begin{equation}\label{e:137}
\left\{    \begin{array}{rcclc}
      \mv z^{(j)}(x)\dt\mv y^{(i)}(x)&=&0, && 0\le i+j\le n-1 \\[1ex]
      \mv z^{(j)}(x)\dt\mv y^{(i)}(x)&=&(-1)^jW_{\mv y}(x),&& i+j=
      n\,.
    \end{array}\right.
\end{equation}
\end{lemma}

\noindent\textit{Proof.} By (\ref{e:136}) and Lemma~\ref{l:03}, we
immediately get that $\mv z(x)\dt\mv y^{(j)}(x)=0$ for all
$j\in\{0,\dots,n-1\}$. On differentiating the latter equations we
obviously obtain the first set of equation in (\ref{e:137}). Now we
compute $\mv z(x)\dt\mv y^{(n)}(x)$:
\begin{eqnarray*}
  \mv z(x)\dt\mv y^{(n)}(x) &\stackrel{\eqref{e:136}}{=}&
  \Big(\mv y(x)\we\mv y'(x)\we\ldots\we\mv y^{(n-1)}(x)\Big)^\perp\dt\mv y^{(n)}(x) \\
   &\stackrel{\eqref{e:031}}{=}&
   \Big(\mv i\dt(\mv y(x)\we\mv y'(x)\we\ldots\we\mv y^{(n-1)}(x))\Big)\dt
   \mv y^{(n)}(x) \\
   &\stackrel{\eqref{e:030}}{=}&
   \Big(\mv e_0\we\mv e_1\we\ldots\we\mv e_n\Big)\dt
   \Big(\mv y(x)\we\ldots\we\mv y^{(n-1)}(x)\we\mv y^{(n)}(x)\Big) \\
   &\stackrel{\eqref{e:027}}{=}&
   \det\Big(\mv e_i\dt\mv y^{(j)}(x)\Big)_{0\le i,j\le n}=
   \det\Big(y^{(j)}_i(x)\Big)_{0\le i,j\le n}=W_{\mv y}(x)\,,
\end{eqnarray*}
where $\mv e_0,\dots,\mv e_n$ is the standard basis of $\V$. This
shows the $j=0$ equation of the second set of equations in
(\ref{e:137}). Then we proceed by induction. Assume that the second
set of inequalities of (\ref{e:137}) holds for $j=j_0\le n-1$. Then
differentiating $\mv z^{(j_0)}(x)\dt\mv y^{(n-j_0-1)}(x)=0$ we get
$$
\begin{array}{rcl}
0 & = & \mv z^{(j_0+1)}(x)\dt\mv y^{(n-j_0-1)}(x)+\mv
z^{(j_0)}(x)\dt\mv y^{(n-j_0)}(x)\\[1ex]
& =& \mv z^{(j_0+1)}(x)\dt\mv y^{(n-j_0-1)}(x)+(-1)^{j_0}W_{\mv
y}(x)\,.
\end{array}
$$
This implies (\ref{e:137}) for $j=j_0+1$ and thus completes the
proof Lemma~\ref{l:25}.\proofend

\begin{lemma}\label{l:26}
Let $\mv y$ and $\mv z$ be as above. Then for all $x$,
$|W_{\mv z}(x)|\ge |W_{\mv y}(x)|^n$.
\end{lemma}

\noindent\textit{Proof.} By (\ref{e:137}) and (\ref{e:027}), it is
easy to see that the inner product in
$\We^{n+1}(\R^{n+1})$
\begin{equation}\label{e:138}
(\mv z(x)\we\mv z'(x)\we\dots\we\mv z^{(n)}(x))\dt(\mv y(x)\we\mv
y'(x)\we\dots\we\mv y^{(n)}(x))
\end{equation}
is the determinant of an $(n+1)\times(n+1)$ triangle matrix with
$\pm W_{\mv y}(x)$ on the diagonal and is equal to $(-1)^{[n/2]}
W_{\mv y}(x)^{n+1}$. Further, recall that
$$
 |\mv y(x)\we\dots\we\mv y^{(n)}(x)|=|W_{\mv y}(x)|\qquad\text{and}\qquad
 |\mv z(x)\we\dots\we\mv z^{(n)}(x)|=|W_{\mv z}(x)|.
$$
Then, applying the Cauchy-Schwarz inequality to the inner product
(\ref{e:138}) gives
$$
  \big|W_{\mv y}(x)\big|^{n+1} \le
  |\mv z(x)\we\dots\we\mv z^{(n)}(x)|\cdot|\mv y(x)\we\dots\we\mv y^{(n)}(x)|
 =  |W_{\mv z}(x)|\cdot|W_{\mv y}(x)|
$$
further implying $|W_{\mv z}(x)|\ge |W_{\mv y}(x)|^n$.\proofend

\mysubsection{Proof of Theorem~\ref{t:08}}

Clearly, Lemma~\ref{l:20} can be used in the context of
Theorem~\ref{t:08}. Then we can assume that $B_0$ is a sufficiently
small interval centred at an arbitrary point $x_0$ in $U$. We may
assume without loss of generality that $B_0$ is compact. Further,
shrink $B_0$ if necessary to ensure that Theorem~\ref{t:06} is
applicable. Next, let $B$ be an arbitrary interval in $B_0$ and let
$\psi$ and $Q$ satisfy the conditions of Theorem~\ref{t:08}, where
$C_0$ and $Q_0$ are sufficiently large constants. Furthermore, in
view of Theorem~\ref{t:01}, without loss of generality we may assume
that
\begin{equation}\label{e:139}
C_0Q^{-3/(2n-1)}<\psi<Q^{-1/n}.
\end{equation}
Let $\mv z$ be dual to $\mv y$ (see \S\ref{dual}). Since $B_0$ is
compact, there is a constant $K_1>1$ such that
\begin{equation}\label{e:140}
    |\mv z^{(i)}(x)|\le K_1\qquad\text{for all $x\in B_0$ and all
    $i\in\{0,\dots,n\}$}.
\end{equation}

Let $c_0=c_0(B_0)>1$ be the constant arising from Theorem~\ref{t:05}
and let $K_0$, $\alpha$ and $\delta=\delta(\frac12B)$ be the
constants arising from Theorem~\ref{t:07}. Set
\begin{equation}\label{e:141}
    \kappa= (2K_1(4K_0)^{\frac{1}{\alpha}})^{-n-1}.
\end{equation}
Obviously, $0<\kappa<1$ and is independent of $B$. Define
$\psi_*,Q_*,\delta_0$ and $\rho$ by (\ref{e:126}) assuming that
$d=1$ and $m=n-1$.

Then, using (\ref{e:126}) and (\ref{e:139}) one readily verifies
(\ref{e:063}) and (\ref{e:079}) provided that $C_0$ and $Q_0$ are
sufficiently large. Therefore, Theorem~\ref{t:05} is applicable and
so (\ref{e:080}) is satisfied.

Take any point $\vv x\in \tfrac12B\setminus\cG_{\vv
f}(Q_*,\psi_*,\kappa)$. Then, by definition, there is a non-zero
integer solution $\mv r$ to system (\ref{e:060}). Observe that $\mv
g=(\mv y\we\mv y')^\perp$. Then, by (\ref{e:137}) and
Lemma~\ref{l:03}, we get $\mv z^{(i)}\in\cV(\mv g)$ for
$i=0,\dots,n-2$. Therefore, by Lemma~\ref{l:08}, (\ref{e:060})
implies that
\begin{equation}\label{e:142}
    |\mv z^{(i)}(x)\dt\mv r|\le K_1\psi_*\qquad\text{for
$i=0,\dots,n-2$}.
\end{equation}
Again, by (\ref{e:137}) and Lemma~\ref{l:03}, $\mv
z^{(n-1)}\in\cV(\mv y^\perp)$. Therefore, we get
\begin{equation}\label{e:143}
\begin{array}[b]{rcl}
 |\mv z^{(n-1)}(x)\dt\mv r| & \stackrel{{\rm Lemma~\ref{l:08}}}{\le}
 &\displaystyle |\mv z^{(n-1)}(x)|\,|\mv y(x)^\perp\dt\mv r|\,
 |\mv y(x)^\perp|^{-1}\\[1ex]
 &\stackrel{(\ref{e:140})}{\le} &K_1\,|\mv y(x)^\perp\dt\mv r|\,
 |\mv y(x)^\perp|^{-1}\\[1ex]
 &\stackrel{(\ref{e:034})}{\le} &K_1\,|\mv y(x)\we\mv r|\,
 |\mv y(x)|^{-1} \stackrel{(\ref{e:039})}{\le} K_1\,|\mv y(x)\we\mv r|.
\end{array}
\end{equation}
Here we have also used the fact that the Hodge operator is an isometry. Using
(\ref{e:126}) and the r.h.s.\! of (\ref{e:139}) we get that
$\psi_*<(\psi_*^{n-1}Q_*)^{-1}$. Therefore, applying
Lemma~\ref{l:10} and (\ref{e:060}) to (\ref{e:143}) further gives
\begin{equation}\label{e:144}
 |\mv z^{(n-1)}(x)\dt\mv r| \le K_1\,(\psi_*+(\psi^{n-1}Q_*)^{-1})\le 2 K_1 (\psi^{n-1}Q_*)^{-1}.
\end{equation}
Finally, arguing the same way as in Step 1 of the proof of
Theorem~\ref{t:04} (see \S\ref{GCT}), one easily verifies that $|\mv
r|^2\le 1+Q_*+\kappa^2Q_*^2$, whence $|\mv r|\le 2\kappa
Q_*$ when $Q_*$ is sufficiently large. Now we trivially get
\begin{equation}\label{e:145}
    |\mv z^{(n)}(x)\dt\mv r|\le |\mv z^{(n)}(x)|\,|\mv r|\le 2K_1\kappa
    Q_*.
\end{equation}

Let $G$ be the Wronski matrix of the dual map $\mv z$. For $W_{\mv
y}(x)\not=0$ for all $x\in U$, by Lemma~\ref{l:26}, $W_{\mv
z}(x)\not=0$ for all $x\in U$, that is $G:U\to\GL_{n+1}(\R)$. The
inequalities (\ref{e:142}), (\ref{e:144}) and (\ref{e:145}) are
equivalent to $x\in\cA(G,\bmtheta)$ with
\begin{equation}\label{e:146}
  \theta_1  = \dots=\theta_{n-1 }=K_1\psi_*,\qquad
  \theta_{n} = 2K_1(\psi_*^{n-1}Q_*)^{-1},\qquad
  \theta_{n+1} = 2K_1\kappa Q_*.
\end{equation}
Thus we have shown that $\tfrac12B\smallsetminus\cG_{\vv f}(Q_*,\psi_*,\kappa)\subset\tfrac12B\cap \cA(G,\bmtheta)$. Hence, $\tfrac12B\setminus \cA(G,\bmtheta) \ \subset \ \tfrac12B\cap\cG_{\vv f}(Q_*,\psi_*,\kappa)$. By (\ref{e:080}),
\begin{equation}\label{e:147}
 \tfrac12B\setminus \cA(G,\bmtheta)
 \ \subset \ \Delta^{\delta_0}(Q,\psi,B,\rho)\cap B.
\end{equation}

By the r.h.s.\! of (\ref{e:139}), $\theta_{i}\le
\theta_{i+1}$ for all $i=1,\dots,n$. Further, by the l.h.s.\! of (\ref{e:139}), $\theta_n\ll Q^{1/2}$.
Further, $\theta_{n+1}\asymp Q$. Then, by Lemma~\ref{l:19},
$\widetilde\Theta\ll Q^{-1/(2n+2)}<\delta$ for sufficiently large
$Q$. Theorem~\ref{t:07}, (\ref{e:141}) and (\ref{e:146}) imply that $\mu_d(\tfrac12B\cap\cA(G,\bmtheta))\le
2K_0\theta^\alpha\le\tfrac12\mu_d(\tfrac12B)$. Combining this with
(\ref{e:147}) gives the required result.

\section{Final comments}\label{final}

In view of the results of this paper, establishing upper bounds for $N(Q,\varepsilon)$ becomes a very topical problem. Unfortunately, non-degeneracy alone is not enough to reverse (\ref{e:005}). A counterexample can be easily constructed by considering the manifolds
$\cM_k=\big\{(x_1,\dots,x_{d-1},x_d,x_d^{k+1},x_d^{k+2},\dots,x_d^{k+m}):
\max_{1\le i\le d}|x_i|<1\big\}$. Nevertheless, requiring that for every $\vv x\in U$ there exists $l\in\{1,\dots,\ccc\}$ such that $\hess f_l(\vv
x)\not=0$ is possibly enough to reverse (\ref{e:004}),
where $\hess f(\vv x)$ denotes the Hessian matrix of $f(\vv x)$.
Any progress with this would have obvious implications for the theory of Diophantine approximation on manifolds, where the following two conjectures are now of extremely high interest.

\begin{conjecture}\label{conj2}
Any analytic non-degenerate submanifold of\/ $\R^n$ is of Khintchine
type for convergence.
\end{conjecture}

\begin{conjecture}\label{conj3}
Let $\cM$ be a non-degenerate analytic submanifold of\/ $\R^n$,
$\ddd=\dim\cM$ and $\ccc=\codim\cM$. Let $\psi:\N\to\Rp$ be a
monotonic function. If\/ $\tfrac{\ccc}{\ccc+1}{\ddd}<s<\ddd$ then
\begin{equation}\label{e:150}
 \cH^s(\cS_n(\psi)\cap\cM)=0 \qquad \text{if} \qquad \sum_{q=1}^\infty \
 q^n\Big(\frac{\psi(q)}q\Big)^{s+\ccc }<\infty.
\end{equation}
\end{conjecture}

In the case of $\cM=\R^n$ Conjecture~\ref{conj3} together with
Theorem~\ref{t:03}\, coincides with \Jarnik's theorem, or rather the
modern version of \Jarnik's theorem -- see
\cite{Beresnevich-Dickinson-Velani-06:MR2184760}. Therefore,
Conjecture~\ref{conj3} can be regarded as a \Jarnik-type theorem for
convergence for manifolds. In turn, Theorem~\ref{t:03}\, can be
regarded as a \Jarnik-type theorem for divergence for manifolds. Conjecture~\ref{conj3} combined with Theorem~\ref{t:03} would also imply that (\ref{e:013}) is an equality.

If Theorem~\ref{t:04} was used to its `full potential' then one would be able to prove (\ref{e:012}) for $s\in(d/2,d)$. This naturally suggests the following problem: \emph{Describe analytic non-degenerate manifolds $\cM$ for which $(\ref{e:012})$ and/or
$(\ref{e:150})$ hold for $s\in(d/2,d)$.}

Note that within this paper there are two instances when (\ref{e:012}) is established for $s\in(d/2,d)$: hypersurfaces and curves. It is quite possible that for these types of manifolds (\ref{e:150}) also holds for $s\in(d/2,d)$.
However, note that for manifolds other than curves and hypersurfaces establishing (\ref{e:150}) for $s\in(d/2,d)$ is in general impossible unless extra constraints are added. This can be shown by considering $\cM$ as in Example~\ref{examp1}.

The main results of this paper are established in the case of analytic manifolds. However, within this paper the analyticity assumption is only used in establishing Theorem~\ref{t:06}. More precisely, the analyticity assumption is used to verify condition $(i)$ of Theorem~KM. A natural challenging question is then \emph{to what extent the analyticity assumption can be relaxed within Theorem~\ref{t:06} and consequently within all the main results of this paper}.

Recall that the analyticity assumption is not present in the planar curves results \cite{Beresnevich-Dickinson-Velani-07:MR2373145, Vaughan-Velani-2007}. Even though, there is a minor disagreement in the smoothness conditions imposed in convergence and divergence results: the divergence results deal with $C^{(3)}$ non-degenerate planar curve only. In general, the non-degeneracy of planar curves requires $C^{(2)}$. This raises the following intriguing question: \emph{Are $C^{(2)}$ non-degenerate planar curves of Khintchine type for
divergence?}

\bigskip

\noindent{\it Acknowledgements.} I would like to thank Vasili~Bernik and Maurice~Dodson for their valuable comments on the preliminary version of this paper.

\bigskip

{\footnotesize

\begin{minipage}{0.9\textwidth}
\footnotesize{\sc University of York, Heslington, York, YO10 5DD,
England}\\[0.5ex]
{\it E-mail address}\,:~~ \verb|vb8@york.ac.uk|\\[0.5ex]
\end{minipage}



\begin{thebibliography}{10}

\bibitem{Badziahin-Levesley-07:MR2347267}
{\sc D.~Badziahin and J.~Levesley}, {\em A note on simultaneous and
  multiplicative {D}iophantine approximation on planar curves}, Glasg. Math.
  J., 49 (2007), pp.~367--375.

\bibitem{Baker-1976}
{\sc R.~C. Baker}, {\em Metric {Diophantine} approximation on manifolds}, J.
  Lond. Math. Soc., 14 (1976), pp.~43--48.

\bibitem{Beresnevich-02:MR1905790}
{\sc V.~Beresnevich}, {\em A {G}roshev type theorem for convergence on
  manifolds}, Acta Math. Hungar., 94 (2002), pp.~99--130.

\bibitem{Beresnevich-Bernik-96:MR1387861}
{\sc V.~Beresnevich and V.~I. Bernik}, {\em On a metrical theorem of {W}.
  {S}chmidt}, Acta Arith., 75 (1996), pp.~219--233.

\bibitem{Beresnevich-Bernik-Dodson-02:MR2069553}
{\sc V.~Beresnevich, V.~I. Bernik, and M.~M. Dodson}, {\em On the {H}ausdorff
  dimension of sets of well-approximable points on nondegenerate curves}, Dokl.
  Nats. Akad. Nauk Belarusi, 46 (2002), pp.~18--20, 124.

\bibitem{Beresnevich-Dickinson-Velani-06:MR2184760}
{\sc V.~Beresnevich, D.~Dickinson, and S.~Velani}, {\em Measure theoretic laws
  for lim sup sets}, Mem. Amer. Math. Soc., 179 (2006), pp.~x+91.

\bibitem{Beresnevich-Dickinson-Velani-07:MR2373145}
\leavevmode\vrule height 2pt depth -1.6pt width 23pt, {\em Diophantine
  approximation on planar curves and the distribution of rational points}, Ann.
  of Math. (2), 166 (2007), pp.~367--426.
\newblock With an Appendix II by R. C. Vaughan.

\bibitem{Beresnevich-Vaughan-Velani-08-Inhom}
{\sc V.~Beresnevich, R.~Vaughan, and S.~Velani}, {\em Inhomogeneous
  {D}iophantine approximation on planar curves}.
\newblock Preprint, \url{http://arxiv.org/abs/0903.2817}, 2009.

\bibitem{Beresnevich-Velani-07:MR2285737}
{\sc V.~Beresnevich and S.~Velani}, {\em A note on simultaneous {D}iophantine
  approximation on planar curves}, Math. Ann., 337 (2007), pp.~769--796.

\bibitem{Beresnevich-Velani-09}
\leavevmode\vrule height 2pt depth -1.6pt width 23pt, {\em Ubiquity and a
  general logarithm law for geodesics}, S\'eminaires \& Congr\`es, 22 (2009),
  pp.~21--36.

\bibitem{Beresnevich-Bernik-Kleinbock-Margulis-02:MR1944505}
{\sc V.~V. Beresnevich, V.~I. Bernik, D.~Y. Kleinbock, and G.~A. Margulis},
  {\em Metric {D}iophantine approximation: the {K}hintchine-{G}roshev theorem
  for nondegenerate manifolds}, Mosc. Math. J., 2 (2002), pp.~203--225.
\newblock Dedicated to Yuri I. Manin on the occasion of his 65th birthday.

\bibitem{Bernik-73:MR0337863}
{\sc V.~I. Bernik}, {\em Asymptotic behavior of the number of solutions of
  certain systems of inequalities in the theory of {D}iophantine approximations
  of dependent variables}, Vesc\=\i\ Akad. Navuk BSSR Ser. F\=\i z.-Mat. Navuk,
   (1973), pp.~10--17, 134.

\bibitem{Bernik-77:MR0480402}
\leavevmode\vrule height 2pt depth -1.6pt width 23pt, {\em An analogue of
  {H}in\v cin's theorem in the metric theory of {D}iophantine approximations of
  dependent variables. {I}}, Vesc\=\i\ Akad. Navuk BSSR Ser. F\=\i z.-Mat.
  Navuk,  (1977), pp.~44--49, 141.

\bibitem{Bernik-1979}
\leavevmode\vrule height 2pt depth -1.6pt width 23pt, {\em On the exact order
  of approximation of almost all points on the parabola}, Mat. Zametki, 26
  (1979), pp.~657--665.
\newblock (In Russian).

\bibitem{Bernik-1983a}
\leavevmode\vrule height 2pt depth -1.6pt width 23pt, {\em An application of
  {Hausdorff} dimension in the theory of {Diophantine} approximation}, Acta
  Arith., 42 (1983), pp.~219--253.
\newblock (In Russian). English transl. in {\it Amer. Math. Soc. Transl.}
  {\bf140} (1988), 15--44.

\bibitem{Bernik-03:MR2163817}
\leavevmode\vrule height 2pt depth -1.6pt width 23pt, {\em The {K}hinchin
  transference principle and lower bounds for the number of rational points
  near smooth manifolds}, Dokl. Nats. Akad. Nauk Belarusi, 47 (2003),
  pp.~26--28.

\bibitem{BernikDodson-1999}
{\sc V.~I. Bernik and M.~M. Dodson}, {\em Metric {Diophantine} approximation on
  manifolds}, vol.~137 of Cambridge Tracts in Mathematics, Cambridge University
  Press, Cambridge, 1999.

\bibitem{Bernik-Kleinbock-Margulis-01:MR1829381}
{\sc V.~I. Bernik, D.~Kleinbock, and G.~A. Margulis}, {\em Khintchine-type
  theorems on manifolds: the convergence case for standard and multiplicative
  versions}, Internat. Math. Res. Notices,  (2001), pp.~453--486.

\bibitem{Bocher-00:MR1503482}
{\sc M.~B{\^o}cher}, {\em The theory of linear dependence}, Ann. of Math. (2),
  2 (1900/01), pp.~81--96.

\bibitem{Budarina-Dickinson-07:MR2397137}
{\sc N.~Budarina and D.~Dickinson}, {\em Simultaneous {D}iophantine
  approximation on surfaces defined by polynomial expressions {$x\sp d\sb
  1+\dots+ x\sp d\sb m$}}, in Analytic and probabilistic methods in number
  theory/{A}naliziniai ir tikimybiniai metodai skai\v ci\polhk u teorijoje,
  TEV, Vilnius, 2007, pp.~17--23.

\bibitem{DickinsonDodson-2000a}
{\sc H.~Dickinson and M.~M. Dodson}, {\em Extremal manifolds and {Hausdorff}
  dimension}, Duke Math. J., 101 (2000), pp.~271--281.

\bibitem{DickinsonDodson-2001a}
\leavevmode\vrule height 2pt depth -1.6pt width 23pt, {\em Diophantine
  approximation and {Hausdorff} dimension on the circle}, Math. Proc. Cambridge
  Phil. Soc., 130 (2001), pp.~515--522.

\bibitem{DodsonRynneVickers-1989b}
{\sc M.~M. Dodson, B.~P. Rynne, and J.~Vickers}, {\em Metric {Diophantine}
  approximation and {Hausdorff} dimension on manifolds}, Math. Proc. Cam. Phil.
  Soc., 105 (1989), pp.~547--558.

\bibitem{DodsonRynneVickers-1991a}
\leavevmode\vrule height 2pt depth -1.6pt width 23pt, {\em Kchintchine-type
  theorems on manifolds}, Acta Arithmetica, 57 (1991), pp.~115--130.

\bibitem{Dodson-Rynne-Vickers-1996}
\leavevmode\vrule height 2pt depth -1.6pt width 23pt, {\em Simultaneous
  {Diophantine} approximation and asymptotic formulae on manifolds}, J. Number
  Theory, 58 (1996), pp.~298--316.

\bibitem{Drutu-05:MR2195121}
{\sc C.~Dru{\c{t}}u}, {\em Diophantine approximation on rational quadrics},
  Math. Ann., 333 (2005), pp.~405--469.

\bibitem{Elkies-2000}
{\sc N.~Elkies}, {\em Rational points near curves and small nonzero {$\vert
  x\sp 3-y\sp 2\vert $} via lattice reduction}, in Algorithmic number theory
  (Leiden, 2000), vol.~1838 of Lecture Notes in Comput. Sci., Springer, Berlin,
  2000, pp.~33--63.

\bibitem{Gorodnik-Shah-08}
{\sc A.~Gorodnik and N.~Shah}, {\em Khinchin theorem for integral points on
  quadratic varieties}, Preprint. arXiv:0804.3530.

\bibitem{Harman-03:MR1979906}
{\sc G.~Harman}, {\em Simultaneous {D}iophantine approximation and asymptotic
  formulae on manifolds}, Acta Arith., 108 (2003), pp.~379--389.

\bibitem{Huxley-1994-rational_points}
{\sc M.~N. Huxley}, {\em The rational points close to a curve}, Ann. Scuola
  Norm. Sup. Pisa Cl. Sci. (4), 21 (1994), pp.~357--375.

\bibitem{Huxley-96:MR1420620}
\leavevmode\vrule height 2pt depth -1.6pt width 23pt, {\em Area, lattice
  points, and exponential sums}, vol.~13 of London Mathematical Society
  Monographs. New Series, The Clarendon Press Oxford University Press, New
  York, 1996.
\newblock , Oxford Science Publications.

\bibitem{Jarnik-1931}
{\sc A.~Jarnik}, {\em Uber die simultanen diophantischen {Approximationen},},
  Math. Z., 33 (1931), pp.~505--543.

\bibitem{Khintchine-1924}
{\sc A.~Khintchine}, {\em Einige {S}\"atze \"uber {K}ettenbr\"uche, mit
  {Anwendungen} auf die {Theorie} der {Diophantischen} {Approximationen}},
  Math. Ann., 92 (1924), pp.~115--125.

\bibitem{Khintchine-1925}
\leavevmode\vrule height 2pt depth -1.6pt width 23pt, {\em Zwei {Bemerkungen}
  zu einer {Arbeit} des {Herrn} {Perron}}, Math. Zeitschr., 22 (1925),
  pp.~274--284.

\bibitem{Khintchine-1926}
\leavevmode\vrule height 2pt depth -1.6pt width 23pt, {\em Zur metrischen
  {Theorie} der diophantischen {Approximationen}}, Math. Zeitschr., 24 (1926),
  pp.~706--714.

\bibitem{Kleinbock-03:MR1982150}
{\sc D.~Kleinbock}, {\em Extremal subspaces and their submanifolds}, Geom.
  Funct. Anal., 13 (2003), pp.~437--466.

\bibitem{Kleinbock-04:MR2094125}
\leavevmode\vrule height 2pt depth -1.6pt width 23pt, {\em Baker-{S}prind\v zuk
  conjectures for complex analytic manifolds}, in Algebraic groups and
  arithmetic, Tata Inst. Fund. Res., Mumbai, 2004, pp.~539--553.

\bibitem{Kleinbock-Lindenstrauss-Weiss-04:MR2134453}
{\sc D.~Kleinbock, E.~Lindenstrauss, and B.~Weiss}, {\em On fractal measures
  and {D}iophantine approximation}, Selecta Math. (N.S.), 10 (2004),
  pp.~479--523.

\bibitem{Kleinbock-Tomanov-07:MR2314053}
{\sc D.~Kleinbock and G.~Tomanov}, {\em Flows on {$S$}-arithmetic homogeneous
  spaces and applications to metric {D}iophantine approximation}, Comment.
  Math. Helv., 82 (2007), pp.~519--581.

\bibitem{Kleinbock-Margulis-98:MR1652916}
{\sc D.~Y. Kleinbock and G.~A. Margulis}, {\em Flows on homogeneous spaces and
  {D}iophantine approximation on manifolds}, Ann. of Math. (2), 148 (1998),
  pp.~339--360.

\bibitem{Mahler-1932b}
{\sc K.~Mahler}, {\em \"uber das {M}a\ss der {M}enge aller ${S}$-{Z}ahlen},
  Math. Ann., 106 (1932), pp.~131--139.

\bibitem{Mashanov-87:MR888593}
{\sc V.~I. Mashanov}, {\em On a problem of {B}aker in the metric theory of
  {D}iophantine approximations}, Vests\=\i\ Akad. Navuk BSSR Ser. F\=\i z.-Mat.
  Navuk,  (1987), pp.~34--38, 125.

\bibitem{Mazur-04:MR2058289}
{\sc B.~Mazur}, {\em Perturbations, deformations, and variations (and
  ``near-misses'') in geometry, physics, and number theory}, Bull. Amer. Math.
  Soc. (N.S.), 41 (2004), pp.~307--336 (electronic).

\bibitem{Melnichuk-1979a}
{\sc Y.~Melnichuk}, {\em Diophantine approximation on a circle and {Hausdorff}
  dimension}, Mat. Zametki, 26 (1979), pp.~347--354.
\newblock English transl. in Math. notes {\bf26} (1980), 666--670.

\bibitem{Pyartli-1969}
{\sc A.~Pyartli}, {\em Diophantine approximation on submanifolds of euclidean
  space}, Funkts. Anal. Prilosz., 3 (1969), pp.~59--62.
\newblock (In Russian).

\bibitem{Schmidt-64:MR0171753}
{\sc W.~M. Schmidt}, {\em Metrische {S}\"atze \"uber simultane {A}pproximation
  abh\"angiger {G}r\"ossen}, Monatsh. Math., 68 (1964), pp.~154--166.

\bibitem{Schmidt-1980}
\leavevmode\vrule height 2pt depth -1.6pt width 23pt, {\em Diophantine
  {Approximation}}, Springer-Verlag, Berlin and New York, 1980.

\bibitem{Sprindzuk-1969-Mahler-problem}
{\sc V.~Sprind\v{z}uk}, {\em Mahler's problem in the metric theory of numbers},
  vol.~25, Amer. Math. Soc., Providence, RI, 1969.
\newblock Translations of Mathematical Monographs.

\bibitem{Sprindzuk-1980-Achievements}
\leavevmode\vrule height 2pt depth -1.6pt width 23pt, {\em Achievements and
  problems in {Diophantine} approximation theory}, Russian Math. Surveys, 35
  (1980), pp.~1--80.

\bibitem{Vaughan-Velani-2007}
{\sc R.~C. Vaughan and S.~Velani}, {\em Diophantine approximation on planar
  curves: the convergence theory}, Invent. Math., 166 (2006), pp.~103--124.

\bibitem{Weyl-Weyl-38:MR1503422}
{\sc H.~Weyl and J.~Weyl}, {\em Meromorphic curves}, Ann. of Math. (2), 39
  (1938), pp.~516--538.

\bibitem{Whitney}
{\sc H.~Whitney}, {\em Geometric integration theory}, Princeton University
  Press, Princeton, N. J., 1957.

\end{thebibliography}

\def\cprime{$'$} \def\cprime{$'$} \def\cprime{$'$} \def\cprime{$'$}
  \def\cprime{$'$} \def\cprime{$'$} \def\cprime{$'$} \def\cprime{$'$}
  \def\cprime{$'$} \def\cprime{$'$} \def\cprime{$'$} \def\cprime{$'$}
  \def\cprime{$'$} \def\cprime{$'$} \def\cprime{$'$}
  \def\polhk#1{\setbox0=\hbox{#1}{\ooalign{\hidewidth
  \lower1.5ex\hbox{`}\hidewidth\crcr\unhbox0}}} \def\cprime{$'$}
  \def\cprime{$'$}

}

\end{document}